\documentclass[reqno]{amsart}
\usepackage{amssymb,amscd,verbatim, amsthm,graphicx, color}
\usepackage{fancybox}
\usepackage{longtable}

\newcommand \fk[1]{{{\mathfrak #1}}}
\newcommand \C[1]{{\mathcal #1}}
\newcommand \ovl[1]{{\overline {#1}}}
\newcommand \ch[1]{{\check{#1}}}
\newcommand \bb[1]{{\mathbb #1}}

\newcommand \wti[1]{{\widetilde {#1}}}
\newcommand \wht[1]{{\widehat {#1}}}

\newcommand \unb[2]{\underset{#1}{{\underbrace{#2}}}}

\newcommand\fg{\mathfrak g}
\newcommand\fn{\mathfrak n}
\newcommand\fp{\mathfrak p}
\newcommand\fm{\mathfrak m}

\newcommand \bA{{\mathbb A}}
\newcommand \bC{{\mathbb C}}
\newcommand \bF{{\mathbb F}}
\newcommand \bH{{\mathbb H}}
\newcommand \bR{{\mathbb R}}
\newcommand \bZ{{\mathbb Z}}
\newcommand \bG{{\mathbb G}}
\newcommand \bB{{\mathbb B}}
\newcommand \bO{{\mathbb O}}

\newcommand\one{1\!\!1}

\newcommand\CA{{\C A}}

\newcommand\CH{{\C H}}
\newcommand\CI{{\C I}}
\newcommand\CO{{\C O}}

\newcommand\CP{{\C P}}

\newcommand\ic{infinitesimal character }

\newcommand\ie{{\it i.e.~}}
\newcommand\cf{{\it cf.~ }}
\newcommand\eg{{\it e.g.~ }}

\newcommand\ep{{\epsilon}}

\newcommand\om{{\omega}}
\newcommand\al{{\alpha}}

\newcommand \vH{{\check H}}
\newcommand \vR{{\check R}}

\newcommand \vG{{^\vee\! G}}

\newcommand\fa{{\mathfrak a}}
\newcommand\fh{{\mathfrak h}}

\newcommand\fz{{\mathfrak z}}

\newtheorem*{theorem}{Theorem}
\newtheorem*{theorem 2}{Theorem 2}
\newtheorem*{corollary}{Corollary}
\newtheorem*{lemma}{Lemma}
\newtheorem*{proposition}{Proposition}
\newtheorem*{definition}{Definition}
\newtheorem*{remark}{Remark}
\newtheorem*{example}{Example}

\newcommand\Hom{\operatorname{Hom}}
\newcommand\Ind{\operatorname{Ind}}
\newcommand\rank{\operatorname{rank}}

\numberwithin{equation}{subsection}

\begin{document}

\today

\bigskip
\title[Unitarizable minimal principal series]{Unitarizable minimal
  principal series of reductive groups}
\author{Dan Barbasch}
       \address[D. Barbasch]{Dept. of Mathematics\\
               Cornell University\\Ithaca, NY 14850}
       \email{barbasch@math.cornell.edu}

\author{Dan Ciubotaru}
        \address[D. Ciubotaru]{Dept. of Mathematics\\ University of
          Utah\\ Salt Lake City, UT 84112}
        \email{ciubo@math.utah.edu}

\author{Alessandra Pantano}
      \address[A. Pantano]{Dept. of Mathematics\\University of
        California\\ Irvine, CA}
        \email{apantano@math.uci.edu}

%{\it Dedicated to B. Casselman and D. Mili\v ci\'c.}

\begin{abstract}
The aim of this paper is to give an exposition of recent progress on
the determination of the unitarizable Langlands quotients of
minimal principal series for reductive groups over the real or
$p$-adic fields in characteristic $0$.
\end{abstract}

\maketitle

{\it Dedicated to B. Casselman and D. Mili\v ci\'c.}

\setcounter{tocdepth}{1}

\begin{small}
\tableofcontents
\end{small}

\section{Introduction}\label{sec:1}

\subsection{}\label{sec:1.1}The aim of this paper is to give an exposition of some  recent
progress on the computation of the unitary dual of a reductive group
over a local field of characteristic $0$. The results are
for Langlands quotients of minimal principal series, and form the
subject of \cite{Ba1}, \cite{BP}, and \cite{BC}.

Let $\bF$ be the real field or a $p$-adic field of characteristic
zero, and let $|\ |$ denote the absolute value, respectively the
$p$-adic norm. When $\bF$ is $p$-adic, we let
\begin{equation}
\bO=\{x\in\bF: |x|\le 1\},\text{ and } \CP=\{x\in\bF: |x|<1\},
\end{equation}
be the ring of integers, respectively the unique prime ideal. Let
$\bO^\times=\bO\setminus \CP$ be the set of invertible elements in
the ring of integers. We fix an uniformizer $\varpi\in \C P,$
and then $\CP=\varpi\bO.$ The quotient $\bO/\CP$ is isomorphic to
a finite field  $\bF_q$ of characteristic $p.$

Let $G(\bF)$ be the $\bF$-points of a linear connected reductive
group $G$ defined over $\bF$.  Assume that $G(\bF)$ is split.
This means that $G$ has a Cartan subgroup $H$ such
that $H(\bF)$ is isomorphic to a product of $r:=\rank(G)$ copies of
$\bF^\times.$  In this paper we will deal almost exclusively with split groups.
Since $G(\bF)$ is split, $G$ also has a Borel subgroup which is defined over
$\bF.$ Choose  such a Borel subgroup $B$,  and a split Cartan subgroup
$H\subset B.$

The reductive connected algebraic group $G$ is
determined by the root datum $(\C X,R,\C Y,\ch R),$ where $\C X$ is the
lattice of algebraic characters of $H,$ $\C Y$ is the lattice of
algebraic cocharacters of $H,$ $R$ is the set of roots of
$H$ in $G$ and $\check R$ are the coroots.  Then
$H(\bF)=\C Y\otimes_\bZ \bF^\times.$ The choice of Borel $B$
determines the positive roots $R^+$ and positive coroots $\check R^+.$
Also let $\fh=\C Y\otimes_\bZ\bC$ be the complex Lie algebra of
$H(\bF).$

Let $\vG$ be the
\textit{dual} complex connected group corresponding to $(\C Y,\ch R,\C
X,R).$ Denote by $\vH:=\C X\otimes_\bZ\bC^\times$.  This is a Cartan
subgroup of $\vG.$ Its Lie algebra is $\check\fh=\C X\otimes_\bZ\bC.$
The torus $\vH$ has a polar decomposition
\begin{equation}
\vH=\vH_e\cdot \vH_h,\text{ where } \vH_e=\C X\otimes_\bZ S^1 \text{ and } \vH_h=\C X\otimes_\bZ\bR_{>0}.
\end{equation}
A semisimple element $s\in\vG$ is called {elliptic} (respectively,
hyperbolic) if $s$ is conjugate to an element of $\vH_e$
(respectively, $\vH_h$). Similarly, the Lie algebra $\check\fh$ of $\vH$ has a decomposition into a real (hyperbolic) part
$\check\fh_\bR=\C X\otimes_\bZ\bR$ and an imaginary (elliptic) part
$\check\fh_{i\bR}=\C X\otimes_\bZ i\bR.$

\smallskip

We fix a maximal compact subgroup $K$ in $G(\bF)$ in the real case,
and $K=G(\bb O)$ in the p-adic case.
Set  $^0\!H=H(\bF)\cap K.$ In
the real case $^0\!H$ is isomorphic to $(\rank G)$-copies of
$\bZ/2\bZ.$ The familiar notation for $^0 H$ in this case is $M.$  In the $p$-adic case, $^0\!H$ is isomorphic with $(\rank
G)$-copies of $\bO^\times$.

For a character $\chi$ of $H(\bF)$, the {\it minimal
principal series} $X(\chi)$ is defined to be
\begin{equation}
X(\chi)=\Ind_{B(\bF)}^{G(\bF)}(\chi\otimes 1),
\end{equation}
where $\Ind$ means unitary induction.
 Each such principal series
has finite composition series, and in particular a canonical
completely reducible subquotient denoted $L(\chi).$ This paper
is concerned with the problem of when the constituents of $L(\chi)$
are unitary.

A character $\chi \colon H(\bF)\to \bC$ is called {\it spherical} (or
{\it unramified}) if its restriction to $^0\!H$ is trivial, \ie
$^0\!\chi:=\chi|_{^0\!H}=triv.$ It is called {\it nonspherical}
(or {\it ramified}) otherwise.

There are two main themes in this paper: the first is to describe what
is known about the unitarizability of the unramified subquotient
$L(\chi)$, and the second is to relate the unitarizability of the
Langlands subquotients of ramified minimal principal series $X(\chi)$ of
$G(\bF)$ with the
spherical  unitary dual for certain
groups $G(^0\!\chi)$ attached to $^0\!\chi$ (see definition in \ref{sec:1.8}).

\subsection{}\label{sec:1.2} We begin with the case when the character
$\chi$ of $H(\bF)$ is unramified.

A representation $(\pi,V)$ is called {\it spherical} if the set $V^K$
of vectors fixed by $K$ is nontrivial. By Frobenius reciprocity
$X(\chi)^K$ is 1-dimensional when $\chi$ is unramified. Thus $X(\chi)$
is spherical, and has a unique irreducible subquotient which is
spherical. In the unramified case, this is precisely the canonical
subquotient $L(\chi)$ alluded to earlier.

The basic example, going back to \cite{Bar} is  $G=SL(2,\bR).$ Then
$H(\bR)=\bR^\times.$ The unramified character $\chi$ can be written as
\begin{align}\label{eq:1.2.1}
\chi=\chi_\nu:H(\bR)=\bR^\times\to \bC,\ \chi_\nu(z)=|z|^\nu,\ z\in
\bR^\times, \text{ for some }\nu\in\bC.
\end{align}
When $\nu$ is purely imaginary, $X(\chi_\nu)$ is
irreducible (and unitary). If $\nu$ is real, then $X(\nu)$ is
reducible if and only if $\langle\check\al,\nu\rangle\in 2\bZ+1.$ When
$\langle\check\al,\nu\rangle=1,$ respectively $-1$, then the trivial
representation is a quotient, respectively a submodule of $X(\chi).$
The spherical Langlands
subquotient $L(\chi_\nu)$  is unitary  for
\begin{align}\label{eq:1.2.2}
&\nu\in i\bR,\quad
\text{
  and},\\\notag
&\nu\in\bR, \text{ such that } -1\le\langle\check\al,\nu\rangle\le 1.
\end{align}

\subsection{}\label{sec:1.3} Now we consider the similar example
for $SL(2,\bF),$ when $\bF$ is $p$-adic. This is again well-known
(\cite{Sa}), and we refer the reader to chapter 9 of \cite{Cas2} for a
detailed treatment. The unramified characters of
$\bH(\bF)=\bF^\times$ can be described similarly to (\ref{eq:1.2.1}),
but since now $|\ |$ is the $p$-adic norm, the actual
parameterization is
\begin{align}\label{eq:1.3.1}
\chi=\chi_\nu:H(\bF)=\bF^\times\to \bC,\ \chi_\nu(z)=|z|^\nu,\ z\in
\bF^\times, \text{ for some }\nu\in\bC/(2\pi i/\log q)\bZ.
\end{align}
When $\nu$ is real, $X(\chi_\nu)$ is reducible if and only if
$\langle\check\al,\nu\rangle=\pm 1,$ which is the value of the
parameter where the trivial representation is a subquotient. But when
$\nu$ is purely imaginary, there are two cases one needs to
consider. If $\nu\neq \pi i/\log q,$ then $X(\chi_\nu)$ is irreducible
(and unitary). If $\nu=\pi i/\log q$, then $X(\chi_\nu)$ decomposes
into a sum of two irreducible submodules, both of which are
unitary. It is worth noting that the Weyl group $W=S_2$ acts by
$\nu\to -\nu$,  and $\nu=\pi i/\log q$ is the only
  nontrivial parameter so that $\chi_\nu$ is fixed by $W.$

The spherical Langlands subquotient $L(\chi_\nu)$ is unitary for
\begin{align}\label{eq:1.3.2}
&\nu\in i(\bR/(2\pi/\log q)\bZ),\ \nu\neq \pi i/\log q,\\\notag
&\nu=\pi i/\log q,\quad {and}\\\notag
& \nu\in\bR, \text{ such that } -1\le\langle\check\al,\nu\rangle\le 1.
\end{align}

What should be noted from equations (\ref{eq:1.2.1}), (\ref{eq:1.2.2}) and
(\ref{eq:1.3.1}), (\ref{eq:1.3.2}) is that, if $\nu$ is assumed real, then the
parameter sets for the spherical representations {\it and } the
unitary sets are {\it the same} for both real and $p$-adic case. This
suggests the following natural questions:
\begin{enumerate}
\item[(a)] Can the determination of the unitary representations be reduced
  to the case when the parameter is (in some technical sense) real?
\item[(b)] If the parameter is assumed real, can one generalize the
  identification of spherical unitary duals between the real and
  $p$-adic cases to all split groups?
\end{enumerate}

When $\bF=\bR,$ the answer to (a) is the well-known reduction to
real infinitesimal character to Levi subgroups for unitary representations (see
Chapter XVI in \cite{Kn} for example). In the particular case of
spherical representations, the
reduction has a simple form. Assume
$\nu=\Re\nu+{{i}}\Im\nu$ with
$\Im\nu\neq 0.$  $\Im\nu$  defines a proper parabolic subgroup
$P_{\Im\nu}=L_{\Im\nu} U_{\Im\nu}$ of $G$ (similarly to the
construction in section \ref{sec:2.9}).  Then
\begin{equation}\label{eq:1.3.3}
L(\chi_\nu)=\Ind_{P_{\Im\nu}}^G(L_{M_{\Im\nu}}(\chi_{\Re\nu})),
\end{equation}
and since this is unitary induction, $L(\chi_\nu)$ is unitary for $G$ if and
only if $L_{M_{\Im\nu}}(\chi_{\Re\nu})$ is unitary for $M_{\Im\nu}$ (a
proper Levi).

When $\bF$ is $p$-adic, the situation is more complicated as already
seen in the case $\nu=\pi i/\log q$ in $SL(2)$. The answer is positive
nevertheless for the representations which appear in unramified
principal series, and it is the subject of \cite{BM2}.

The answer to (b) is known to be positive, at least when $G$ is of
classical type (\cite{Ba1}). But in order to explain this (and the
reduction to ``real infinitesimal character'' in the $p$-adic case),
we need to introduce the Iwahori-Hecke algebra and its graded
versions.

\subsection{}\label{sec:1.4}
Let $\bF$ be a $p$-adic field.
In this case, the set of spherical irreducible representations
are part of a larger class, the Iwahori-spherical representations.
Let
$K_1\subset K$ be the subgroup of $K$ of elements congruent to  $Id$
modulo $\CP.$ Then there is an exact sequence
\begin{equation}\label{eq:1.4.1}
  1\longrightarrow K_1\longrightarrow K\overset{\pi}{\longrightarrow}
K/K_1\cong  G(\bF_q)\longrightarrow 1.
\end{equation}
The group $B(\bF_q)$ is a Borel subgroup of $G(\bF_q).$
Then $\C I:=\pi^{-1}(B(\bF_q))$ is an open compact subgroup of $K,$
called an {\it Iwahori subgroup}. The {\it Iwahori-Hecke algebra},
denoted $\C H$ is the algebra of locally constant compactly supported $\C
I$-biinvariant functions under convolution.  If $(\pi,V)$ has $\C
I$-fixed vectors, then $\C H$ acts
on $V^\C I.$ Let $C(\C I,triv)$ be  the category of admissible
representations such that each subquotient is generated by its  $\C
I$-fixed
vectors. A central theorem in the representation theory of reductive
$p$-adic groups is the following result of Borel and Casselman.
\begin{theorem}[Borel-Casselman]\label{t:bc}
The association $V\mapsto V^\C I$ is an equivalence of categories
between $C(\C I,triv)$ and the category of finite dimensional
representations of $\C H.$
\end{theorem}

\subsection{}\label{sec:1.5}The Hecke algebra $\C H$ also has a $*$ operation
$f^*(g):=\ovl{f(g^{-1})}.$ So it makes sense to talk about hermitian
and unitary $\C H$-modules. The next theorem transforms the analytic
problem of classifying $\C I$-spherical unitary representation of the group
$G(\bF)$ ($\bF$ $p$-adic), which is about infinite dimensional
representations, to the algebraic problem of classifying finite
dimensional unitary representations of $\C H.$
\begin{theorem}[\cite{BM1},\cite{BM2}]\label{t:bm}
$(\pi, V)$ is a unitary representation of $G(\bF)$, $\bF$ $p$-adic, if
  and only if $V^{\C I}$ is a unitary module of the Iwahori-Hecke
  algebra $\C H.$
\end{theorem}
 The proof of this theorem requires the classification of
  the irreducible finite dimensional representations of $\C H$ from
  the work of Kazhdan-Lusztig (\cite{KL}). We give an outline of the machinery
  involved, more details are available in later sections.
The algebra $\C H$ has two well-known descriptions in terms of
generators and relations. The original description is in terms of the
affine Weyl group (\cite{IM}). We will use a second one due to Bernstein,
because it is better suited for our purposes (see section \ref{sec:3.0}).
Schur's lemma for $\C H$ holds, so the center of $\C H$ acts by scalars on any
irreducible finite dimensional representation. By the aforementioned
result of Bernstein, the characters of the
center of $\C H$, which we call
{\it infinitesimal characters} in analogy with the real case, can be
identified with semisimple orbits of $\vG.$

Therefore, they are also in 1-1 correspondence with $W$-orbits of
elements of $\vH.$ Suppose $s\in\vH$ denotes an infinitesimal
character. Decompose $s=s_e\cdot s_h,$ where $s_e$ is elliptic, and
$s_h$ hyperbolic.

\begin{definition}\label{d:1.5}The infinitesimal character $s$ is
called \textit{real}, if the $W$-orbit of $s_e$ has only one element.
\end{definition}
\cite{BM1} proves theorem \ref{t:bm} in the
case when the infinitesimal character is real. Then \cite{BM2} extends
the result to all infinitesimal characters. In the process, \cite{BM2}
shows that it is sufficient to determine the unitary dual of the
graded version of the Iwahori-Hecke algebra defined in \cite{L1},
denoted $\bH,$  and in addition it is enough to
consider modules with real infinitesimal character only.

The affine graded Hecke
algebra is the associated graded object to a filtration of ideals of
$\C H$ defined relative to a ($W$-orbit) of an elliptic element $s_e.$
It has an explicit description in terms of generators and relations,
similar to the Bernstein presentation (\cite{L1}).
For example when $s_e=1$, as a vector space, $\bH$ is generated by a
copy of $\bC[W],$ the group algebra of $W$,
and a copy of the symmetric algebra $\bA=Sym(\fh)$, with a nontrivial
commutation relation between $\bC[W]$ and $\bA$. (See \ref{sec:3.1}
for the precise definition.) The center of $\bH$ is $\bA^W,$ so the
central (``infinitesimal'') characters are parameterized by $W$-orbits
in $\check\fh.$ We call an infinitesimal character for $\bH$ {\it
  real} if it is an orbit of elements in $\check\fh_\bR.$

When $s_e\neq 1,$ we refer the reader to \cite{BM2}, and we will only
give the details (which are representative) for the case of
$SL(2,\bF).$
\cite{BM2} proves that the unitarity of a $\C H$-module with
infinitesimal character $s=s_es_h$ is equivalent with the unitarity of
the corresponding $\bH_{s_e}$-module with (real) infinitesimal character
$\log s_h$, where $\bH_{s_e}$ is the algebra obtain by ``grading'' at
$W\cdot s_e.$

\subsection{}\label{sec:1.6} We explain the example of
$SL(2,\bF),$ $\bF$ $p$-adic, in this setting, and illustrate the role
of the Hecke algebra.  The dual complex group is $\vG=PGL(2,\bC)$. The
algebra $\C H$ is generated (in Bernstein's presentation) by $T$
corresponding to the nontrivial element of $W$, and $\theta$
corresponding to the unique coroot,  subject to the relations
\begin{equation}
  \label{eq:1.6.1}
  T^2=(q-1)T +q,\qquad T\theta =\theta^{-1}T+(q-1)(\theta+1).
\end{equation}
The center of $\C H$ is generated by $\theta + \theta^{-1}.$
Homomorphisms of the center correspond to $W$-orbits of semisimple
elements $s=s_es_h\in\vH.$ In the correspondence \ref{t:bc}, if $V$ is
a subquotient of $X(\chi_\nu)$ (as in (\ref{eq:1.3.1})), then the
center of $\C H$ acts on   $V^{\C I}$ by
$q^\nu+q^{-\nu}.$ The infinitesimal character $s\in \vH$ is (the
$W$-orbit of)
$\left(\begin{matrix}a &0\\ 0&a^{-1}\end{matrix}\right),$ where
$a=\nu\log q/\pi.$
As noted in section \ref{sec:1.3}, the principal series $X(\chi_\nu)$
is irreducible for $\nu\neq 1,\nu\neq\pi i/\log q.$ When it is reducible,
the
constituents are: the trivial module ($triv$), the
Steinberg module ($St$), which is the unique submodule (discrete
series) at $\nu=1,$ the two components $\overline X_{sph}$ and
$\overline X_{nonsph}$ at $\nu=\pi i/\log q.$
They
correspond to the four one-dimensional Hecke
modules as follows:

\bigskip

\noindent\begin{tabular}{|c|c|c|c|c|}
\hline
$SL(2)$-mod &$\C H$-mod &Action of $Z(\C H)$ &Inf. char. $s$ &$s_e$\\
\hline
$triv$ &$T=q,\ \theta=q$ &$q+q^{-1}$ & $\left(\begin{matrix} \log
  q/\pi&0\\0&\pi/\log q\end{matrix}\right)$ &$\left(\begin{matrix} 1&0\\0&1\end{matrix}\right)$\\
$St$   &$T=-1, \ \theta=q^{-1}$ &$q+q^{-1}$ &&\\\hline
$\overline X_{sph}$ &$T=q,\ \theta=-1$ &$-2$ & $\left(\begin{matrix} i&0\\0&-i\end{matrix}\right)$& $\left(\begin{matrix} i&0\\0&-i\end{matrix}\right)$
\\
$\overline X_{nonsph}$ &$T=-1,\ \theta=-1$ &$-2$ &&\\\hline
\end{tabular}

\bigskip

Relative to $s_e$, there are three cases as follows.

\subsubsection{}\label{sec:1.6.1} $s_e=\left(\begin{matrix} 1&0\\0&1\end{matrix}\right)$.
The infinitesimal character $s$ is real. The affine graded Hecke algebra is generated by
$t,\ \om$ subject to the relations
\begin{equation}
  \label{eq:hck1}
  t^2=1,\qquad t\om +\om t=1.
\end{equation}
%The center is generated by $\om^2.$
This case is relevant for the unitarity of the complementary series
in $SL(2,\bF).$

\subsubsection{}\label{sec:1.6.2}\   $s_e=\left(\begin{matrix}
    i&0\\0&-i\end{matrix}\right)$. The infinitesimal character $s$ is real. Because $\theta(s_e)=-1,$ the affine graded Hecke algebra
is the group algebra of the affine group generated by $t,\ \om$ with
relations
\begin{equation}
  \label{eq:hck2}
t^2=1,\qquad  t\om +\om t=0.
\end{equation}
The finite dimensional representation theory of this algebra is also
well known, but quite different from case \ref{sec:1.6.1}. In particular, at real
infinitesimal character the only unitary representations occur at
$s_h=1.$ They correspond to the two modules $\overline X_{sph}$ and
$\overline X_{nonsph}.$

\subsubsection{} \label{sec:1.6.3} $s_e=\left(\begin{matrix}\zeta&0\\ 0&\zeta^{-1}\end{matrix}\right)$.
When $\zeta^2\ne \pm 1,$ the infinitesimal characters are not real, so there is a
reduction to a smaller algebra. The affine graded Hecke algebra $\bH_{s_e}$ in \cite{BM2} is generated
by $E_\zeta,\ E_{\zeta^{-1}},\ t,\ \om$ satisfying the following  relations:
  \begin{enumerate}
  \item $E_\zeta^2=E_\zeta,\ E_{\zeta^{-1}}^2=E_{\zeta^{-1}},\
    E_\zeta\cdot E_{\zeta^{-1}}=0,\ E_\zeta +E_{\zeta^{-1}}=1,$ in
    other words they are projections,
\item $tE_\zeta=E_{\zeta^{-1}}t,\qquad E_\zeta t=tE_{\zeta^{-1}},$
\item $t^2=1,\qquad t\om +\om t=1.$
  \end{enumerate}
Let $\C M_2$ be the algebra of $2\times 2$ matrices with complex
coefficients and the usual
basis $E_{ij}.$ Let $\bA$ be the polynomial algebra generated by
$\om.$ Then theorem 3.3 in \cite{BM2} states that the map $m\otimes
a\mapsto \Psi(m)\cdot a$ from $\C M_2\otimes_\bC\bA$ to $\bH_{s_e}$
\begin{equation}
  \label{eq:hck3}
  \begin{aligned}
&\Psi:\C M_2\longrightarrow \bH_{s_e}\\
& E_{11}\mapsto E_\zeta,\quad E_{12}\mapsto E_\zeta t,\quad
E_{21}\mapsto tE_\zeta,\quad E_{22}\mapsto tE_\zeta t=E_{\zeta^{-1}}
  \end{aligned}
\end{equation}
is an algebra isomorphism. Therefore $\bH_{s_e}$ is Morita equivalent to
$\bA$. The equivalence also preserves unitarity.
It follows that the only unitary representation with such
infinitesimal character is the trivial one. This case corresponds to
the unitarity of $X(\chi_\nu),$ $\nu$ purely imaginary, but $\nu\neq
0,\pi i/\log q.$

\subsection{}\label{sec:1.7} Now we can describe the role of the
graded algebra $\bH$ in the determination of the spherical unitary
dual of $G(\bR).$ The $SL(2)$ examples suggest that, in this
correspondence, under the appropriate technical assumption (``real
infinitesimal character''), there should be a matching of the unitary
representations. The technical notion for proving this correspondence is
that of {\it   petite   $K$-types} for real groups, which were defined in
\cite{Ba1} and studied further in \cite{Ba2} and \cite{BP}.

\smallskip

In theorem \ref{t:bc}, the unramified principal series
$X(\chi)$  correspond to the induced Hecke modules
$\bH\otimes_{\bA}\bC_\chi$,
where $\chi$ can be identified with an element of $\check\fh.$ The action
of $\bH$ on this induced module, which we will denote by $X(\chi)$ as well, is by
multiplication on the left. Therefore as a $\bC[W]$ module it is
isomorphic to the left regular representation,
\begin{equation}
  \label{eq:i1}
  \sum_{(\psi,V_\psi)\in \wht W} V_\psi\otimes V_\psi^*.
\end{equation}
A module $(\pi,V)$  of $\bH$ is called spherical if it has nontrivial
$W$-fixed vectors, \ie $V^W\ne (0).$ It is clear from (\ref{eq:i1})
that the $\bH$-module $X(\chi)$ contains the trivial $W$-type with
multiplicity $1$, so it has a unique spherical subquotient, denoted
again by $L(\chi).$
In theorem \ref{t:bc}, spherical modules for $\bH$
match spherical modules for $G(\bF).$

Assume that $\chi$ is real and dominant,  and such that
  the spherical quotient $L(\chi)$ is hermitian. Then $X(\chi)$ has
an invariant hermitian form so that the radical is the maximal proper invariant
subspace (so that the quotient is $L(\chi)$). Each space $V_\psi ^*$
inherits a hermitian form $\CA^\bH_\psi(\chi)$ depending continuously on
$\chi$ so that $L(\chi)$ is unitarizable if and only if all
$\CA^\bH_\psi(\chi)$, $\psi\in\wht W$, are positive semidefinite.
The set of \textit{relevant} $W$-types defined in \cite{Ba1},
\cite{Ba2}, \cite{BC2} is a minimal set of $W$-representations with
the property that  an $L(\chi)$ is unitary if and only if the form
$\CA_\psi$ is positive definite for all $\psi$ relevant.

\begin{example}For $G=SL(n)$, the
representations $\sigma$ of the Weyl
group $W=S_n$ are parameterized by partitions, and the set of relevant
$W$-types are formed only of the partitions with at most two parts
(\cite{Ba1}). In $G=E_8,$ there are $112$ $W$-types, and only nine are
called relevant (\cite{BC2}).
\end{example}

Now consider the case $\bF=\bR$ and $\chi$ a real unramified
character.   Assume that there is $w\in W$ so that
$w\chi=\chi^{-1}$. This is the condition $\chi$ must satisfy so that $X(\chi)$
admits an invariant hermitian form. Assume further that $\chi$ is
dominant so that $L(\chi)$ is the quotient by the maximal proper
invariant subspace. The maximal proper invariant subspace is also the
radical of the hermitian form, so that $L(\chi)$ inherits a
nondegenerate hermitian form. Then for every $K$-type $\mu,$ there is
a hermitian form $\CA^\bR_\mu(\chi)$ on $\Hom_K[\mu,X(\chi)].$ The
module $L(\chi)$ is unitary if and only if $\CA^\bR_\mu(\chi)$ is
positive semidefinite for all $\mu.$
By Frobenius reciprocity, we can interpret $\CA^\bR_\mu(\chi)$ as a
hermitian form on $(V_\mu^{^0\!H})^*$. Moreover $(V_\mu^{^0\!H})^*$ is a
  representation of $W$ (not necessarily irreducible). We denote it
  by $\psi_\mu.$

The petite $K$-types are $K$-types such that the real operator
$\CA^\bR_\mu(\chi)$ coincides with the Hecke operator
$\CA^\bH_{\psi_\mu}(\chi)$. The petite spherical $K$-types for split
$G(\bR)$ are studied in \cite{Ba1,Ba2}.

\begin{example} For the spherical principal series of $SL(2,\bR),$ (where
$K=SO(2)$ and $W=\bZ/2\bZ$) the only petite $SO(2)$-types appearing in
$X(\chi)$ are $(0)$, and $(\pm 2),$ which correspond to the trivial, and the
sign $W$-representations respectively.
\end{example}

The main result is that, for
every simple split group $G(\bR)$, every
relevant $W$-type occurs in $(V_\mu^{^0\!H})^*$ for a petite $\mu.$  The
consequence is that the set of unitary spherical representations with
real infinitesimal character for a
real split group $G(\bR)$  is contained in the set of unitary
spherical parameters with real infinitesimal character for
the corresponding $p$-adic group $G(\bF).$ In fact, for the split
classical groups,
\cite{Ba1} proves that these sets are equal. (This is false for
nonsplit groups, see \cite{Ba3} for complex groups, \cite{Ba4} and
\cite{BC3} for unitary groups.)

\subsection{}\label{sec:1.7a}The description of the spherical unitary
dual (for real infinitesimal character) when $G(\bF)$ is a simple
split group has a particularly nice form. The details are in section
\ref{sec:4}, where we also  explain  the main theorem.

View the spherical parameter $\chi$ as an element of
$\check\fh_{\bR}.$ Then, motivated by the results in \cite{KL}, one
can attach to $\chi$ a nilpotent $\vG$-orbit $\check\CO(\chi)$ in
$\check\fg$ as follows. Note first that {because $\chi\in
  \check\fg$ is semisimple}, the centralizer $\vG(\chi)$ is connected.
In fact, it is a Levi subgroup of $\vG.$ Define
\begin{equation}
\check\fg_1=\{X\in \fg: [\chi,X]=X\}.
\end{equation}
This vector space consists of nilpotent elements of $\check\fg.$ The
group $\vG(\chi)$ acts on $\check\fg_1$ via the adjoint representation with
finitely many orbits. Therefore there exists a unique orbit which is
open (dense) and we define $\check\CO(\chi)$ to be its $\vG$-saturation.

For every nilpotent orbit $\check\CO$ in $\check\fg$, we fix a Lie triple
$\{\check e,\check h,\check f\}$ (see \cite{CM}), and we denote by
$\fz(\check\CO)$
the centralizer in $\check\fg$ of this Lie triple. Then, in particular, any
parameter $\chi$ for which $\check\CO(\chi)=\check\CO$ can be written
(up to $W$-conjugacy) as
\begin{equation}\label{eq:1.6.2}
\chi=\check h/2+\nu,\quad \nu\in\fz(\check\CO).
\end{equation}
In fact,  according to \cite{BM1}, $\check\CO$ is the
unique maximal nilpotent orbit for which $\chi$ can be written in this way.

\begin{definition}\label{d:cs}
We define the {\it complementary series attached to } $\check\CO$,
denoted by $CS_{\check\fg}(\check\CO)$, to be the set of spherical
parameters $\chi$ such that $L(\chi)$ is unitary and
$\check\CO(\chi)=\check\CO.$ Clearly, the spherical unitary dual of
$G(\bF)$ is the disjoint union of all $CS_{\check\fg}(\check\CO).$
\end{definition}

When $\check\CO=0$, \ie the trivial nilpotent orbit in $\check\fg$,
the set of parameters $\chi$ such that $\check\CO(\chi)=0$ corresponds
to those characters $\chi$ for which $X(\chi)$ is
irreducible. By the results of \cite{V4} in the real case, \cite{BM4}
  in the adjoint $p$-adic case, these are precisely the
spherical principal series which are generic, \ie admit Whittaker models.
So $CS_{\check\fg}(0)$ consists of the unitary generic
spherical parameters for $G(\bF).$

The reducibility of the spherical principal series $X(\chi)$ is
well-known. $X(\chi)$ is reducible if and only if
\begin{align}\label{eq:hyp}
&\langle\check\al,\chi\rangle\in 2\bZ+1, \quad\text{in the real case},\\\notag
&\langle\check\al,\chi\rangle=1, \quad\text{in the $p$-adic case},
\end{align}
for all $\al\in R.$ Assuming $\chi$ is dominant,  $CS_{\check\fg}(0)$
is necessarily a subset of the complement in the dominant Weyl chamber
of $\check\fh_\bR$ of the arrangement of hyperplanes given
by (\ref{eq:hyp}).  The region of $\check\fh_\bR$ in the dominant Weyl
chamber on which all coroots
$\check\al$ are strictly less than $1$ is called the \textit{fundamental
alcove}. Moreover, any
region in $\check\fh_\bR$ which is conjugate under the affine Weyl group
to the fundamental alcove is called an alcove.  Since
the spherical principal series $X(\chi)$ is irreducible at $\chi=0,$
by unitary induction and a well-known deformation argument, the
parameters in the fundamental alcove must be in $CS_{\check\fg}(0).$

\medskip

We can now list the main results of \cite{Ba1} and
\cite{BC}. Earlier, for split $p$-adic groups, the 
spherical unitary dual in  type $A$ was determined in \cite{Ta},
for types $B,C,D$ in \cite{BM3}, and for $G_2$ in \cite{Mu}. For
split real groups, types $A$ and $G_2$ are part of \cite{V2} and
\cite{V3}, respectively.

\begin{theorem}[\cite{Ba1},\cite{BC}] Assume $\bF$ is a $p$-adic
  field. Recall that $G(\bF)$ is a simple split group, $\check\CO$ is
  a nilpotent   orbit in $\check\fg$, with a fixed Lie triple $\{\check e,\check
  h,\check f\}$ in $\check\fg$, whose centralizer is
  $\fz(\check\CO)$. The spherical
  complementary series $CS_{\check\fg}(\check\CO)$ are defined in
  definition \ref{d:cs}.

\begin{enumerate}

\item A spherical parameter $\chi=\check h/2+\nu,$ $\nu\in
  \fz(\check\CO)$ is in $CS_{\check\fg}(\check\CO)$ if and only if
  $\nu$ is in $CS_{\fz(\check\CO)}(0).$ There is one exception to this
  rule in type
  $F_4,$ one in type $E_7$, and six in type $E_8$ (tabulated in
  section \ref{sec:5}).

\item $CS_{\check\fg}(0)$ is the disjoint union of $2^\ell$ alcoves in
  $\fh_\bR$, where $\ell$ and the explicit description of the alcoves
  are listed in section \ref{sec:zero}.
\end{enumerate}

\end{theorem}

\begin{theorem}[\cite{Ba1}]
When $\bF=\bR$, the same
description of the spherical unitary dual holds for split classical
groups.
\end{theorem}

In sections \ref{sec:4.4} and \ref{sec:4.5}, we give an explicit description of the spherical
unitary parameters for classical groups in terms of the
Zelevinski-type strings introduced in \cite{BM3} and \cite{Ba1}. We
then explain, following \cite{Ba1}, how this allows one to write any given infinitesimal character
$\chi$ in the form (\ref{eq:1.6.2}), and to check if it parameterizes a
unitary representation as in the theorem above.

\smallskip

For the split exceptional real groups, the same theorem is expected
to hold, but at this point, by the correspondence relevant
$W$-types/petite $K$-types, we only know that
$CS_{\check\fg}(\check\CO)$ in the real case is a subset of
$CS_{\check\fg}(\check\CO)$ in the $p$-adic (\cite{Ba2,BC2}).

\subsubsection{} The  pictures below show the set of spherical unitary
parameters for $G(\bR)=SO(3,2)$ ($\bH=\bH(C_2)$) and
$G(\bR)=Sp(4,\bR)$ ($\bH=\bH(B_2)$) respectively. In both cases, the
set is partitioned into a disjoint union of complementary   series
attached to nilpotent $\check G$-orbits  in $\check\fg$; the orbits
are listed in the adjoint  tables. Note that the complementary series
attached to the trivial nilpotent orbit ($\check\CO=(1,\dots,1)$) is
the fundamental alcove. The complementary series attached to
the regular  nilpotent orbit consists of only one point, the trivial
representation.
%(Recall that this orbit is distinguished.)

%%%%%%%%%%%%%%%%%%%%%%%%%%%%%%%%%%%%%%%%%%%%%%%%%%%%%%%%%%%%%%%%%%%%%%%%%%%%%%%%%%%%%%%%%%%%%%%%%%%%%%%%%%%%%%%%%%%%%%%%%%%%%%%%%%%%%%%%%%%%%%%%%%%%%%%%%%%%%%%%%%%%%%%%%%%

\begin{figure}[h]

 \begin{picture}(200,350)(0,-200)

%%%%%%%%%%%%%%%%%%%%%%%%%%%%%%%%%%%%%%%%%%%%%%%%%%%%%%%%%%%%%%%%%%%%%%%%%%%%%
  \put(-67,100){{\color{black}{$\boxed{SO(3,2)}$}}}
  \put(-61,0){\tiny{0}}
\put(-60,10){\line(1,0){110}} \put(30,1){\tiny{$\nu_2=0$}}
\put(-60,10){\line(1,1){72}} \put(15,80){\tiny{$\nu_1=\nu_2$}}
\put(-30,10){\line(1,1){52}}
\put(-15,25){{\color{black}{\circle*{4}}}}
 \put(-30,10){\line(-1,1){15}}
  \put(-31,0){\tiny{1}}
 \put(-1,0){\tiny{2}}
 \put(0,10){\circle*{1}}
 \put(-45,10){\line(0,1){15}}
\put(-46,10){{\color{black}{\line(0,1){14}}}}
\put(-47,10){{\color{black}{\line(0,1){13}}}}
\put(-48,10){{\color{black}{\line(0,1){12}}}}
\put(-49,10){{\color{black}{\line(0,1){11}}}}
\put(-50,10){{\color{black}{\line(0,1){10}}}}
\put(-51,10){{\color{black}{\line(0,1){9}}}}
\put(-52,10){{\color{black}{\line(0,1){8}}}}
\put(-53,10){{\color{black}{\line(0,1){7}}}}
\put(-54,10){{\color{black}{\line(0,1){6}}}}
\put(-55,10){{\color{black}{\line(0,1){5}}}}
\put(-56,10){{\color{black}{\line(0,1){4}}}}
\put(-57,10){{\color{black}{\line(0,1){3}}}}
\put(-58,10){{\color{black}{\line(0,1){2}}}}
\put(-59,10){{\color{black}{\line(0,1){1}}}}
\put(-45,25){\line(1,0){85}}
 \put(-47,0){\tiny{$\frac 12$}}
\put(-62,29){\tiny{$(\frac 12,\frac 12)$}}\put(-20,15){\tiny{$(\frac
32,\frac 12)$}}

%%%%%%%%%%%%%%%%%%%%%%%%%%%%%%%%%%%%%%%%%%%%%%%%%%%%%%%%%%%%%%%%%%%%%%%%%%%%

\put(147,42){\tiny$(\frac 12,0)$} \put(146.5,70){\tiny$(\frac
12,\frac 12)$}
 \put(148.5,51.5){{\color{black}{\circle*{3}}}}
\put(146.5,62){{\color{black}{$\circ$}}}
\put(148,51){{\color{black}{\line(0,1){12}}}}
\put(148.5,51){{\color{black}{\line(0,1){12}}}}
\put(125,38){{\color{black}{\line(1,0){140}}}}

\put(125,80){{\color{black}{\line(1,0){140}}}}

\put(148,93){{\color{black}{\circle*{4}}}} \put(146,98){\tiny$(\frac
12,\frac 12)$} \put(148,130){{\color{black}{\circle*{4}}}}
\put(146,135){\tiny$(\frac 32,\frac 12)$}

\put(125,115){{\color{black}{\line(1,0){140}}}}

\put(125,150){{\color{black}{\line(1,0){140}}}}

\put(185,55){{\color{black}{$\check{\mathcal{O}}= (2,1,1)$}}}

\put(185,92){{\color{black}{$\check{\mathcal{O}}= (2,2)$}}}

\put(185,127){{\color{black}{$\check{\mathcal{O}}= (4)$}}}

 \put(153,0){\tiny$(\frac 12,0)$}
\put(152.5,28){\tiny$(\frac 12,\frac 12)$}
  \put(139,0){\tiny{0}}
\put(140,10){\line(1,0){14}} \put(140,10){\line(1,1){13.5}}
\put(153,7){{\color{black}{$\circ$}}}
\put(153,22){{\color{black}{$\circ$}}}
\put(155,16){{\color{black}{\line(0,1){2}}}}
\put(155,20){{\color{black}{\line(0,1){2}}}}
\put(155,12){{\color{black}{\line(0,1){2}}}}
\put(152.5,10.5){{\color{black}{\line(0,1){12}}}}
\put(150,10){{\color{black}{\line(0,1){10}}}}
\put(148,10){{\color{black}{\line(0,1){8}}}}
\put(146,10){{\color{black}{\line(0,1){6}}}}
\put(144,10){{\color{black}{\line(0,1){4}}}}
\put(142,10){{\color{black}{\line(0,1){2}}}}

\put(185,10){{\color{black}{$\check{\mathcal{O}}= (1,1,1,1)$}}}

\put(175,-3){{\color{black}{\line(0,1){153}}}}

\put(125,-3){{\color{black}{\line(0,1){153}}}}

\put(265,-3){{\color{black}{\line(0,1){153}}}}

\put(125,-3){{\color{black}{\line(1,0){140}}}}

%\end{picture}

%\[\]\[\]

 %\noindent\hskip-5cm
 %\begin{picture}(0,-40)(180,210)
   \put(-63,-70){{\color{black}{$\boxed{Sp(4,\bR)}$}}}

\put(-60,-160){\line(1,0){110}} \put(30,-170){\tiny{$\nu_2=0$}}
\put(-60,-160){\line(1,1){72}} \put(15,-90){\tiny{$\nu_1=\nu_2$}}
\put(-30,-160){\line(1,1){52}}
 \put(-30,-160){\line(-1,1){15}}
  \put(-60,-170){\tiny{0}}
  \put(-30,-170){\tiny{1}}
 \put(0,-170){\tiny{2}}
 \put(0,-160){\circle*{1}}
 \put(-30,-160){\line(0,1){30}}
\put(-46,-160){{\color{black}{\line(0,1){14}}}}
\put(-47,-160){{\color{black}{\line(0,1){13}}}}
\put(-48,-160){{\color{black}{\line(0,1){12}}}}
\put(-49,-160){{\color{black}{\line(0,1){11}}}}
\put(-50,-160){{\color{black}{\line(0,1){10}}}}
\put(-51,-160){{\color{black}{\line(0,1){9}}}}
\put(-52,-160){{\color{black}{\line(0,1){8}}}}
\put(-53,-160){{\color{black}{\line(0,1){7}}}}
\put(-54,-160){{\color{black}{\line(0,1){6}}}}
\put(-55,-160){{\color{black}{\line(0,1){5}}}}
\put(-56,-160){{\color{black}{\line(0,1){4}}}}
\put(-57,-160){{\color{black}{\line(0,1){3}}}}
\put(-58,-160){{\color{black}{\line(0,1){2}}}}
\put(-59,-160){{\color{black}{\line(0,1){1}}}}

\put(-30,-130){\line(1,0){65}}
\put(0,-130){{\color{black}{\circle*{4}}}}
\put(-68,-145){\tiny{$(\frac 12,\frac 12)$}}

\put(-2,-137){\tiny{$(2,1)$}}

\put(-45,-160){{\color{black}{\line(0,1){15}}}}
\put(-44,-160){{\color{black}{\line(0,1){14}}}}
\put(-43,-160){{\color{black}{\line(0,1){13}}}}
\put(-42,-160){{\color{black}{\line(0,1){12}}}}
\put(-41,-160){{\color{black}{\line(0,1){11}}}}
\put(-40,-160){{\color{black}{\line(0,1){10}}}}
\put(-39,-160){{\color{black}{\line(0,1){9}}}}
\put(-38,-160){{\color{black}{\line(0,1){8}}}}
\put(-37,-160){{\color{black}{\line(0,1){7}}}}
\put(-36,-160){{\color{black}{\line(0,1){6}}}}
\put(-35,-160){{\color{black}{\line(0,1){5}}}}
\put(-34,-160){{\color{black}{\line(0,1){4}}}}
\put(-33,-160){{\color{black}{\line(0,1){3}}}}
\put(-32,-160){{\color{black}{\line(0,1){2}}}}
\put(-31,-160){{\color{black}{\line(0,1){1}}}}

%\end{picture}

 %\noindent\hskip5.5cm

 %\begin{picture}(0,-40)(55,255)

%%%%%%%%%%%%%%%%%%%%%%%%%%%%%%%%%%%%%%%%%%%%%%%%%%%%%%%%%%%%%%%%%%%%%%%%%%%%%%%

\put(175,-183){{\color{black}{\line(0,1){152}}}}

\put(125,-183){{\color{black}{\line(0,1){152}}}}

\put(265,-183){{\color{black}{\line(0,1){152}}}}

\put(125,-31){{\color{black}{\line(1,0){140}}}}

%figura 1
\put(148,-50){{\color{black}{\circle*{4}}}}
\put(146,-45){\tiny{$(2,1)$}}

\put(185,-50){{\color{black}{$\check{\mathcal{O}}= (5)$}}}

% figura 2
 \put(148,-87){{\color{black}{\circle*{4}}}} \put(146,-82){\tiny{$(1,0)$}}

\put(185,-89){{\color{black}{$\check{\mathcal{O}}= (3,1,1)$}}}

\put(125,-65){{\color{black}{\line(1,0){140}}}}

% figura 3

\put(157,-138){\tiny{$(1,0)$}} \put(145,-110.5){\tiny{$(\frac
12,\frac 12)$}} \put(158.5, -132.8){{\color{black}{$\circ$}}}
\put(147,-117){{\color{black}{\circle*{4}}}}
\put(159.3,-130){{\color{black}{\line(-1,1){11}}}}

\put(185,-130){{\color{black}{$\check{\mathcal{O}}= (2,2,1)$}}}

\put(125,-102){{\color{black}{\line(1,0){140}}}}

%figura 4

\put(125,-143){{\color{black}{\line(1,0){140}}}}

 \put(157,-178){\tiny{$(1,0)$}}
\put(148,-153){\tiny{$(\frac 12,\frac 12)$}}
  \put(139,-178){\tiny{0}}
\put(140,-170){\line(1,0){21.5}}
 \put(140,-170){\line(1,1){11}}
\put(152,-170){{\color{black}{\line(0,1){10.5}}}}
\put(150,-170){{\color{black}{\line(0,1){10}}}}
\put(148,-170){{\color{black}{\line(0,1){8}}}}
\put(146,-170){{\color{black}{\line(0,1){6}}}}
\put(144,-170){{\color{black}{\line(0,1){4}}}}
\put(142,-170){{\color{black}{\line(0,1){2}}}}
\put(154,-170){{\color{black}{\line(0,1){10}}}}
\put(156,-170){{\color{black}{\line(0,1){8}}}}
\put(158,-170){{\color{black}{\line(0,1){6}}}}
\put(160,-170){{\color{black}{\line(0,1){4}}}}
\put(150,-161){{\color{black}{$\circ$}}}
\put(160.5,-172){{\color{black}{$\circ$}}}

\put(185,-168){{\color{black}{$\check{\mathcal{O}}= (1,1,1,1,1)$}}}

\put(125,-183){{\color{black}{\line(1,0){140}}}}

\end{picture}
\caption{Spherical unitary parameters for $SO(3,2)$ and $Sp(4,\bR)$}

\end{figure}
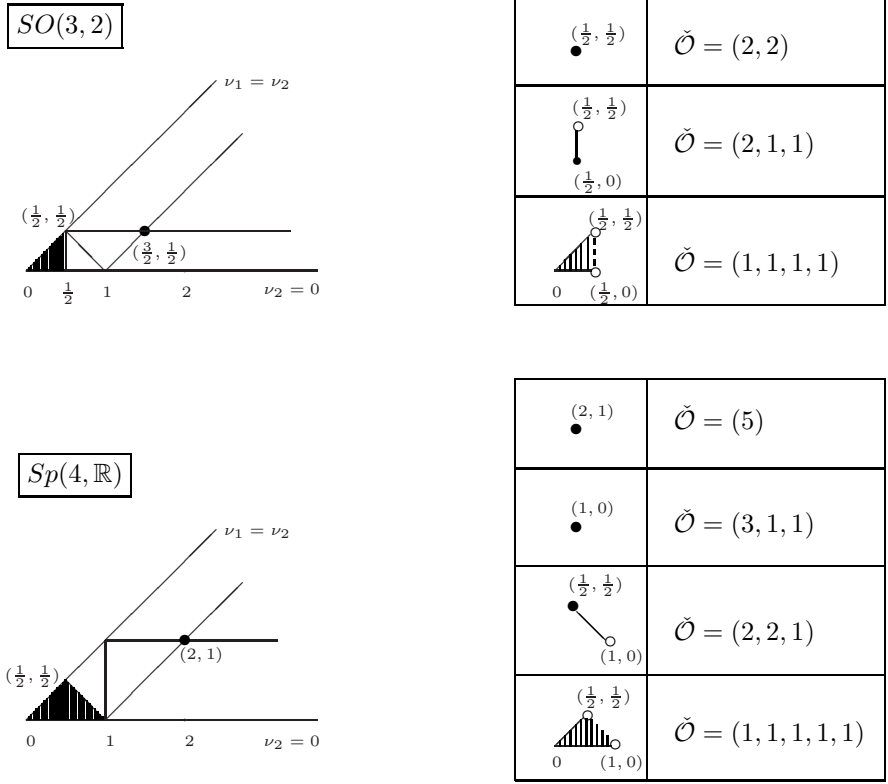

%%%%%%%%%%%%%%%%%%%%%%%%%%%%%%%%%%%%%%%%%%%%%%%%%%%%%%%%%%%%%%%%%%%%%%%%%%%%%%%%%%%%%%%%%%%%%%%%%%%%%%%%%%%%%%%%%%%%%%

Note that $Sp(4,\bR)\cong Spin(3,2),$ and indeed the pictures are the
same, up to the change of coordinates $(\nu_1,\nu_2)\mapsto
(\frac{\nu_1+\nu_2}2, \frac{\nu_1-\nu_2}2)$ (from the $Sp(4,\bR)$
parameters to the $SO(3,2)$ parameters).

\subsection{}\label{sec:1.8} Now we consider the ramified principal series.
If $\chi$ is a ramified character of $H(\bF)$, recall that $^0\!\chi$ is
the restriction to the compact part $^0\!H$.
To $^0\!\chi$, one associates
the homomorphisms
\begin{align}
&\hat\chi: \bO^\times\to \vH,& \text{when }\bF\text{ is
  $p$-adic},\\\notag
&\hat\chi: \bZ/2\bZ=\{1,-1\}\to \vH,& \text{when }\bF=\bR,
\end{align}
defined by the property that
\begin{equation}
\lambda\circ \hat\chi(z)=~^0\!\chi\circ\lambda(z),\  z\in\bO^\times\
(\text{or}\ z\in\bZ/2\bZ\subset \bR^\times), \text{ for all }\lambda\in\C Y,
\end{equation}
where in the left hand side $\lambda$ is regarded as a character of
$\vH$ and in right hand side as a cocharacter of $H(\bF).$
In the real case, the image of $\hat\chi$ can clearly be identified
with a single semisimple (elliptic) element of order $2$ in $\vH.$
\begin{comment}In
this case, if $\check\al\in \C Y$ is a coroot of $G$, the classical
notation (see section \ref{sec:2.1a}) is $\check\al(-1)=m_\al\in M$ ($M=~^0\!H$).
\end{comment}

Define
\begin{equation}
\vG(\hat\chi)=\text{ the centralizer of the image of }\hat\chi\text{ in }\vG.
\end{equation}
In the $p$-adic case, it is proved in \cite{Ro} that there exists an
elliptic semisimple element $s$ of $\vG$ such that the the identity
component of $\vG(\hat\chi)$ equals the identity component of the
centralizer of $s$ in $\vG.$ Moreover, the component group
$\check R_{^0\!\chi}$ of $\vG(\hat\chi)$ is abelian. In the real case this is a
classical result (Knapp-Stein): the group $\vG(\hat\chi)$ is the one
defined by the {\it good roots} with respect to $^0\!\chi$, while
$\check R_{^0\!\chi}$ is the dual of the $R$-group. (For a definition of good
(co)roots, see section \ref{sec:fine} if $\bF$ is real, and section
\ref{sec:3.0a} if $\bF$ is $p$-adic.)

Let $G(^0\!\chi)$ be the connected split subgroup of $G(\bF)$ whose
complex dual is the identity component of $\vG(\hat\chi)$.
The methods described in this paper give a way to compare the unitary
Langlands quotients of the ramified principal series
$X(\chi)=\Ind_{B(\bF)}^{G(\bF)}(\chi\otimes 1)$ with the spherical
unitary dual  of $G(^0\!\chi).$

\medskip

{\noindent{\bf Remark.} One may consider the extension
$G(^0\!\chi)'$ of $G(^0\!\chi)$ by the (abelian) $R$-group $R_{\
{^0}{\chi}}\cong \vR_{\ ^0\!\chi},$ and say that a representation
$(\pi,V)$ of $G(^0\!\chi)'$ is \textit{quasi-spherical} if the
restriction of $\pi$ to $G(^0\!\chi)$ is spherical. A refinement
would be to explain the unitarizable Langlands quotients of the
ramified principal series via the quasi-spherical unitary dual of
$G(^0\!\chi)'$, rather then the spherical unitary dual of
$G(^0\!\chi).$}

\subsection{}\label{sec:1.10} We give some examples illustrating the
construction of $G(^0\!\chi)$ when $\chi$ is ramified.

\subsubsection{}\label{sec:1.10.1} The first example is the
nonspherical principal series of $SL(2,\bR).$ Let $^0\!\chi$ be the sign
character of $\bZ/2\bZ.$ We can write $\chi=~^0\!\chi\otimes |\cdot |^\nu$, and
we assume that  $\nu$ is real.   It is well-known that $X(\chi)$ is
reducible when $\nu=0$, and then it splits into a sum of two unitary
representations $L(\chi)=X(\chi)=L(\chi)_1\oplus L(\chi)_{-1}$ (the
limits of discrete series). Moreover, $\nu=0$ is
the only value of the parameter $\nu$ for which $L(\chi)$ is
unitary. In this case, the nontrivial element in the image of
$\hat\chi$ is $\left(\begin{matrix} 1&0\\0&-1\end{matrix}\right)$ in
  $PGL(2,\bC)$. The centralizer $\vG(\hat\chi)$ has two components:
  the identity component is the diagonal $\vH,$ and the other is
  formed by the
  matrices with $0$ on the diagonal.  Then $G(^0\!\chi)$ is the diagonal
  torus $H(\bR)=\bR^\times$, and the only unramified real unitary
  character is the trivial one. This can be formulates as: $L(\chi)$
  is unitary if and only if $\chi$ is an unramified unitary parameter
  of $G(^0\!\chi).$ One refines this statement by introducing
  $G(^0\!\chi)'.$
The R-group is
$\bZ/2\bZ$, and  $G(^0\!\chi)'=H(\bR)\rtimes \bZ/2\bZ=O(1,1),$ where
the nontrivial $\bZ/2\bZ$ element acts by flipping the diagonal
entries. There are two quasi-spherical unitary representations of $G(^0\!\chi)'$
corresponding to the trivial, respectively the sign representations of
$Z/2\bZ.$

\subsubsection{} An  example of a nontrivial character $^0\!\chi$ for
  $p$-adic $SL(2,\bF)$, $p$ odd, is as follows. First set
\begin{equation}\label{eq:1.10.1}
sgn:\bF_q^\times\to \{\pm1\},\quad  sgn(a)=\left\{\begin{matrix} 1,&
a\in (\bF_q^\times)^2,\\
-1, &a\notin (\bF_q^\times)^2\end{matrix}\right.
\end{equation}
Then ${^0}{H}\cong \bO^\times,$  and the analog of (\ref{eq:1.4.1}) is
\begin{equation}\label{eq:1.10.2}
1\longrightarrow 1+\C P\longrightarrow \bO^\times\longrightarrow
\bF_q^\times\longrightarrow 1.
\end{equation}
The character $^0\!\chi$ we consider is the pull back of $sgn$ from
$\bF_q$ to $\bO^\times.$

\smallskip

 Later in the paper, in section \ref{sec:3.0a}, we will define a
  correspondence between the minimal principal series for
$G(\bR)$ and certain minimal principal series of $p$-adic $G(\bF)$,
by sending the sign character of $\bZ/2\bZ$ to a fixed nontrivial
quadratic character $^0\!\chi$ of $\bO^\times$.

\subsubsection{} When $G$ has connected center, then $\vG$ has
simply-connected derived
subgroup, and $\vG(\hat\chi)$ is connected by a well known theorem of
Steinberg. Using the results of
\cite{Ro} and \cite{KL}, the Langlands classification for the ramified principal
series in the $p$-adic case  was obtained in this case in
\cite{Re}. When $G$ has connected center, $G(^0\!\chi)'=G(^0\!\chi)$
is the (endoscopic) group of $G(\bF)$ we consider, and the intent is
to match the unitary Langlands quotients
of the ramified principal series $\Ind_{B(\bF)}^{G(\bF)}(\chi\otimes
1)$ with the spherical unitary dual of $G(^0\!\chi).$

For example, when
$G=SO(2n+1)$, the groups which appear as $G(^0\!\chi)$ in our cases are precisely the
elliptic (in the sense of not being contained in any Levi subgroup)
endoscopic split groups $SO(2(n-m)+1,\bF)\times SO(2m+1,\bF).$

\subsubsection{}On the other hand, when $G=Sp(2n),$ the groups which
appear are the split groups
$G(^0\!\chi)=Sp(2n-2m,\bF)\times SO(2m,\bF),$ and
$G(^0\!\chi)'=Sp(2n-2m,\bF)\times O(2m,\bF).$
When $\bF=\bR$, the case of nonspherical minimal principal series for $Sp(2n,\bR)$
is presented in more detail in sections \ref{sec:2.5} and
\ref{sec:3a.7.1}.

\medskip

For example, consider $G=Sp(4).$ If $\bF=\bR$, then $^0\!H=(\bZ/2\bZ)^2$,
so there are four characters $^0\!\chi$ of $^0\!H$: $^0\!\chi=triv\otimes triv,$
$triv\otimes\chi_0,$ $\chi_0\otimes triv,$ and $\chi_0\otimes \chi_0$,
where $\chi_0=sgn.$  Since $triv\otimes\chi_0$ and $\chi_0\otimes
triv$ are conjugate, there are really only three cases.   When $\bF$
is $p$-adic, $^0\!H=(\bO^\times)^2.$ We fix a quadratic, nontrivial
character $\chi_0$ of $\bO^\times$, and only consider the similar three
$^0\!\chi$'s: $triv\otimes triv,$ $triv\otimes \chi_0,$ and
$\chi_0\otimes\chi_0$ respectively. Then we have:

\begin{center}
\begin{tabular}{|c|c|c|}
\hline
$^0\!\chi$ &$G(^0\!\chi)$ &$G(^0\!\chi)'$\\
\hline
$triv\otimes triv$ &$Sp(4,\bF)$ &$Sp(4,\bF)$\\
$triv\otimes \chi_0$ &$Sp(2,\bF)\times SO(2,\bF)$ &$Sp(2,\bF)\times O(2,\bF)$\\
$\chi_0\otimes\chi_0$ &$SO(4,\bF)$ &$O(4,\bF)$\\
\hline
\end{tabular}
\end{center}

The first case is the unramified one. The notation in the
table may be confusing when $\bF=\bR$: we mean the split orthogonal
groups, in classical notation, $SO(n,n)$ and $O(n,n).$

\subsection{}
Similarly to the unramified case,  we define  $\bH(^0\!\chi)$ to be the
graded Hecke algebra for  $\vG(^0\!\chi)$,
   and we denote by $\bH'(^0\!\chi)$ the extension of
$\bH(^0\!\chi)$ by the dual R-group. (These extensions of the graded
   Hecke algebra are defined in section
\ref{sec:3.7}).

Notice that if  $^0\!\chi$ is trivial, \ie if the principal series
$X(\chi)$ is spherical, then the R-group is
trivial. In this case $\bH'(^0\!\chi)=\bH(^0\!\chi)$ is the (usual) Hecke
algebra $\bH$ attached to the coroot system of $G(\mathbb{F})$.
Also, when $G$ has connected center, or more generally whenever
the R-group of $^0\!\chi$ is trivial, $\bH'(^0\!\chi)=\bH(^0\!\chi).$

In the $p$-adic case, at least when $G$ has connected center, the
generalization of theorem \ref{t:bm} and of the \cite{BM2} translation
of unitarity to the graded Hecke algebra appears to hold. We will not consider
this problem in the present paper, but we hope to pursue it in future work.

\smallskip

The generalization of the notion of petite $K$-types in \cite{BP}
provides a similar matching of intertwining operators for the ramified
principal series of the real split group on one hand, with
intertwining operators for the Hecke algebra $\bH(^0\!\chi)$
on the other. For example, for the nonspherical principal series of $SL(2,\bR),$
the only petite $SO(2)$-types are the two {\it fine} $K$-types
(section \ref{sec:fine}), $(\pm1).$  The difficulty of defining and
computing the petite $K$-types for nonspherical principal series in
general lies in the fact that one needs to capture at the same time the specifics of
both the spherical and nonspherical principal series of $SL(2,\bR).$
This is realized in \cite{BP}, and we postpone the technical details
until section \ref{sec:3a}.  Instead we present the example of $Sp(4,\bR).$

\subsection{The example of $\mathbf{Sp(4,\bR)}$}\label{sec:1.12} In this example, we
use the more customary notation
$\delta$ for $^0\!\chi$, and $\hat\delta$ for $\hat\chi$ (defined in
\ref{sec:1.8}). {There are four minimal principal series $X(\delta,\nu)$
($M=~^0\!H=(\bZ/2\bZ)^2$): $\delta_0=triv\otimes triv,$
$\delta_1^+=triv\otimes sgn$, $\delta_1^-=sgn\otimes triv$, and
$\delta_2=sgn\otimes sgn.$} The maximal compact group is $K=U(2)$, whose
representations are parameterized by pairs of integers $(a,b),$ $a\ge
b.$  The ramified character $\chi$ is written $\chi=\delta\otimes\nu,$
 where $\nu=(\nu_1,\nu_2),$ $\nu_1,\nu_2\in\bR.$ We assume that
  $\nu_1\ge\nu_2\ge 0$, \ie that the parameter is (weakly) dominant.

As in the spherical case, it is well known when
  $X(\delta,\nu)$ admits an invariant hermitian form. Since the
  parameter $\nu$ is assumed real, the condition is
  that there exists $w\in W$ such that
\begin{equation}\label{eq:1.12.1}
w\cdot \delta\cong\delta,\text{ and }w\cdot \nu=-\nu.
\end{equation}
{In the $Sp(4,\bR)$ cases, all quotients $L(\delta,\nu)$ are hermitian.}

We denote the Weyl groups of $G(\delta)$ and $G(\delta)'$, by
$W^0_{\delta}$ and $W_{\delta}$ respectively. For every $K$-type
$\mu,$ similarly to the spherical case, one has an operator
$\CA_\mu^\bR(\delta,\nu)$ on the
space \begin{equation}\label{eq:hom}\Hom_{M}[\mu|_M:\delta],\end{equation} whose
signature we would like to compute. This space carries a
representation of $W_\delta$ (and $W_\delta^0$), which we denote as
before by $\psi_\mu.$ The conditions in section
\ref{sec:3a} used to define petite $K$-types are such that this
operator is the same as the ``$p$-adic'' operator on $\psi_\mu$ for the
Hecke algebra $\bH(\delta)$ (or more precisely, the extended
version $\bH'(\delta)).$

The tables with examples of petite $K$-types $\mu$, their operators
$\CA^\bR_\mu(\nu_1,\nu_2)$,  and the set of unitary Langlands
subquotients of minimal principal series are next. In the tables
below, by ``mult.'' we mean the dimension of the space (\ref{eq:hom}).

\subsubsection{$\delta_0=triv\otimes triv$}\
This is the spherical principal series.  In this case
$\vG(\hat\delta_0)=SO(5,\bC)$. We may assume $\nu_1\ge\nu_2\ge 0.$ Then
$W_\delta=W_\delta^0=W=W(B_2).$ The representations of $W(B_2)$ are
parameterized by pairs of partitions of total sum $2$. There are $4$
one dimensional representations, and one two-dimensional, labeled
$(1)\times (1).$

\bigskip
\begin{tabular}{|c|c|c|c|}
\hline
$\mu$ &mult.  &$\psi_\mu\in
  \widehat{W}$ &$\CA^\bR_\mu(\nu)$\\
\hline
$(0,0)$ &$1$ &$(2)\times (0)$ &$1$   \\\hline
$(1,-1)$ &$1$ &$(1,1)\times (0)$
&$\frac{(1-(\nu_1+\nu_2))(1-(\nu_1-\nu_2))}{(1+(\nu_1+\nu_2))(1+(\nu_1-\nu_2))}$\\\hline
$(2,2)$ &$1$ &$(0)\times (2)$
&$\frac{(1-\nu_1)(1-\nu_2)}{(1+\nu_1)(1+\nu_2)}$\\\hline
$(2,0)$ &$2$ &$(1)\times (1)$
&$m$\\
\hline
\end{tabular}

\bigskip
\noindent where $m=\frac
1{c(\nu_1,\nu_2)}\begin{tiny}\left(\begin{matrix}(1+\nu_2)[(1+\nu_1)+(1-\nu_1)(\nu_1^2-\nu_2^2)]
  &2\nu_1(1-\nu_2^2)\\ 2\nu_1(1-\nu_2^2) &(1-\nu_2)[(1-\nu_1)+(1+\nu_1)(\nu_1^2-\nu_2^2)]      \end{matrix}
\right)\end{tiny},$ and $c(\nu_1,\nu_2)=(1+\nu_1)(1+\nu_2)(1+(\nu_1-\nu_2))(1+(\nu_1+\nu_2)).$

The spherical quotient $L(\nu)$ is unitary if and only if
$\nu_1+\nu_2\le 1,$ or $(\nu_1,\nu_2)=(2,1).$

\subsubsection{$\delta^+_1=triv\otimes sgn$
and $\delta_1^-=sgn\otimes triv$ }\
In this case,
$\vG(\hat\delta^+_1)\cong S[O(3,\bC)\times O(2,\bC)],$ and
$\vG(\hat\delta^-_1)\cong S[O(2,\bC)\times O(3,\bC)],$
$W_{\delta_1^\pm}^0=W(A_1)$ and
$W_{\delta_1^\pm}=W_{\delta^\pm_1}^0\times W(A_1).$

The two principal series are $X(\delta_1^+,
(\nu_1,\nu_2))$ and $X(\delta_1^-,(\nu_1,\nu_2))$ with $\nu_1\ge\nu_2\ge 0.$ If
$\nu_1=\nu_2,$ they have the same
Langlands quotient. There is a difference however: in the case of
$\delta_1^+,$ the parameter $\nu_1$ corresponds to the $SL(2)$ in
$G(\delta_1^+)$, and $\nu_2$ to the $SO(1,1),$ while in the case of
$\delta_1^-$, $\nu_1$ corresponds to $SO(1,1)$ and $\nu_2$ to
$SL(2)$ in $G(\delta_1^-).$

\medskip

1) We consider first $\delta_1^+$ and $\nu_1\ge\nu_2\ge 0.$
We need to
distinguish between two cases:

\medskip
\noindent (i) $\nu_2>0.$ The Langlands quotient $L(\delta_1^+,\nu)$ is
irreducible. The operators are as in the following table.

\bigskip
\begin{tabular}{|c|c|c|c|c|}
\hline
$\mu$ &mult. &$\psi_\mu\in \widehat{W_\delta^0}$ &$\psi_\mu\in
  \widehat{W_\delta}$ &$\CA^\bR_\mu(\nu)$\\
\hline
$(1,0)$ &$1$ &$triv$ &$triv\times triv$ &$+1$\\\hline
$(0,-1)$ &$1$ &$triv$ &$triv\times sgn$ &$-1$\\\hline
$(2,1)$ &$1$ &$sgn$ &$sgn\times triv$ &$\frac{1-\nu_1}{1+\nu_1}$\\\hline
$(-1,-2)$ &$1$ &$sgn$ &$sgn\times sgn$ &$-\frac{1-\nu_1}{1+\nu_1}$\\

\hline
\end{tabular}

\medskip
Thus $L(\nu)$ is \textbf{not unitary} for $\nu_2>0.$

\bigskip

\noindent (ii) $\nu_2=0.$ Then $X(\delta_1^+,(\nu_1,0))$ is a direct sum
of two modules, $X_{(1,0)}(\delta_1^+,(\nu_1,0))$ and
$X_{(0,-1)}(\delta_1^+,(\nu_1,0))$. The Langlands quotient is a direct
sum of two irreducible modules, $L_{(1,0)}(\delta_1^+,(\nu_1,0))$ and
$L_{(0,-1)}(\delta_1^+,(\nu_1,0)),$ distinguished by the fact that the
former contains the $K$-type $(1,0)$ (and $(2,1)$), and the latter
contains $(0,-1)$ (and $(-1,-2)$). For $\nu_1\ge 0,$ $L_{(1,0)}$ is the
unique irreducible quotient of $X_{(1,0)},$ and $L_{(0,-1)}$ is the
unique irreducible quotient of $X_{(0,-1)}.$ The intertwining operator
is as in the table above. Thus the
Langlands quotients are \textbf{not unitary} for $1<\nu_1.$ On the
other hand, $X_{(1,0)}=L_{(1,0)}$ and $X_{(0,-1)}=L_{(0,-1)}$ for
$0\le\nu_1<1.$ At $\nu_1=0,$ $X(\delta_1^+,(0,0))$ is unitarily induced
from a unitary character $\delta_1^+$, and equal to the direct sum of $L_{(1,0)}$
and $L_{(0,-1)}.$ Thus $L_{(1,0)}(0,0)$ and $L_{(0,-1)}(0,0)$ are
both unitary. It follows that $L_{(1,0)}$ and $L_{(0,-1)}$ are
unitary for $0\le \nu_1\le 1$ by the continuity of the hermitian form
in the parameter $\nu_1.$

In this case, $G(\delta_1^+)=SL(2,\bR)\times SO(1,1).$ A
spherical parameter for this group, $(\nu_1,\nu_2)$ with $\nu_1$ on
the $SL(2,\bR)$ and $\nu_2$ on the $SO(1,1)$, is hermitian if and only
if $\nu_2=0.$ The  representations $L(\delta_1^+,(\nu_1,0))$ are unitary
if and only if the corresponding spherical $L(\nu_1)$  on $SL(2,\bR)$ is unitary.

\bigskip

{

2) Now consider $\delta_1^-,$ and $\nu_1\ge\nu_2\ge 0.$ The intertwining
operators are as in the table for $\delta_1^+$ above. But in this case,
the Langlands quotient $L(\delta_1^-,(\nu_1,\nu_2))$ is irreducible if
and only if $\nu_1>0,$ which implies that unless
$(\nu_1,\nu_2)=(0,0),$ the Langlands quotient is not unitary. At
$(0,0)$, there are two Langlands quotients, same as for $\delta_1^+,$ and
they are both unitary.

The group $G(\delta^-_1)$ is $SO(1,1)\otimes SL(2,\bR),$ with $\nu_1$
corresponding to  $SO(1,1)$ and $\nu_2$ to $SL(2,\bR).$

\medskip

The conclusion is that a hermitian parameter for $\delta_1^\pm$
\textit{with $\nu_1\ge \nu_2\ge 0$} is unitary if and only if the
corresponding parameter for $G(\delta_1^\pm)$ is unitary.
}

\subsubsection{$\delta_2=sgn\otimes sgn$}\label{sec:del2} \
In this case, $\vG(\hat\delta_2)=S[O(4,\bC)\times O(1,\bC)],$ and
$W_\delta^0=W(A_1)\times W(A_1),$ $W_\delta=W_\delta^0\rtimes
(\bZ/2\bZ)\cong W(B_2).$ We may assume $\nu_1\ge\nu_2\ge 0.$ There are two cases:

\bigskip
\noindent (i) $\nu_2>0$. The Langlands quotient is irreducible, and the
operators are as in the table

\bigskip

\begin{tabular}{|c|c|c|c|c|}
\hline
$\mu$ &mult. &$\psi_\mu\in \widehat{W_\delta^0}$ &$\psi_\mu\in
  \widehat{W_\delta}$ &$\CA^\bR_\mu(\nu)$\\
\hline
$(1,1)$ &$1$ &$triv\otimes triv$ &$(2)\times (0)$ &$1$\\\hline
$(-1,-1)$ &$1$ &$triv\otimes triv$ &$(0)\times (2)$ &$1$\\\hline
$(2,0)$ &$1$ &$sgn\otimes sgn$ &$(1,1)\times (0)$
&$\frac{(1-(\nu_1+\nu_2)(1-(\nu_1-\nu_2)}{(1+(\nu_1+\nu_2)(1+(\nu_1-\nu_2)}$
\\\hline
$(0,-2)$ &$1$ &$sgn\otimes sgn$ &$(0)\times (1,1)$
&$\frac{(1-(\nu_1+\nu_2)(1-(\nu_1-\nu_2)}{(1+(\nu_1+\nu_2)(1+(\nu_1-\nu_2)}$
\\\hline
$(1,-1)$ &$2$ &$triv\otimes sgn$ &$(1)\times (1)$
&$\begin{pmatrix}\frac{1-(\nu_1-\nu_2)}{1+(\nu_1-\nu_2)} &0\\ 0&\frac{1-(\nu_1+\nu_2)}{1+(\nu_1+\nu_2)}\end{pmatrix}$\\
         &    &$+sgn\otimes triv$ &  &\\
\hline
\end{tabular}

\bigskip
The Langlands quotient is unitary if and only if $0<\nu_2\le 1-\nu_1.$

\bigskip

\noindent (ii) $\nu_2=0.$ In this case $X(\delta_2,(\nu_1,0))$ is the
direct sum of
modules $X_{(1,1)}(\nu_1,0)$ and $X_{(-1,-1)}(\nu_1,0)$.
The Langlands subquotient is the direct sum of two irreducible modules
$L_{(1,1)}(\nu_1,0)$ and
$L_{(-1,-1)}(\nu_1,0)$, characterized by the fact
that the former contains
the $K$-type $(1,1)$ (and $(2,0)$, and a copy of $(1,-1)$), and the
latter contains $(-1,-1)$ (and $(0,-2)$, and a copy of
$(1,-1)$). For $\nu_1\ge 0,$ $L_{(1,1)}(\nu_1,0)$ is the
unique irreducible quotient of $X_{(1,1)},$ and
$L_{(-1,-1)}(\nu_1,0)$ the unique irreducible subquotient
of $X_{(-1,-1)}.$ The intertwining operator is obtained by taking $\nu_2\to
0$ in the table above. This means in particular that the operator on
 $(1,-1)$ is $\frac{1-\nu_1}{1+\nu_1}Id.$ Thus both
 $L_{(1,1)}(\nu_1,0)$ and
 $L_{(-1,-1)}(\nu_1,0)$ are \textbf{ not unitary} for  $1<\nu_1.$
At $\nu_1=0,$ the module $X(\delta_2,(0,0))$ is unitarily induced from
the unitary character $\delta_2,$ so both $L_{(1,1)}(0,0)$ and
$L_{(-1,-1)}(0,0)$ are unitary. Since both $X_{(1,1)}$ and $X_{(-1,-1)}$
are irreducible for $0\le\nu_1<1,$ the unitarity of $L_{(1,1)}$ and
$L_{(-1,-1)}$ for $0\le \nu_1\le 1$ follows from the continuity of the
hermitian form in the parameter $\nu_1$.

 In this case, $G(\delta)=SO(2,2).$ The unitary
  representations $L(\delta_2,(\nu_1,\nu_2))$ match the spherical
  unitary dual of this group, which in the
coordinates used in this subsection, consists of parameters
$(\nu_1,\nu_2)$ satisfying $0\le\nu_2\le 1-\nu_1.$

\subsubsection{}

% We give a more detailed explanation of the role of the
%  group $G(\delta)'.$

 As already alluded to in the
  cases above, every principal series $X(\delta_i,\nu)$ contains
  certain special $K$-types, called {\it fine} (see section
  \ref{sec:fine}). They are as follows:

\bigskip

\begin{center}
\begin{tabular}{|c|c|}
\hline
$\delta_0$ &$(0,0)$\\
\hline
$\delta_1^\pm$ &$(1,0),\ (0,-1)$\\
\hline
$\delta_2$ &$(1,1),\ (-1,-1)$\\
\hline
\end{tabular}
\end{center}

\bigskip

In conclusion, the matching for unitary nonspherical principal series
of $Sp(4,\bR)$ can be formulated as follows.

\begin{proposition}\label{r:m}
Assume that $\nu=(\nu_1,\nu_2)$ is weakly
dominant, and $(\delta,\nu)$  is hermitian. Fix a fine $K$-type $\mu$ and
let $L_\mu(\delta,(\nu_1,\nu_2))$ denote the
  irreducible Langlands quotient which contains the fine $K$-type  $\mu.$

Then $L_\mu(\delta,(\nu_1,\nu_2))$
  is unitary if and only if $(\nu_1,\nu_2)$
parameterizes a unitary spherical representation for the split group
$G(\delta).$
\end{proposition}

\subsection*{Remarks}\

\noindent (1) The spherical unitary dual for $G(\delta)$ can be read from
the theorems in \ref{d:cs}.

\smallskip

\noindent (2) Even though $\delta_1^+$ and $\delta^-_1$ are conjugate, the
corresponding unitary duals are not the same. This should not be
surprising in view of the Langlands classification of admissible representations. The fact that in the classification,
we must assume $(\nu_1,\nu_2)$ is dominant implies that
$(\delta_1^+,(\nu_1,\nu_2))$ is conjugate to
$(\delta^-_1,(\nu_1',\nu_2'))$ if and only if $\nu_1=\nu_2=\nu_1'=\nu_2'.$

\smallskip

{
\noindent (3) The correspondence $(\delta, (\nu_1,\nu_2))\mapsto
  (\nu_1,\nu_2)$ from parameters of $G$ to spherical parameters of $G(\delta)$
does not necessarily take hermitian parameters to hermitian parameters.
Again in the $\delta_1^-$ example,  any
parameter $(\nu_1,\nu_2)$ is hermitian for $G$. For the factor
$SO(1,1)=\bR^\times,$ the hermitian dual of a real character $z\mapsto
|z|^{\nu_1}$ is $z\mapsto |z|^{-\nu_1}$, so the only hermitian
character is the unitary one, \ie for $\nu_1=0.$ One of the roles of
the group $G(\delta)'$ is to fix this problem, as well as to take into
account that the Langlands quotient is not always irreducible. See also the
remark in \ref{sec:1.8}. For this example,
$G(\delta_1^-)'=O(1,1)\times SL(2,\bR).$ The maximal compact group of
$O(1,1)$ is $O(1)\times O(1)$, and for $SO(1,1)$ it is $S[O(1)\times
O(1)].$  We describe the  quasi-spherical
representations of $O(1,1)=\bR^\times\rtimes \bZ/2\bZ$. Explicitly we
realize $O(1,1)$ as the subgroup preserving the form $J=\left(\begin{matrix}
  0&1\\1&0\end{matrix}\right)$. The quasispherical representations
have restrictions to $O(1)\times O(1)$ which
consist of $triv\otimes triv$ and $det\otimes det$. The analogues of
the $X(\delta,\nu)$ are the two dimensional  representations
$\kappa_{\nu_1}$, given by the formulas
\begin{equation}
\kappa_{\nu_1}:\ \left(\begin{matrix}z&0\\0&z^{-1}\end{matrix}\right)\mapsto
\left(\begin{matrix}|z|^{\nu_1}&0\\0&|z|^{-\nu_1}\end{matrix}\right), \text{
  and }\left(\begin{matrix}0&z\\z^{-1}&0\end{matrix}\right)\mapsto
\left(\begin{matrix}0&|z|^{\nu_1}\\|z|^{-\nu_1}&0\end{matrix}\right),\
z\in\bR^\times.
\end{equation}
The restriction of $\kappa_{\nu_1}$ to $SO(1,1)$ is the sum of the
characters $|\ |^{\nu_1}$ and $|\ |^{-\nu_1}$, and the
restriction to $O(1)\times O(1)$ is the sum of
$triv\otimes triv$ and $det\otimes det$.
The representation $\kappa_{\nu_1}$ is irreducible for $\nu_1\ne 0,$
and decomposes into
the sum of the \textit{trivial} and the \textit{determinant}
representation for $\nu_1=0.$ So the $\kappa_{\nu_1}$ for $\nu_1\ne 0$
and the $trivial$ and $determinant$ representations for $\nu_1=0$ can
be viewed as the analogues of the $L_\mu(\delta,\nu).$  We can make
this precise by pairing the two fine $K$-types of $X(\delta_1^-,\nu)$,
$(1,0)$ and $(0,-1)$ with
$triv\otimes triv$ and $det\otimes det$ of $O(1)\times O(1)$
respectively. Furthermore every $\kappa_{\nu_1}$ is hermitian, but not
unitary if $\nu_1>0.$ Thus we have set up a correspondence between
parameters $L_\mu(\delta_1^-,(\nu_1,\nu_2))$ for $G$ with parameters
$L_{\phi(\mu)}((\nu_1,\nu_2))$  which is 1-1 onto parameters for
$G(\delta_1^-)'$ satisfying $\nu_1\ge \nu_2$ which matches hermitian
parameters with hermitian parameters, and unitary parameters with
unitary parameters. The only drawback is that there is a  choice in
the matching of fine $K$-types which is not canonical.

A similar result can be stated for $\delta_1^+$ and $\delta_2.$
We summarize this discussion in a proposition.

\begin{comment}
Note that the admissible parameters $\nu$ for
  $X(\delta_1^-,\nu)$ satisfy $\nu_1\ge\nu_2\ge 0,$ while the
  quasi-spherical parameters for $G(\delta_1^-)'=O(1,1)\times
  SL(2,\bR)$  satisfy $\nu_1\ge 0,\nu_2\ge 0.$ But if one assumes that
  $\nu_1\ge\nu_2\ge 0$ in both cases (that is, $\nu$ is such that it is weakly
  dominant for $G$),
 then there is a 1-1 correspondence between the hermitian parameters as
well as the unitary parameters for $\delta_1^-$ and $G(\delta_1^-)'$.
\end{comment}

\begin{proposition}Let $\phi$ denote the
(non-canonical) correspondence for fine $K$-types. Assume that $\nu$ is weakly
dominant  for $G$. Let $L_{\phi(\mu)}(G(\delta)',\nu)$ be
the unique irreducible quasi-spherical representation of $G(\delta)'$ which
contains $\phi(\mu).$ Then:
\begin{enumerate}
\item $L_\mu(\delta,\nu)$ is hermitian if and only if
  $L_{\phi(\mu)}(G(\delta)',\nu)$ is hermitian.
\item $L_\mu(\delta,(\nu_1,\nu_2))$
  is unitary if and only if  $L_{\phi(\mu)}(G(\delta)',\nu)$ is unitary.
\end{enumerate}
\end{proposition}
}

The two propositions \ref{r:m} are indicative of the general
  case. The calculation of nonspherical petite $K$-types in
  \cite{BP}, given here in section \ref{sec:3a}, implies the
 ``only if'' statement of the proposition for all nonspherical
  principal series of classical split simple groups.  For exceptional
  groups there is at least one minimal
  principal series for $G=F_4$, the one labeled $\delta_3$ in table
  \ref{table:fine} in section \ref{sec:fine1}, for which the unitary
  set is larger than the   spherical unitary dual for the
  corresponding $G(\delta_3).$

\subsection{}We give an outline of the paper. Section 2 presents basic facts and
examples about   minimal principal series of quasisplit real groups,
intertwining operators, fine types, and R-groups. Section 3 is
concerned with minimal principal series of split $p$-adic groups and
 affine Hecke algebras. Section 4 presents the idea and the
construction  of petite $K$-types, and the relation between
intertwining operators for minimal
 principal series of real split groups and intertwining operators for
 (extended) graded Hecke algebras. In  section 5, we record the main
 elements involved in the determination of the spherical unitary dual
 for split real and $p$-adic classical groups and for split $p$-adic
 exceptional groups. Section 6 has lists of parameters for unitary
 spherical representations of split groups.

\bigskip

Some of the results on which this exposition is based were presented
in June 2006, at the Snowbird conference on ``Representations of real
Lie groups'', in honor of B. Casselman and D. Mili\v ci\'c. We would
like to thank the organizers, particularly P. Trapa, for the
invitations to attend the conference and for their efforts in creating
a successful meeting.  { This research was supported by the
  NSF grants DMS-0300172 and FRG-0554278.

}

\section{Minimal principal series for real groups}\label{sec:2}

In this section, and in section \ref{sec:3a}, we will use the
classical notation from the theory of reductive Lie groups. For example,
if $G$ is a real group, the Lie algebra of
$G$ is denoted by $\fg_0$ and the complexification by $\fg,$ and
a minimal parabolic subgroup $P$ has the decomposition $P=MAN$.
The set of inequivalent irreducible representations of a group $H$ is
denoted by $\widehat H.$

\subsection{Minimal Principal Series}\label{sec:2.1} The definitions and
basic properties of minimal principal series can be found for
example in chapter 4 of \cite{green}.

Let $G$ denote a quasisplit real linear reductive group in the sense
of \cite{green}. The main case of
interest in this paper is when $G$ is the real points of a linear
connected reductive group {\it split} over $\bR.$
 Let $K$ be a maximal compact subgroup corresponding to
the Cartan involution $\theta$.

Let $P=MAN$ be a minimal parabolic
subgroup of $G$ (a Borel subgroup), and $H=MA$ be the Cartan
subgroup. The group $M$ is abelian, because $G$ is quasisplit. (If
$G$ is split, then $M$ is a finite abelian $2$-group.) The Iwasawa decomposition
is $G=KAN.$ One identifies $\fa^*$($\cong\check\fa$) with $\widehat
A$ via
\begin{equation}\label{d:2.1}
\nu\mapsto e^\nu,\ (e^\nu)(a)=\exp(\nu(\log a)),\ \nu\in\fa^*, a\in
A.
\end{equation}

 For every irreducible representation
$(\delta,V^\delta)$   of $M$, and every element $\nu$ of
$\check\fa$, we denote the {\it minimal
  principal series} by
\begin{equation}\label{eq:2.1.2}
X(\delta,\nu):=\Ind_{MAN}^G(\delta\otimes e^\nu\otimes 1).
\end{equation}
($\Ind$ denotes normalized induction.) The $K$-structure of
$X(\delta,\nu)$ is easy to describe:
\begin{enumerate}
\item The restriction of $X(\delta,\nu) $ to $K$ is \begin{equation}X(\delta,\nu)|_K=\Ind_M^K(\delta).\end{equation}
\item For every $K$-type $(\mu,E_\mu)$, one has Frobenius reciprocity
\begin{equation}\label{eq:2.1.3}
\Hom_K\left(\mu,X(\delta,\nu)\right)=\Hom_M\left(\mu|_M,\delta\right).
\end{equation}
\end{enumerate}

\begin{definition}\label{d:2.1.1} A representation $(\pi,V)$ of $G$ is called {\it
    spherical} if $\pi$ has fixed vectors under $K.$ The minimal
  principal series $X(\nu):=X(triv,\nu)$ is called the {\it spherical
    principal series}.
\end{definition}

%\smallskip

\subsection{}\label{sec:2.1a}
%{The \boldmath{$SL(2)$}-subgroup attached to a root}

Denote by
\begin{align}
&\Delta(\fg,\fh)        &\text{the set of roots of }\fh\text{ in }\fg,\\\notag
&\Delta^+(\fg,\fh)        &\text{the set of roots in
  }\Delta(\fg,\fh)\text{ whose root spaces lie in }\fn,\\\notag
&\Delta &\text{ the set of {\it restricted roots} of
  }\fa\text{ in }\fg, \\\notag
&\bar\Delta &\text{ the {\it reduced roots} in }\Delta,\\\notag
&\Delta^+\subset\Delta &\text{ the restricted roots whose  root
spaces lie in }\fn_0.
\end{align}
If $G$ is split, then $\bar\Delta=\Delta.$ The Weyl group of
$\Delta$ is
\begin{equation}\label{eq:2.1.1}
W=N_G(A)/A=N_K(A)/M.
\end{equation}

A root in $\Delta(\fg,\fh)$ is called {\it real} if $\theta\al=-\al.$
It is called complex otherwise (recall that $G$ is assumed
quasisplit). A root in $\Delta$ is called {real} if it is a
restriction of a real root in $\Delta(\fg,\fh)$), or otherwise it is
called {complex}. When $G$ is split, all roots in $\Delta$ are real.

 For every $\al\in\Delta$, define
\begin{align}\label{eq:2.1.4}
&\fa_0^\al=\text{the kernel of }\al\text{ on }\fa_0, \text{ with }\\\notag
  &A^\al=\text{the corresponding connected subgroup of }G,\\\notag
&M^\al A^\al=Z_G(A^\al), \text{ where the Lie algebra of }M^\al\text{
    is }\\\notag
&\fm_0^\al=\text{the algebra spanned by }\fm_0\text{ and
  }\{E_{c\al}\text{ root vectors}:~c\ge 1\},\\\notag
 &K^\al=K\cap M^\al.
\end{align}

When the root $\al\in\bar\Delta$ is real, this gives rise an
$SL(2)$-subgroup
\begin{equation}
G^\al=K^\al A^\al N^\al\cong SL(2,\bR)
\end{equation}
 as follows.
 Let
\begin{equation}\label{eq:2.1.5}
\phi_\al\colon sl(2,\bR)\to \fg_0
\end{equation}
be a root homomorphism such   that $H_\al=\phi \left(\begin{matrix}1 &0\\0
    &-1\end{matrix}\right) \in \fa_0,$ and $E_\al:=\phi \left( \begin{matrix}0 &1\\0
    &0\end{matrix}\right)$  is an $\al$-root vector. Then $G^\al$ is
    the connected subgroup of $G$ with Lie algebra $\phi(sl(2,\bR)).$

We will need the following notation later:
\begin{align}\label{eq:2.1.6}
&Z_\al:=E_\al+\theta(E_\al)\in \fk k_0^\al\\\notag
&\sigma_\al:=\exp\left(\frac\pi 2 Z_\al\right)\in K^\al \text{ is a
  representative of }s_\al\text{ in } N_K(A)\\\notag
&m_\al:=\sigma_\al^2\in M \text{ an element of order }2.
\end{align}
When $G$ is split semisimple, the elements $m_\al$ generate $M$.

\begin{comment}
 If $G$ is semisimple and has a
complexification $G^{\mathbb{C}}$, then the group $M$ is generated by
the identity component $M^0$ and by the elements $m_{\alpha}$'s for
$\alpha$ real. Also notice that, in this case,  $m_{\alpha} =
\exp(2\pi i \|\alpha\|^{-2}H_{\alpha}).$
\end{comment}

%\smallskip

\subsection{Fine \boldmath{$K$}-types}\label{sec:fine}
We will only  recall the properties of  {\it fine $M$-types} and
{\it fine $K$-types} that we need later (for a general definition see chapter 4 of
\cite{green}).

 The number of fine $M$-types is directly related to
the disconnectedness of $M.$ The two extreme cases are as follows.
\begin{enumerate}
\item If $G$ is a
split semisimple group, then any $M$-type is fine. In this case, a
$K$-type is fine if and only if
\begin{equation}\label{eq:2.1.7}
\text{the eigenvalues of }\mu(iZ_\al)\text{ are in }\{0,\pm
1\}\text{ for all (real) }\al.
\end{equation}
\item if $G$ is a complex semisimple group, then only the trivial
  $M$-type and the trivial $K$-type are fine.
\end{enumerate}
If $\mu$ is a fine $K$-type containing the fine $M$-type $\delta$,
then $\mu$ is {\it minimal} (in the sense of \cite{green}) in the
principal series $X(\delta,\nu).$

The number of fine $K$-types containing a given fine $M$-type is
related to the {\it R-group} of $\delta$, which we now recall.

\begin{definition}\label{d:2.1.2} Define $\Delta_\delta$ to be the
  subset of roots $\al$ in $\bar\Delta$ such that either $\al$ is complex,
  or, if $\al$ is real, then $\delta(m_\al)=1.$ One calls $\Delta_\delta$
  the {\it good roots with respect to $\delta$}. Let $^\vee
  \Delta=\{\check \al:\al\in\Delta\}$ be
  the set of {\it good coroots}.
\end{definition}
If $\delta$ is fine, then $\Delta_\delta$ is a root system. We
denote the corresponding Weyl group by
$W_\delta^0=W(\Delta_\delta)$. We also set
\begin{equation}
W_\delta=\{w\in W\colon  w\cdot\delta\simeq \delta\}
\end{equation}
(the stabilizer of $\delta$ in $W$).  Let
$\Delta_\delta^+=\Delta^+\cap\Delta_\delta,$ and similarly define
$^\vee\Delta_\delta^+.$ Then
\begin{enumerate}
\item $W_\delta^0$ is a normal subgroup of $W_\delta$, and
\[
W_{\delta} = W_{\delta}^{0} \rtimes R_{\delta}^{c}
\]
with  $ R_{\delta}^{c} = \{w \in W \colon
w({}^{\vee}\Delta_{\delta}^{+})= {}^{\vee}\Delta_{\delta}^{+}  \}. $
Note that $R_{\delta}^{c} $ is contained in the Weyl group generated
by the reflections through the roots perpendicular to
$\rho(\Delta_{\delta})$. Since these roots are strongly orthogonal,
$R_{\delta}^{c} $ is an abelian 2-group.

\item  The quotient
\begin{equation}R_\delta=W_\delta/W_\delta^0
\end{equation}
 is a finite abelian  $2$-group (isomorphic to $R_{\delta}^{c}$). We call $R_\delta $ the {\it  R-group} of $\delta$.
\item If $G$ is connected  semisimple (and has a complexification),
then we can identify $W_{\delta}$ with the stabilizer of the good
co-roots for $\delta$:
\[
W_{\delta} = \{w \in W \colon w({}^{\vee}\Delta_{\delta})=
{}^{\vee}\Delta_{\delta} \}.
\]

\end{enumerate}

\begin{theorem}[\cite{green}, 4.3.16]\label{p:2.2}
Let $\delta$ be a fine $M$-type. Set
\begin{equation}
A(\delta)=\{\mu\text{ fine }K\text{-type such that
}\Hom_M(\mu|_M,\delta)\neq \{0\}\}.
\end{equation}
Then
\begin{enumerate}
\item $A(\delta)\neq\emptyset.$
\item If $\mu\in A(\delta)$, then $\mu|_M=\displaystyle{\bigoplus_{\delta'\in W\cdot \delta} \delta'}.$ In
  particular, $\mu$ appears with multiplicity one in $X(\delta,\nu).$
\item There is a natural simply transitive action of $\widehat{R_\delta}$ on
  $A(\delta).$
\end{enumerate}
\end{theorem}
\noindent($\widehat{R_\delta}$ is a group, because $R_\delta$ is abelian.)
For the action of $\widehat{R_\delta}$ on $\widehat K$, and in
particular on $A(\delta)$, we refer the reader to 4.3.47 in
\cite{green}.
Notice that the cardinality of $R_\delta$  equals the number of fine
$K$-types containing  $\delta$ (as well as the number of fine
$K$-types contained in $X(\delta,\nu)$).

\medskip

\noindent{\bf Examples.} (a) If $G=SL(2,\bR),$ then $K=SO(2)$ and
$M=\bZ/2\bZ$. There are two $M$-types $triv$ and $sgn$, and they are
both fine. \\If $\delta$ is trivial, then $\Delta_{\delta}=\Delta$,
$R_{\delta}=\{1\}$ and $X(\delta,\nu)$ contains  a  unique
 fine $SO(2)$-type ($\mu = (0)$).   If $\delta$ is $sign$, then
  $\Delta_{\delta}=\emptyset$,  $R_{\delta}=\bZ/2\bZ$ and $X(\delta,\nu)$  contains two
fine  $SO(2)$-types,   namely $(+1)$ and $(-1)$.

(b) If $G=SL(2,\bC)$, then $K=SU(2),$ and
$M=\left\{\left(\begin{matrix}
    e^{i\theta} &0\\0
    &e^{-i\theta}\end{matrix}\right)\colon \theta\in\bR \right\}\cong U(1).$ So
$\widehat M\cong \bZ,$ but (by definition) only the trivial $M$-type
is fine.
If $\delta= triv$, then   $\Delta_{\delta}=\Delta$ and $ \#
R_{\delta}=\#A(\delta)=1$.

(c) If $G=SU(2,1),$ then $K=S(U(2)\times U(1)),$ and $M=U(1)$. Again
$\widehat M\cong \bZ,$ but (by definition) only the trivial $M$-type
is fine.    If $\delta= triv$, then   $\#R_\delta=\#A(\delta)=1$.

\subsection{Subquotients of Principal Series}\label{subquotients}
 Let $\delta$ be a fine
$M$-type and let $A(\delta)$ be the set of fine $K$-types in
$X(\delta,\nu)$ from theorem \ref{p:2.2}. Since every $\mu\in
A(\delta)$ appears in $X(\delta,\nu)$ with multiplicity one, there
is a unique  irreducible subquotient  of $X(\delta,\nu)$
 which contains $\mu$. We set
\begin{equation}
L(\delta,\nu)(\mu)=\text{ the unique irreducible subquotient
 of } X(\delta,\nu)
  \text{ containing }\mu.
\end{equation}
In general, $L(\delta,\nu)(\mu)$ may contain other fine $K$-types
(other than $\mu$).
\begin{definition}
Define
\begin{align}
&W(\nu)=\{w\in W\colon  w\nu=\nu\},\\\notag &W_\delta(\nu)=W(\nu)\cap
W_\delta,\quad
  W_\delta^0(\nu)=W_\delta^0\cap W(\nu),\\\notag
&R_\delta(\nu)=W_\delta(\nu)/W_\delta^0(\nu)\subset R_\delta,\\\notag
&R_\delta^\perp(\nu)=\{\chi \in \widehat{R_\delta}\colon \chi(r)=1,
\, \text{for all}\ r\in
  R_\delta(\nu)\}.
\end{align}

\end{definition}

\begin{theorem}[\cite{green}] \label{t:2.4} Assume $\mu\in A(\delta).$

\begin{enumerate}

\item Two subquotients $L(\delta,\nu)(\mu)$ and $L(\delta,\nu')(\mu)$
  are equivalent if and only if there exists $w\in W_\delta$ such that
  $\nu'=w\nu.$

\item  Every irreducible
  $(\fg,K)$-module containing the $K$-type $\mu$ is equivalent to a
  subquotient $L(\delta,\nu)(\mu)$ for some $\nu\in \fa^*.$

\item Let $\mu'\in A(\delta)$ be another fine $K$-type. Then
  $L(\delta,\nu)(\mu)=L(\delta,\nu)(\mu')$ if and only if $\mu$ and
  $\mu'$ are conjugate by an element of $R_\delta^\perp(\nu).$

\end{enumerate}

\end{theorem}

\begin{comment}
(These results were known before their appearance in \cite{green};
however, the proof presented in \cite{green} is different from the
original one.)
\end{comment}

\subsection{} A parameter $\nu\in\fa^*$ is called {\it weakly dominant} (respectively,
{\it strictly dominant}) if
\begin{equation}\label{eq:dom}
\langle\Re(\nu),\check\al\rangle\ge 0 \ (\text{respectively},
\langle\Re(\nu),\check\al\rangle>0), \text{ for all }\al\in\Delta^+.
\end{equation}

An important consequence of theorem \ref{t:2.4} is:

\begin{corollary}
If $\Re(\nu)$ is strictly dominant, then there exists a unique
irreducible  subquotient which contains all fine $K$-types in $
A(\delta)$.
\end{corollary}

On the other hand, when $\nu$ is purely imaginary ($\Re(\nu)=0$), every
fine $K$-type in $A(\delta)$ parameterizes one constituent of
$X(\delta,\nu)$ (a result obtained first by Knapp-Stein):

\begin{theorem}[Knapp-Stein]
If $\nu\in \fa^*$ is purely imaginary, then
\begin{equation}
X(\delta,\nu)=\bigoplus_{i=1}^{|R_\delta(\nu)|} X^i(\delta,\nu),
\end{equation}
where $X^i(\delta,\nu)$ are irreducible and inequivalent
subrepresentations.
 Moreover, each $X^i(\delta,\nu)$ contains some
fine $K$-type $\mu\in A(\delta).$
\end{theorem}

The trivial $K$-type is the only fine $K$-type in $A(triv)$, so
every spherical principal series $X(\nu)=X(triv,\nu)$ has a unique
irreducible spherical subquotient $L(\nu)=L(triv,\nu)(0)$. In particular,
the spherical principal series $X(\nu)$ is irreducible for all $\nu$
purely imaginary. The original proof of this result is due to Kostant.

\subsection{An example: $Sp(2n,\mathbb{R})$}\label{sec:2.5}
Let  $G$ be the split group $ Sp(2n,\,  \mathbb{R})$. Then  $K
\simeq U(n)$ and    $ M \simeq \mathbb{Z}_{2}^{n}$. Let
\begin{equation}
 \Delta^{+} = \{  \epsilon_i \pm \epsilon_{j\,} \colon
i,\,j= 1 \dots n,\, i < j \} \cup  \{  2 \epsilon_{l\,} \colon  l= 1
\dots n \}\end{equation} be the set of positive roots. Note that
every root is real. The Weyl group  $W= \mathcal{S}_n \rtimes
\left(\mathbb{Z}/2\mathbb{Z}\right)^n$ acts as the group of all
permutations
 and sign changes on  $\{\epsilon_1 ,\,\epsilon_2 \dots \epsilon_n
 \}$.

 The abelian
group $M$ has $2^n$ irreducible inequivalent  representations, all
of the form \begin{equation} \delta_{S}  \colon \text{diag} (
\lambda_1, \,\lambda_2, \dots , \,\lambda_n,\,\lambda_1,
\,\lambda_2, \dots , \,\lambda_n  ) \mapsto \prod\limits_{j \in S}
\lambda_j
\end{equation}
for some $S \subset \{1,\dotsb, n \}.$ Because $G$ is split and
semisimple, every $M$-type is fine. The Weyl group partitions
$\widehat{M} $ into $(n+1)$ conjugacy classes;  we
 choose representatives:
\begin{align}
\delta_{0} &=\delta_{\emptyset} \quad \text{and}\\
 \delta_p & = \delta_{\{n-p+1,n-p+2,\dots, n\}}\quad \text{for } 1\leq p\leq n.
\end{align}
We compute the number of irreducible subquotients of
$X(\delta_i,\nu)$, for all $i=0\dots p$.

  $\delta_{0} $ is the
trivial representation of $M$, so
$\Delta_{\delta_{\emptyset}}=\Delta$,
$W_{\delta_{\emptyset}}=W_{\delta_{\emptyset}}^0=W$ and
$R_{\delta_{\emptyset}}=\{1\}$. Since there is  a unique fine
$K$-type containing $\delta_0$, the principal series
$X(\delta_0,\nu)$ has   a unique subquotient
$L(\nu)=L(\delta_0,\nu)(0)$.

Now assume   $1\leq p\leq n$. To identify the  good roots for
$\delta_{p}$, we need to evaluate $\delta_p$ on the elements
$m_{\alpha}$'s. For every  (positive) root  $\alpha$, write
\[  m_{\alpha}= \exp\left( 2 \pi i \| \alpha \|^{-2}
H_{\alpha} \right)   = \text{diag}(d_1,\,d_2,\dots ,d_n,\,
d_1,\,d_2,\dots, d_n). \]
 Notice that the entries of $m_{\alpha}$
are 1, with the exception of $d_i=d_j=-1$  if $\alpha = \epsilon_i
\pm \epsilon_j$, and
 $d_k = -1$   if $\alpha =
2\epsilon_k$. Hence we find:
\begin{equation}
\delta_{p}(m_{\epsilon_i \pm \epsilon_j})= \begin{cases} +1, \text{
if either }  1 \leq i  < j \leq n-p \,
\text{ or }\,n-p+1  \leq i <  j  \leq n  \\
-1, \text{ otherwise}
\end{cases}
\end{equation}
and
\begin{equation}
\delta_{p}(m_{2\epsilon_{k}})= \begin{cases} +1, \text{ if } 1\leq k
\leq n-p
\\
-1, \text{ otherwise.}
\end{cases}
\end{equation}

The good roots for $\delta_p$ are
\begin{equation}
\Delta_{\delta_p} =
 \{ \pm \epsilon_i \pm \epsilon_j
\}_{1 \leq i <  j  \leq n-p}\cup
 \{ \pm 2\epsilon_k
\}_{1 \leq k \leq n-p}    \cup \{ \pm \epsilon_i \pm \epsilon_j
\}_{n-p+1 \leq i < j \leq n}.
\end{equation}
For brevity of notation, set $C_0$=$D_0$=$ D_1$=$ \emptyset$  and
$C_1 $=$ A_1$. Then
\begin{equation}
W_{\delta_p}^0= C_{n-p}\times D_p \quad \, \forall \, p=1\dots n.
\end{equation}

 Next, we compute the stabilizer of
$\delta_p$.  Recall that $W_{\delta_{p}}$  can be identified  with
the subgroup of Weyl group elements that preserve the good coroots
for $\delta_p$. Every permutation and sign change  on the sets
$\{\epsilon_1 \dots \epsilon_{n-p} \}$ and $\{\epsilon_{n-p+1} \dots
\epsilon_{n} \}$ has this property.
  Hence
\begin{equation}
W_{\delta_p}= C_{n-p}\times C_p \quad \, \text{ for all } p=1,\dotsb, n.
\end{equation}
Note that, for all $p=1\dots n$, $W_{\delta_{p}}^0$ is a normal
subgroup of $\,W_{\delta_{p}}$ of index 2. \\
Because  $R_{\delta_{p}}\simeq
 \mathbb{Z}/2\mathbb{Z}$,
 $\delta_{p}$ is contained into  two distinct fine $K$-types (notably $\Lambda^p(\mathbb{C}^n)$ and its dual), which we denote by $\mu_{\delta_p}^+$ and
 $\mu_{\delta_p}^-$.

Finally, we give the $R_{\delta_p}(\nu)$-group. Write
$\nu=(a_1,a_2,\dots ,a_n)$ with
\[
a_1 \geq a_2 \geq \dots \geq a_{n-p} \geq a_{n-p+1} \geq  \dots \geq
a_{n} \geq 0.
\]
Then
\begin{equation}
R_{\delta_p}(\nu) = \begin{cases} \{1\} & \text{ if $0 \not \in
\{a_{n-p+1},\dots ,a_n\} \Leftrightarrow a_n \neq 0$}
\\
\;\mathbb{Z}/2\mathbb{Z}& \text{ if $0  \in \{a_{n-p+1},\dots ,a_n\}
\Leftrightarrow a_n = 0$.}
\end{cases}
\end{equation}
We conclude that
\begin{itemize}

\item If the last entry of $\nu$ is  nonzero, then
 the principal series $X(\delta_p,\nu)$
has   a unique irreducible subquotient
$L(\delta,\nu)(\mu_{\delta_p}^+)=L(\delta,\nu)(\mu_{\delta_p}^-).$

\item If the last  entry of $\nu$  is
zero, then the principal series $X(\delta_p,\nu)$ has   two distinct
subquotients $L(\delta,\nu)(\mu_{\delta_p}^+)\neq
L(\delta,\nu)(\mu_{\delta_p}^-).$
\end{itemize}

\subsection{Lists of fine \boldmath{$K$}-types for split groups}\label{sec:fine1}
For the convenience of the reader, we record a list of examples of
fine $M$-types, fine $K$-types,  R-groups and sets of good roots for
split simple linear groups. For brevity, we do not include the trivial $M$-type,
and we only give  one fine $M$-type for each   orbit of $W$ in
$\widehat{M}$. We list, instead, all the fine $K$-types containing a
given fine $M$-type.

\begin{small}
\begin{center}
\begin{longtable}{|c|c|c|c|c|}
\caption{Table of fine types}\label{table:fine}\\
\hline \multicolumn{1}{|c|}{Group} & \multicolumn{1}{c|}{Fine
$M$-type $\delta$} &  \multicolumn{1}{c|}{Fine $K$-types}
&\multicolumn{1}{c|}{$R_{\delta}$-group} &
\multicolumn{1}{c|}{$\Delta_\delta$}
\\ \hline \hline
\endfirsthead

\begin{comment}
\multicolumn{3}{c}%
{{  \tablename\ \thetable{} -- continued from previous page}} \\
\hline \multicolumn{1}{|c|}{Group} & \multicolumn{1}{c|}{Fine
$M$-type $\delta$}&  \multicolumn{1}{c|}{Fine $K$-types} &
\multicolumn{1}{c|}{$R_{\delta}$-group} &
\multicolumn{1}{c|}{$\Delta_\delta$}
 \\ \hline  \hline
\endhead
\hline \hline
\endlastfoot
\end{comment}

$SL(2n+1,\bR)$ &$\delta_p,~1\le p\le n$ &    $\Lambda^p(\mathbb{C}^{2n+1})$     &$1$ &$A_{p-1}+A_{2n-p}$ \\
\hline \hline $SL(2n,\bR)$   &$\delta_{p},~1\le p<n$ &
$\Lambda^p(\mathbb{C}^{2n})$ &$1$
&$A_{p-1}+A_{2n-p-1}$\\
           \hline    &$\delta_{n}$            &  $\Lambda^n(\mathbb{C}^{2n})_{+}$     &$\bZ/2\bZ$ &$A_{n-1}+ A_{n-1}$\\
               &            &     $\Lambda^n(\mathbb{C}^{2n})_{-}$        &&\\
 \hline \hline
$Sp(2n,\bR)$ &$\delta_{p},~1\le p\le n$ &   $\Lambda^p(\mathbb{C}^{n})$,  $\Lambda^p(\mathbb{C}^{n})^{\ast}$        &$\bZ/2\bZ$ &$C_{n-p}+ D_p$  \\
\hline \hline
$SO(n+1,n)$ &$\delta_p,~1\le p\le n$ &$1\otimes\Lambda^p(\bC^n)$ &1 &$B_{n-p}+B_p$\\
\hline\hline
%$SO(n,n)$ &$\delta_p,~1\le p\le n$  &   $1
%\otimes
%\Lambda^p(\mathbb{C}^{n})$        &$\bZ/2\bZ$ &$D_{2n+1-p}+ D_p$     \\
% &&  $\Lambda^p(\mathbb{C}^{n})\otimes 1$         & &    \\
%\hline\hline
$SO(2n+2,2n+1)_0$   & $\delta_{p},~1\le p\le n$  &    $ 1\otimes \Lambda^p(\mathbb{C}^{2n+1})$      &$1$ &$B_{2n+1-p}+ B_p$  \\

\hline\hline

$SO(2n+1,2n)_0$   &$\delta_{p},~1\le p < n$  &   $ 1\otimes \Lambda^p(\mathbb{C}^{2n})$        &$1$ &$B_{2n-p}+ B_p$  \\

\hline
   &$\delta_{n}$  &   $ 1\otimes
\Lambda^n(\mathbb{C}^{2n})_{+}$     &$\bZ/2\bZ$
    &$B_{n}+B_{n}$  \\

      & &     $ 1\otimes \Lambda^n(\mathbb{C}^{2n})_{-}$     &
    &  \\

\hline\hline $SO(2n+1,2n+1)_0$    &$\delta_{p},~1\le p\le n$  &   $1
\otimes
\Lambda^p(\mathbb{C}^{2n+1})$        &$\bZ/2\bZ$ &$D_{2n+1-p}+ D_p$     \\

 &&  $\Lambda^p(\mathbb{C}^{2n+1})\otimes 1$         & &    \\

\hline \hline $SO(2n,2n)_0$         &$\delta_{p},~1\le p<n$ & $1
\otimes \Lambda^p(\mathbb{C}^{2n})$ &
$\bZ/2\bZ$ &$D_{2n-p}+D_{p}$\\

      & &  $\Lambda^p(\mathbb{C}^{2n})\otimes 1$ &
 &\\

   \hline            &$\delta_{n}$  &   $1 \otimes
\Lambda^p(\mathbb{C}^{2n})_{+}$      &      $\bZ/2\bZ\times \bZ/2\bZ$ &$D_{n}+ D_{n}$\\

              & &    $\Lambda^p(\mathbb{C}^{2n})_{+}\otimes 1$       &      & \\

 &  &   $1 \otimes
\Lambda^p(\mathbb{C}^{2n})_{-}$       &      &\\

 &  &   $\Lambda^p(\mathbb{C}^{2n})_{-}\otimes 1$       &      &\\

\hline \hline

$G_2$   &$\delta_3$ &$V_3\otimes\bC$      &$1$ &$A_1+\wti A_1$  \\
\hline \hline
$F_4$   &$\delta_3$  & $(2|0,0,0)$ &$1$ &$C_4$        \\
        &$\delta_{12}$ & $(1|1,0,0)$ &$1$ &$B_3+A_1$\\
\hline \hline
$E_6$   &$\delta_{27}$  & $\omega_2$  &$1$ &$D_5$        \\
        &$\delta_{36}$  & $2\omega_1$  &$1$ &$A_5+A_1$    \\
\hline \hline
$E_7$   &$\delta_{28}$ & $\omega_2$,\, $\omega_6$&$\bZ/2\bZ$ &$E_6$  \\
        &$\delta_{36}$ & $2\omega_1$,\, $2\omega_7$  &$\bZ/2\bZ$ &$A_7$  \\
        &$\delta_{63}$ & $\omega_1+\omega_7$ &$1$ &$D_6+A_1$\\
\hline \hline
$E_8$   &$\delta_{120}$ & $\omega_2$   &$1$ &$E_7+A_1$        \\
        &$\delta_{135}$ &  $2\omega_1$ &$1$ &$D_8$  \\
\hline

\end{longtable}
\end{center}
\end{small}
We explain the notation in the table.

\begin{enumerate}
\item If $G$ is a split  simple
linear group of classical type, the group $M$ consists of diagonal
matrices. We have denoted by $\delta_k$ the character of $M$ that
maps an element $m$ of $M$ into the product of the last $k$ diagonal
entries of $m$.
\item If $G$ is a split simple linear group of exceptional
type,  $\delta_k$ denotes the character of $M$ that has a
$k$-dimensional orbit under the action of the Weyl group.
\end{enumerate}

The group $K$ is  a   maximal compact subgroup of $G$,  as follows:

\[
\begin{array}{|c|c|} \hline G & K \\ \hline \hline
SL(m,\mathbb{R}) & SO(m) \\ \hline
Sp(2n,\mathbb{R}) & U(n)\\ \hline
SO(p,q) &S(O(p)\times O(q)) \\\hline
SO(p,q)_0 & SO(p)\times SO(q)\\

\hline
\end{array}
\quad
\begin{array}{|c|c|} \hline G & K \\ \hline
 \hline
G_2 & SU(2)\times SU(2)/\{\pm I\}\\
\hline F_4 & Sp(1) \times
Sp(3)/\{\pm I\}\\ \hline E_6 &  Sp(4)/\{\pm I\} \\
 \hline E_7 &
SU(8)/\{\pm I\} \\
\hline E_8 &  Spin(16)/\{I,w\}\\ \hline
\end{array}
\]

The notation for the maximal compact subgroup when $G=E_8$ means that
$K$ is a quotient of $Spin(16)$ by a central $\bZ/2\bZ,$ not equal to $SO(16).$  For
all $n$, we denote by $\mathbb{C}^n$   the standard representation
of $SO(n)$, $O(n)$, $SU(n)$  and $U(n)$.  We write $V_3$ for the
three-dimensional irreducible representation of $SU(2)$ (on the
space of homogeneous polynomials of degree 2 in 2 variables). The fine
$K$-types of $E_6$, $E_7$ and $E_8$ are described in terms of
fundamental weights; the ones of $F_4$ are described in terms of
standard (Bourbaki) coordinates.

\begin{remark}
For the purpose of computing  intertwining operators on principal
series, we need to fix \emph{a priori} a fine $K$-type
$\mu_{\delta}$ containing $\delta$. If the cardinality of
$A(\delta)$ is not one, the choice of $\mu_{\delta}$ is not
 canonical;
we choose:
\[
\begin{array}{|c|c|c|} \hline G & \delta
& \mu_{\delta} \\ \hline \hline SL(2n,\mathbb{R}) & \delta_n & \Lambda^n(\mathbb{C}^{2n})_+ \\
\hline SO(2n+1,2n)_0 & \delta_n & 1 \otimes \Lambda^n(\mathbb{C}^{2n})_+\\
\hline
SO(2n+1,2n+1)_0 & \delta_p \,(1\le p\le n) & 1 \otimes \Lambda^p(\mathbb{C}^{2n+1})\\
 \hline
SO(2n ,2n)_0 & \delta_p  \,(1\le p\le n) & 1 \otimes \Lambda^p(\mathbb{C}^{2n}) \\
\hline
SO(2n ,2n)_0 & \delta_n   & 1 \otimes \Lambda^n(\mathbb{C}^{2n})_+  \\
\hline
E_7 & \delta_{28}   & w_2 \\
\hline
E_7 & \delta_{36}   & 2w_1 \\
\hline
\end{array}
\]
This choice of $\mu_{\delta}$ will be
used to define  the full intertwining operator in section
\ref{sec:2.2},  the  $W_{\delta}$-representation in  section
\ref{sec 4.2}, and the bijection between $\widehat{R_{\delta}}$ and
$A(\delta)$ in sections \ref{sec:2.9} and \ref{sec 4.2}.
\end{remark}

\subsection{Intertwining operators}\label{sec:2.2} We recall some basic facts about
intertwining operators for minimal principal series of real groups (see
chapter 7 of \cite{Kn} for more details).

Fix an $M$-type $\delta$ and a character $\nu$ of $A$, and let
$X(\delta,\nu)$ be the minimal principal series induced from
$P=MAN$. Denote by $\bar{P}=MA\bar{N}$  the
 opposite  parabolic. For every $w\in W,$ one
defines a {\it formal intertwining operator}
\begin{align}\label{eq:2.2.1}
&A(w,\delta,\nu) \colon  X(\delta,\nu)\longrightarrow
X(w\delta,w\nu),\\\notag &(A(w,\delta,\nu)F)(g)=\int_{\bar{N}\cap
(wN)} F(gw\bar n)~d\bar
  n,\quad g\in G.
\end{align}
Formally, the integral defines an intertwining operator, but it may not
converge for general $\nu.$

\subsubsection{}
The intertwining operator $A(w,\delta,\nu)$ can be decomposed as
follows. (This is the {\it
Gindikin-Karpelevi\v c  decomposition}.)

\begin{proposition} Assume that $w=w_1w_2$, with
  $\ell(w)=\ell(w_1)\ell(w_2).$ Then
\begin{equation}\label{eq:2.2.3}
A(w,\delta,\nu)=A(w_1,w_2\delta,w_2\nu)\circ A(w_2,\delta,\nu).
\end{equation}
\end{proposition}
\begin{corollary} If  $w=s_1\cdot s_2\cdot\dots\cdot s_m$ is a {\it minimal}
decomposition of $w$ as a product of simple reflections, then
$A(w,\delta,\nu)$ factors as
\begin{equation}\label{eq:2.2.4}
A(w,\delta,\nu)=A(s_1,w_1\delta,w_1\nu)\cdot
A(s_2,w_2\delta,w_2\nu)\cdot\dots\cdot A(s_m,w_m\delta,w_m\nu),
\end{equation}
where $w_k=s_{k+1}\dots s_m,$ for all $1\le k\le m.$
\end{corollary}

If
\begin{equation}\label{eq:2.2.5}
\langle \Re(\nu),\beta\rangle>0, \text{ for every root }\beta\in
\Delta^+\text{ such that }w\beta\notin\Delta^+,
\end{equation}
then the integral in (\ref{eq:2.2.1}) is actually convergent. This
is proved  using the decomposition (\ref{eq:2.2.4}) and an
investigation of the rank one cases (we will make this more precise
below).

\subsubsection{}\label{sec:2.7.2}

For every $K$-type $(\mu,E_\mu)$, the intertwining operator
$A(w,\delta,\nu)$ induces an operator
\begin{equation}\label{eq:2.2.6}
A_\mu(w,\delta,\nu) \colon \Hom_K\left(
\mu,X(\delta,\nu)\right)\longrightarrow
\Hom_K\left(\mu,X(w\delta,w\nu)\right).
\end{equation}
By Frobenius reciprocity  (\ref{eq:2.1.3}), this can
be regarded as an operator
\begin{equation}\label{eq:2.2.7}
A_\mu(w,\delta,\nu) \colon \Hom_M(\mu,\delta)\longrightarrow
\Hom_M(\mu,w\delta).
\end{equation}
Via (\ref{eq:2.2.3}), $A_\mu(w,\delta,\nu)$ also acquires a
decomposition into factors corresponding to simple reflections. Let
$A_\mu(s_\al,\delta,\nu)$ be such a simple reflection factor.

\begin{remark}\label{r:2.2} The operator $A_\mu(s_\al,\delta,\nu)$ for
  $G$ agrees with the operator
  $A_{\mu|_{MK^\al}}(s_\al,\delta,\nu|_{\fa^\al})$ for the real rank
  one group $MG^\al.$
\end{remark}

More precisely, if $G$ is split,
\begin{equation}
\mu=\bigoplus_{m \in \mathbb{Z}}  \mu_m^\al
\end{equation}
is the decomposition of the $K$-type $\mu$ into isotypic components
of $K^\al=SO(2)$ ($\widehat{K^\al}=\bZ$), then the decomposition
\begin{equation}
\Hom_M(\mu,\delta)=\bigoplus_{m \in \mathbb{N}}
\Hom_M(\mu_m^\al+\mu_{-m}^\al,\delta)
\end{equation}
is preserved by $A_\mu(s_\al,\delta,\nu).$  Moreover, the
restriction of $A_\mu(s_\al,\delta,\nu)$ to
$\Hom_M(\mu_m^\al+\mu_{-m}^\al)$ coincides with the operator
$A_{\mu_m^\al+\mu_{-m}^\al}(s_\al,\delta,\nu|_{\fa^\al})$ for
$MG^\al$.

\begin{comment}
(A   list of the  rank one root subgroups that appear can be found
in Appendix C of \cite{Kn2}. The operators for $SL(2,\mathbb{R})$
and $SL(2,\mathbb{C})$ are recorded in section \ref{rank one}.)
\end{comment}

\subsection{Reducibility}\label{sec:2.9} We discuss the reducibility
of the principal series $X(\delta,\nu)$.

Assume that $\Re(\nu)$ is weakly dominant with respect to
$\Delta^+.$ The first instance in which $X(\delta,\nu)$ becomes
reducible is when the R-group $R_\delta(\nu)$ is nontrivial.

 Partition the restricted roots according to their
inner product with $\Re(\nu)$: \[\Delta= \Delta_L \sqcup \Delta_U^+
\sqcup \Delta_U^-\]  with
\begin{itemize}
\item[]$\Delta_L= \{ \alpha \in \Delta \colon \langle
\Re(\nu),\,\alpha \rangle = 0 \}$,

\item[]$\Delta_U^+= \{ \alpha \in \Delta \colon \langle
\Re(\nu),\,\alpha \rangle > 0 \}$,

 \item[]$\Delta_U^-= \{ \alpha \in \Delta \colon \langle
\Re(\nu),\,\alpha \rangle < 0 \}.$
\end{itemize}
The set $\Delta_L$ is a root system, and $\Delta^+=\Delta_L^+\sqcup
\Delta_U^+.$

Denote by $L$ the centralizer of $\Re(\nu)$ in $G$. This is a Levi
subgroup containing the Cartan $MA$, with Lie algebra $\fk l_0=\fk
m_0\oplus\fa_0\oplus\left(\bigoplus_{\al\in\Delta_L}\fg_\al\right).$
Define the parabolic subgroup $Q=LU$ of $G$ with Lie algebra $\fk
q_0=\fm_0\oplus\fa_0\oplus \left(
\bigoplus_{\al\in\Delta_L}\fg_\al\right)\oplus\left(\bigoplus_{\al\in\Delta_U^+}\fg_\al\right).$

The following lemma is {\it induction by stages}.

\begin{lemma}
$X(\delta,\nu)=\Ind_P^G(\delta\otimes\nu)=\Ind_{LU}^G\left(\Ind_P^L(\delta\otimes
  \nu)\otimes 1\right).$
\end{lemma}

Note that $\nu$ is imaginary for $L$, so in parenthesis we have
unitary induction. By the results presented in section \ref{subquotients}, the
unitarily induced module  $\Ind_P^L(\delta\otimes
  \nu)$ decomposes as the direct sum of $\# R_{\delta}(\nu)$
irreducible inequivalent representations of $L$.
 Let
\begin{equation}\label{eq:2.9.1}
X(\delta,\nu)=\bigoplus_{r\in R_\delta(\nu)} X^r(\delta,\nu)
\end{equation}
be the corresponding decomposition of $X(\delta,\nu)$.

  Recall that
$\widehat{R_\delta}$ acts on the set of fine $K$-types $A(\delta)$
  simply transitively (theorem \ref{p:2.2}). Implicit in
  (\ref{eq:2.9.1}) is the fact that we fixed a particular fine
  $K$-type $\mu_\delta$. This gives a bijection between $R_\delta$ and
  $A(\delta)$: $r\in R_\delta\mapsto \mu_{\delta,r},$ where
  $\mu_{\delta,0}=\mu_\delta.$

The fine $K$-types occurring in $X^r(\delta,\nu)$ form an orbit
under the action of $R_\delta^\perp(\nu)$ on $A_\delta(\nu).$

Let $\mu$ be a $K$-type which occurs in $X(\delta,\nu)$. The space
\begin{equation}
\Hom_K\left(\mu,X(\delta,\nu)\right)=\Hom_M(\mu,\delta)
\end{equation}
carries a representation of $W_\delta$. Identify
$\widehat{R_\delta}$ with $R_\delta^c$, so  that
$W_\delta=W_\delta^0\rtimes \widehat{R_\delta}$, and regard
$\Hom_K\left(\mu, X(\delta,\nu)\right)$ as a
$\widehat{R_\delta}$-module (by restriction).\\
Because $R_{\delta}(\nu) \subset R_{\delta}$ and $R_{\delta}$ is
abelian, we can regard $\widehat{R_{\delta}(\nu)}$ as a  subset of
$\widehat{R_{\delta}}$, and restrict $\Hom_K\left(\mu,
X(\delta,\nu)\right)$ to
 $\widehat{R_{\delta}(\nu)}$.
  The  $\widehat{R_{\delta}(\nu)}$-module
structure on $\Hom_K\left(\mu, X(\delta,\nu)\right)$ is compatible
with the action of $R_\delta(\nu)^\perp$ on $A(\delta).$

\begin{proposition}\label{p:2.9} Let $r$ be an element of $R_\delta(\nu),$ and
let
  $\mu$ be a $K$-type containing $\delta$. Identify  $R_\delta(\nu)$ with its double dual.
  Then:
\begin{equation}
\Hom_K\left(\mu,
X^r(\delta,\nu)\right)=\Hom_{\widehat{R_{\delta}(\nu)}}\left(r,\Hom_K\left(\mu,
X(\delta,\nu)\right)\right).
\end{equation}
\end{proposition}

\subsection{} A second way in which $X(\delta,\nu)$ becomes reducible is when the
operator $A(w_0,\delta,\nu)$ has a nontrivial kernel. Let
$w_0=s_1\cdot s_2\cdot\dots\cdot s_m$ be a minimal decomposition of
$w_0$,   and let
\[
A(w_0,\delta,\nu)= \prod_{i=1}^m A(s_i,w_i\delta,w_i\nu)
\]
be the  factorization of
$A(w_0,\delta,\nu)$ (as in section \ref{sec:2.2}). Then
$A(w_0,\delta,\nu) $ has  a nontrivial kernel if and only if one of
its factors has.

\begin{proposition}[\cf \cite{green}, 4.2.25] The long intertwining operator $A(w_0,\delta,\nu)$ has
  a nontrivial kernel if and only if there exists a simple root
  $\al\in\Delta$ such that $Re\langle\check\al,\nu\rangle\neq 0$ and
  $\Ind_{P}^{M^\al A^\al}(\delta,\nu)$ is reducible.
\end{proposition}

This proposition reduces the question of finding the kernel of the
long intertwining operator to a similar question for rank one
groups, where the answer is known.
 If $G$ is split, the only rank
one root subgroup that appears is $SL(2,\bR)$, and the answer is
particularly simple. (The operators for $SL(2,\mathbb{R})$ are given
in section \ref{rank one}.)

\begin{corollary} Assume that $G$ is split. Then the operator
  $A(w_0,\delta,\nu)$ has a kernel if and only if there exists a
 simple (real) root $\al\in\Delta$ such that
\begin{align}
&\langle\check\al,\nu\rangle=k,\text{ for some integer }k\neq 0, \text{
  and}\\\notag
&\delta(m_\al)=(-1)^{k+1}.
\end{align}
The parity condition means that $\langle\check\al,\nu\rangle$ should
be an odd integer if $\al\in\Delta_\delta$, and an even integer
otherwise.
\end{corollary}

 For split groups,  reducibility can only occur in the  two instances
described above. We summarize this in the next statement.

\begin{theorem}
Let $X(\delta,\nu)$ be a minimal principal series for a real split
group. Then $X(\delta,\nu)$ is reducible if and only if
\begin{itemize}
\item[$(i)$] the $R$-group $R_{\delta}(\nu)$ is non-trivial and/or

\item[$(ii)$] there is a simple root  $\alpha$ such that the inner product
$\langle \check{\alpha},\nu \rangle $ is a non-zero integer $k$, and
\[
\delta(m_{\alpha})= (-1)^{k+1} .
\]
\end{itemize}
\end{theorem}

\noindent{\bf Example.} If $G=SL(2,\bR),$ then
\begin{equation}
\text{$X(triv,\nu)$ is
reducible at  $\langle\check\al,\nu\rangle\in\{\pm 1,\pm 3,\pm
5,\dots\},$}
\end{equation}
 while
\begin{equation}
\text{$X(sgn,\nu)$ is reducible at
$\langle\check\al,\nu\rangle\in\{0,\pm 2,\pm 4,\dots\}.$}
\end{equation}
 The point of reducibility $\langle\check\al,\nu\rangle=0$ for
 $X(sgn,\nu)$ comes from the R-group: $R_{sgn}(0)=\bZ/2\bZ.$

\subsection{Langlands quotient}\label{sec:2.10}
We look at the connection between intertwining operators and
Langlands classification. Let $w_0\in W$ be the long Weyl group
element and
\begin{equation}\label{eq:2.4.1}
A(w_0,\delta,\nu)\colon  X(\delta,\nu)\longrightarrow
X(w_0\cdot\delta,w_0\cdot\nu)
\end{equation}
be the long intertwining operator. We choose a fine $K$-type
$\mu_\delta\in A(\delta).$ Since $\mu_\delta$ has multiplicity one
in $X(\delta,\nu)$, the operator $A_{\mu_\delta}(w_0,\delta,\nu)$ is
a scalar function of $\nu.$ We normalize $A(w_0,\delta,\nu)$ by this scalar, and
denote the resulting operators by
\begin{equation}\label{eq:2.4.2}
A'(w_0,\delta,\nu)\text{ and } A'_{\mu}(w_0,\delta,\nu).
\end{equation}
So $A'_{\mu_\delta}(w_0,\delta,\nu)=1.$

\begin{theorem}[Langlands, Mili\v ci\' c]\label{t:2.4.2} Assume $\Re(\nu)$
  is weakly dominant with respect to $\Delta^+.$

\begin{enumerate}
\item The operator $A'(w_0,\delta,\nu)$ has no poles.
\item The (closure of the) image of $A'(w_0,\delta,\nu)$ is the Langlands quotient
  $L(\delta,\nu).$ This is the unique largest completely reducible
  quotient of $X(\delta,\nu).$ When $\Re(\nu)$ is strictly dominant,
  $L(\delta,\nu)$ is irreducible.
\end{enumerate}
\end{theorem}

 The connection between this
classification and fine $K$-types can be formulated as follows.

\begin{corollary} Assume $\Re(\nu)$ is weakly dominant with respect to
  the roots in $\Delta^+.$ Every irreducible summand $L^i(\delta,\nu)$ of the
  Langlands quotient $L(\delta,\nu)$ is of the form
  $L(\delta,\nu)(\mu)$ for some fine $K$-type $\mu\in A(\delta).$

If $\Re(\nu)$ is strictly dominant, then $L(\delta,\nu)|_K$ contains
all the
  $K$-types in $A(\delta)$
  (each with multiplicity one).
\end{corollary}

\subsection{Hermitian forms}\label{sec:Hermitian} Assume for this
 section that $\Re(\nu)$ is strictly dominant with respect to the roots in
$\Delta^+.$ We  use the results in section \ref{t:2.4} to define
Hermitian forms and investigate the unitarity of the (irreducible)
Langlands quotient $L(\delta,\nu).$ The following result is based on
the fact that the Hermitian dual of $L(\delta,\nu)$ is
$L(\delta,-\overline\nu).$

\begin{theorem}[Knapp-Zuckerman]\label{t:2.5} Let $X(\delta,\nu)$ be a minimal principal series induced from $P$.
Assume that  $\Re(\nu)$ is strictly dominant with respect
  to $\Delta^+$ and let   $L(\delta,\nu)$  be the irreducible Langlands quotient
 of $X(\delta,\nu)$. Then
 \begin{enumerate}
 \item
  $L(\delta,\nu)$ admits a nondegenerate invariant Hermitian form if and only if
\begin{equation}\label{eq:Hermitian}
w_0\cdot\delta\cong\delta\text{ and
}w_0\cdot\nu=-\bar{\nu}.
\end{equation}
%($\bar{P}$ is the opposite parabolic, so $w$ must be the long Weyl
%group element $w_0$).
\item If   $L(\delta,\nu)$ is Hermitian, choose an isomorphism $\tau \colon
w_0\cdot\delta\to\delta$, and construct the operator
\begin{equation}\label{eq:2.5.1}
\CA(w_0,\delta,\nu)=\tau\circ A'(w_0,\delta,\nu).
\end{equation}
Let $\langle\ ,\ \rangle_h$ denote the pairing between the module
$L(\delta,\nu)$ and its Hermitian dual. Then the Hermitian form on
$L(\delta,\nu)$ is given by
\begin{equation}\label{eq:2.5.2}
\langle x,y \rangle:=\langle x, \CA(w_0,\delta,\nu) y\rangle_h \quad
\text{for all } x,y\in L(\delta,\nu).
\end{equation}
\end{enumerate}(Note that there is an implicit choice of the normalization,
dictated by  $A'(w_0,\delta,\nu)$.)
\end{theorem}

Fix a  fine $K$-type $\mu_\delta$ containing $\delta$.  We can
identify the representation space $V^\delta$ of $\delta$ with the
isotypic component of $\delta$ in $\mu_{\delta}$, and  define
\begin{equation}
\tau=\mu_\delta(w_0).
\end{equation}

\begin{remark} If $\Re(\nu)$ is dominant and  the  Langlands quotient $L(\delta,\nu)$  is
Hermitian, then   $L(\delta,\nu)$ is unitary (or unitarizable) if
and only if the invariant Hermitian form in (\ref{eq:2.5.2})  is
positive definite.

Because $L(\delta,\nu)$ is the quotient of $X(\delta,\nu)$ modulo
the Kernel of the operator $A'(w,\delta,\nu)$, the form on
$L(\delta,\nu)$ is positive definite if and only if the operator
 $A'(w,\delta,\nu)$  is positive semidefinite on $X(\delta,\nu)$.\end{remark}

  Restricting the
 attention to the various $K$-types in $X(\delta,\nu)$,  one
obtains the following criterion of unitarity.

\begin{corollary}\label{c:2.5}
 If $\Re(\nu)$ is dominant and  the  Langlands quotient $L(\delta,\nu)$  is
Hermitian, then  $L(\delta,\nu)$ is unitary if and only if the
Hermitian operators
\begin{equation}
\CA_{\mu}(w,\delta,\nu)
\end{equation}
are positive semidefinite for all $K$-types $\mu$ appearing in
$X(\delta,\nu).$
\end{corollary}

\begin{remark} If $\Re(\nu)$ is weakly dominant, the Hermitian
condition in (\ref{eq:Hermitian}) can be replaced by:
\begin{equation}
w\cdot Q = \bar{Q},\; w\cdot\delta\cong\delta\text{ and
}w\cdot\nu=-\bar{\nu},
\end{equation}
where $Q$ is the parabolic subgroup of $G$ defined by $\nu$, as in
section \ref{sec:2.9}. Note that the Langlands subquotient
$L(\delta,\nu)$ might be reducible. Equation (\ref{eq:2.5.2})
defines a non degenerate invariant Hermitian form on every
irreducible constituent of $L(\delta,\nu)$, \ie on every irreducible
subquotient of $X(\delta,\nu)$. Note that the form on
$L(\delta,\nu)(\pi)$ is normalized so that it takes the value
$\varrho_{\pi}[w]=\pm 1$ on the fine $K$-type $\pi$ (see section
\ref{sec 4.2}). Every operator $\CA_\mu(w,\delta,\nu)$ has a block
corresponding to  $L(\delta,\nu)(\pi)$. If $\varrho_{\pi}[w]=1 $
(respectively  $\varrho_{\pi}[w]=-1 $), the subquotient
$L(\delta,\nu)(\pi)$ is unitary if and only if this block is
positive semidefinite (respectively negative semidefinite) for all
$K$-types $\mu$.
\end{remark}

\begin{comment}
As stated, the unitarity criterion requires an infinite process,
since there are infinitely many $K$-types in a principal series.
D. Vogan has given a theoretical finite bound on the
$K$-types that one needs to consider, so the determination of
unitarity is a finite process (for each given group $G$).
\end{comment}
 For
practical purposes  the unitarity criterion stated above is used mainly to prove
that a given Hermitian $L(\delta,\nu)$ is {\it not} unitary,
unitarity being proven by other methods.

The idea of petite $K$-types is to give a  small set of
$K$-types on which the operator $\CA_\mu(w,\delta,\nu)$ is
computable, and these computations are sufficient for ruling out all
nonunitary $L(\delta,\nu).$ Before looking at the general case in
section \ref{sec:3a}, we
will present some examples of real rank one cases, where this idea will already
become apparent.

\subsection{Real rank one}\label{rank one}
In view of theorem \ref{t:2.4} and remark \ref{r:2.2}, it is of
particular importance to have the formulas for the intertwining
operators for real rank one groups. These are known (see \cite{JW} and
the references therein). We give three examples. Assume that
$\nu$ is real and weakly  dominant for the spherical cases below, and
strictly dominant for the nonspherical one.

\noindent\begin{enumerate}
\item $G=SL(2,\bR)$.

\begin{enumerate}
\item $\delta=triv$, normalized on $\mu_0$:
\begin{equation}\label{eq:2.3.7}
\CA_{\mu_{2k}}(s_\al,triv,\nu)=\prod_{j=1}^{|k|}\frac{(2j-1)-\langle\nu,\check\al\rangle}{(2j-1)+\langle\nu,\check\al\rangle};
\end{equation}

\item $\delta=sgn$, normalized on $\mu_1$:
\begin{equation}\label{eq:2.3.8}
\CA_{\mu_{2k+1}}(s_\al,sgn,\nu)=sgn(k)\prod_{j=1}^{|k+\frac
12|-\frac
12}\frac{2j-\langle\nu,\check\al\rangle}{2j+\langle\nu,\check\al\rangle}.
\end{equation}
\end{enumerate}
For all $k\in \bZ$, $\mu_{k} $   denotes the character
$e^{i\theta}\mapsto e^{ki\theta}$ of $K=SO(2)\cong U(1)$.

\item $G=SL(2,\bC)$ and $\delta=triv$, normalized on $\mu_0$:
\begin{equation}\label{eq:2.3.9}
\CA_{\mu_{k}}(s_\al,triv,\nu)=\prod_{j=1}^{k}\frac{j-\langle\nu,\check\al\rangle}{j+\langle\nu,\check\al\rangle}.
\end{equation}
For all $k\in \mathbb{N}$, $\mu_{k} $   denotes  the
$(k+1)$-dimensional representation of $SU(2).$
\end{enumerate}

\noindent From corollary \ref{c:2.5} and the calculations above, one concludes
that:

\begin{enumerate}
\item[$\bullet$] If $G=SL(2,\bR)$ and $\delta=triv$, the Langlands
  quotient $L(triv,\nu)$ is unitary if and only if
  $0\le\langle\check\al,\nu\rangle\le 1.$
   The
 spherical petite $K$-types for $G=SL(2,\bR)$ are $\mu_0$ and $\mu_2$. They have the property that $L(triv,\nu)$ is
  unitary if and only if the operators $\CA_{\mu}(s_\al,triv,\nu)$ are
  positive semidefinite for $\mu\in\{\mu_0,\mu_2\}.$
\item[$\bullet$] If $G=SL(2,\bR)$ and $\delta=sgn$, the Langlands
  quotient $L(sgn,\nu)$ is never unitary for $\nu>0.$  The
 nonspherical petite $K$-types for $G=SL(2,\bR)$ are $\mu_1$ and $\mu_{-1}$. They have the property that $L(sgn,\nu)$ is
  unitary if and only if the operators $\CA_{\mu}(s_\al,sgn,\nu)$ are
  positive semidefinite for $\mu\in\{\mu_{1},\mu_{-1}\}.$
\item[$\bullet$] If $G=SL(2,\bC)$ and $\delta=triv$, the Langlands
  quotient $L(triv,\nu)$ is unitary if and only if
  $0\le\langle\check\al,\nu\rangle\le 1.$  The
 spherical petite $K$-types for $G=SL(2,\bC)$ are $\mu_0$ and $\mu_1$. They have the property that  $L(triv,\nu)$ is
  unitary if and only if the operators $\CA_{\mu}(s_\al,triv,\nu)$ are
  positive semidefinite for $\mu\in\{\mu_0,\mu_1\}.$
\end{enumerate}

\section{Graded Hecke algebra and $p$-adic groups}\label{sec:3}

\subsection{Split $p$-adic groups}\label{sec:3.0} As in the
introduction, denote by
\begin{align}
&\bF&\text{ a $p$-adic field of characteristic }0,\\\notag
&\bO&\text{ the ring of integers in }\bF,\\\notag &\CP&\text{ the
unique prime ideal in }\bO,\\\notag
&\bF_q&\text{ the residue field with $q$ elements.}
\end{align}
Let $G$ denote the $\bF$-points of a connected reductive linear
algebraic group with root datum $(\C X,R,\C Y,\check R)$, split over $\bF$, and let $B$ be the Borel
subgroup, $H\cong (\bF^\times)^{\rank G}$ the Cartan subgroup,  $K=\bG(\bO)$ the maximal compact
subgroup, $^0\!H=H\cap K\cong (\bO^\times)^{\rank G}$ the compact part
of $H$ as in section \ref{sec:1.1}. We also fix
an {Iwahori subgroup} $\C I\subset K$  of $G$ as defined in section
\ref{sec:1.3}.

Let $\chi\colon H\to \bC^\times$ be a character.  Recall that $\chi$ is
{unramified} if $\chi|_{^0\!H}=1,$ and otherwise $\chi$ is called
{ramified}.   The minimal principal series $X(\chi)$ is defined
by induction, similarly to the one for real groups. It has a finite
composition series and, if $\chi$ is unramified, it contains a
unique irreducible ($K$-)spherical subquotient.  Every irreducible
spherical module appears in this way. However, unlike the case of
real groups, not every irreducible $G$-representation can be
realized as a subquotient of a minimal (ramified or unramified)
principal series.

A remarkable feature of the representation theory of $p$-adic groups
is the that it is often controlled by affine Hecke algebras.
This approach originated with the seminal work of Iwahori-Matsumoto (\cite{IM}),
Borel (\cite{Bo}),  and Casselman (\cite{Cas}).

\begin{definition}
The {\it Iwahori-Hecke algebra} $\CH(G//\CI)$ is the algebra of $\CI$-biinvariant complex
locally constant functions on $G$ with the convolution.
Define the category $\C C(\CI,triv)$ of (admissible) representations of
$G$ which are generated by their Iwahori fixed vectors.
\end{definition}
Iwahori-Matsumoto have determined the structure of $\CH(G//\CI)$ by
a set of generators in correspondence with the affine Weyl groups. We will instead the set of generators and relations
introduced first by Bernstein. The Iwahori-Hecke algebra is a
specialization of the {\it affine Hecke algebra} $\CH$ which we define now.
Let $z$ be an indeterminate (which can then be specialized to
$q^{1/2}$).  Then $\CH$ is an  algebra over $\bC[z,z^{-1}]$
generated
by $\{ T_w\}_{w\in W}$ and $\{\theta_x\}_{x\in
\C Y}$, subject to the relations
\begin{equation}\begin{aligned}
&T_wT_{w'}=T_{ww'} (l(w)+l(w')=l(ww')),\\
&\theta_x\theta_y=\theta_{x+y},\\
&T_s^2=(z^2-1)T_s+z^2,\\
&\theta_xT_s=T_s\theta_{sx}+(z^2-1)\frac{\theta_x-\theta_{sx}}{1-
\theta_{\alpha}},
\end{aligned} \label{2.1.2}\end{equation}
where $\al$ denotes a simple root, and $s$ is the corresponding simple
reflection.

\begin{theorem}\label{t:3.0}
\begin{enumerate}
\item(Borel) There exists an equivalence of categories
\begin{equation}
\C C(\CI,triv)\longrightarrow \CH(G//\CI)\text{-}modules, \quad V\mapsto
V^\CI.
\end{equation}
\item (Casselman) Every irreducible subquotient of the unramified
  principal series $X(\chi)$ is in $\C C(\CI,triv)$, and every irreducible object in
  $\C C(\CI,triv)$ has this form.
\end{enumerate}
\end{theorem}

As explained in section \ref{sec:1.3}, this approach has proven successful towards the determination of the
unitary dual as well. By \cite{BM1,BM2}, one reduces the
determination of the unitary representations in $\C C(\CI,1)$ to that
of the unitary dual of $\CH(G//\CI),$ and furthermore to the similar
problem for the graded Hecke algebra (section \ref{sec:3.1}) introduced in \cite{L1}.

\subsection{}\label{sec:3.0a}
The Borel-Casselman correspondence was vastly generalized in
\cite{HM}, and in the theory of types (\cite{BK} and the references therein). Hecke algebra
isomorphisms as in part theorem \ref{t:3.0}.(1) are constructed for many categories of representations of
$p$-adic groups. The Iwahori subgroup $\CI$ is replaced by an
arbitrary compact open subgroup $J$, and the trivial character $triv$
of $\C I$ is replaced by a character $\rho$ of $J.$ Then one defines
the category $\C C(J,\rho)$ of representations of $G$ which contain
$\rho$ in their restriction to $J$, and are generated by the
$\rho$-isotypic component. The Iwahori-Hecke algebra is replaced by
the algebra $\CH(G//J,\rho)$ of (locally constant) complex functions
$f$ on $G$, such that $f(j_1gj_2)=\rho(j_1)f(g)\rho(j_2),$ for all
$g\in G,j_1,j_2\in J.$

\medskip

In relation to section \ref{sec:2}, of particular importance to us
will be the case of ramified minimal principal series. The following
definitions should be compared with those for real groups in section
\ref{sec:2}.

Let $\Delta$ denote the root system determined by $\bH$ and $\bG$,
$\Delta^+$ the positive roots given by $\bB$. Let $W$ be the
(finite) Weyl group of $\Delta.$ Fix a character $\chi$ and let
$^0\!\chi=\chi|_{^0\!H} \colon  (\bO^\times)^{\rank G}\to\bC^\times$ be
its restriction to $^0\!H$. Define:
\begin{align}
&W_{^0\!\chi}=\text{ stabilizer of }^0\!\chi \text{ in }W\\\notag
&\Delta_{^0\!\chi}=\{\al\in\Delta\colon
(^0\!\chi)\circ(\check\al|_{\bO^\times})=1\}\\\notag
&\Delta_{^0\!\chi}^+=\Delta_{^0\!\chi}\cap\Delta^+\\\notag
&W_{^0\!\chi}^0=W(\Delta_{^0\!\chi})\\\notag &R^c_{^0\!\chi}=\{w\in
W_{^0\!\chi}\colon w\Delta_{^0\!\chi}^+=\Delta_{^0\!\chi}^+\}.
\end{align}
We call $\Delta_{^0\!\chi}$ the {\it good roots} for $^0\!\chi.$ It is
easy to see that $W_{^0\!\chi}=W_{^0\!\chi}^0\rtimes R_{^0\!\chi}.$

In \cite{Ro}, one associates to $^0\!\chi$ a particular compact open
subgroup $J_{^0\!\chi}$ of $G$ and a character $\rho$ of $J_{^0\!\chi}$.
The pair $(J_{^0\!\chi},\rho)$ satisfies (among other properties)
$J_{^0\!\chi}\cap H=\ ^0\!H$ and $\rho|_{^0\!H}=\ ^0\!\chi.$

\begin{theorem}[\cite{Ro}] Let $G_{^0\!\chi}$ denote the split $p$-adic
  group determined by $\Delta_{^0\!\chi}$ and $\C H(G_{^0\!\chi}//I_{^0\!\chi})$ the
  corresponding Iwahori-Hecke algebra.

\begin{enumerate}
\item There exists a family of $*$-preserving algebra isomorphisms $\CH(G//J_{^0\!\chi},\rho)\cong
  \CH(G_{^0\!\chi}//I_{^0\!\chi},1)\rtimes R_{^0\!\chi}$. (Since $R_{^0\!\chi}$ acts on the
  root datum of $G_{^0\!\chi}$, we can form this extended Hecke algebra.)
\item The group $R_{^0\!\chi}$ is abelian.
\item An irreducible representation of $G$ is in $\C C(J_{^0\!\chi},\rho)$ if and
  only if it is an irreducible subquotient of the ramified principal
  series $X(\chi).$
\end{enumerate}

\end{theorem}

In conclusion, we see that the representation theory of the
unramified principal series is controlled by certain extended
Iwahori-Hecke algebras.

\smallskip

One can make a connection with  table \ref{table:fine} for split
real groups. Recall that, in the real case, every fine $M$-type
$\delta\colon (\bZ/2\bZ)^{\rank G}\to\bC^\times$ can be viewed as a
character  of $(\bZ/2\bZ)^{\rank G}$, or in other words as a $(\rank
G)$-tuple of characters of $\bZ/2\bZ$:
$\delta=(\delta_1,\dots,\delta_{\rank
  G}),$ where $\delta_j$ is either the trivial or the sign character of
$\bZ/2\bZ.$ Similarly, every character $^0\!\chi\colon
^0\!H\to\bC^\times$ can be viewed as  a $(\rank G)$-tuple of
characters of $\bO^\times\colon ^0\!\chi=(^0\!\chi_1,\dots,^0\!\chi_{\rank
G}),$ where $^0\!\chi_j\colon \bO^\times\to\bC^\times.$

Fix a nontrivial quadratic character $^0\!\chi_0$ of $\bO^\times$.
Then define the correspondence
\begin{align}
\delta\longrightarrow \, ^0\!\chi_{\delta},\, \text{ with }
(^0\!\chi_{\delta})_j=
\begin{cases} 1  & \text{ if }\delta_j=triv,\\
^0\!\chi_0         &  \text{ if }\delta_j=sgn.
\end{cases}
\end{align}

Then all the data associated with $\delta$ in the real case is
identical with the data associated to $^0\!\chi_\delta$ in the
$p$-adic case. In principle, by an extension of the results in
\cite{BM1,BM2}, one expects that the unitary representations in the
ramified principal series in $\C C(J_{^0\!\chi},\rho)$ should correspond
to the unitary dual of the extended algebra $\CH(G_{^0\!\chi}//I_{^0\!\chi},1)\rtimes
R_{^0\!\chi}$. The (possibly surprising) observation is that the unitary
representations in the nonspherical minimal principal series of the
real group are closely related to the unitary dual of the same Hecke
algebra. This will become apparent in section \ref{sec:4}.

\subsection{Definitions}\label{sec:3.1} We recall some basic
results on  graded Hecke algebras. The {\it affine graded Hecke
algebra} $\bH$ was introduced in \cite{L1}. We will only need to
consider a special case of the definition. The generators of $\bH$
are the elements $\{t_{s_{\al}}\colon \al\in\Pi\}$  and
$\{\omega\colon \omega\in \fh\}.$ Here $\Pi$ denotes the set of
simple roots. As a $\bC$-vector space,
\begin{equation}\label{eq:3.1.1}
\bH=\bC[W]\otimes \bA,
\end{equation}
where
\begin{equation}\label{eq:3.1.2}
\bA=Sym(\fh).
\end{equation}
The following commutation relations hold:
\begin{equation}
\omega t_{s_\al}=t_{s_\al} s_\al(\omega)+\langle\omega,\al\rangle,
\quad \al\in\Pi,\ \omega\in\fh.
\end{equation}
The center $Z(\bH)$ of $\bH$ consists  of the $W$-invariants in
$\bA$ (\cite{L1}):
\begin{equation}\label{eq:3.1.4}
Z(\bH)=\bA^W.
\end{equation}

 On any irreducible $\bH$-module, which is necessarily
finite dimensional, $Z(\bH)$ acts by a {\it central
  character}. Therefore, the central characters are parameterized by
$W$-conjugacy classes in $\fh$. We will only consider {\it real}
central characters, \ie $W$-conjugacy classes of hyperbolic
semisimple elements in $\check \fh\cong \fh^*.$

We say that a module $V$ of $\bH$ is {\it spherical} if the
restriction of $V$ to $\bC[W]$ contains the trivial
$W$-representation.

\smallskip

Let $\chi\in \check\fh$ be a hyperbolic semisimple element. Denote
by  $\bC_\chi$ the  corresponding  character of $\bA$.
 Then one can form the {\it
  spherical principal series} module
\begin{equation}\label{eq:3.1.5}
X(\chi)=\bH\otimes_{\bA}\bC_\chi.
\end{equation}
It is immediate from the definition that
\begin{equation}\label{eq:3.1.6}
X(\chi)\cong \bC[W]\quad\text{as $W$-modules,}
\end{equation}
hence  $X(\chi)$ contains the trivial $W$-representation with
multiplicity one. Therefore, $X(\chi)$ has a unique spherical
subquotient, $L(\chi).$ The following result is well-known, it is
the Hecke algebra equivalent of a classical result for groups.

\begin{theorem}Let $X(\chi)$ be the spherical principal series defined
  in (\ref{eq:3.1.1}).
\begin{enumerate}
\item If $\chi$ is dominant, respectively antidominant, then $X(\chi)$
  has a unique irreducible quotient, respectively submodule, which
  is $L(\chi).$
\item Every irreducible spherical $\bH$-module is isomorphic to
  $L(\chi),$ for some $\chi\in\check\fh.$
\end{enumerate}

\end{theorem}

\subsection{Operators}\label{sec:3.2} We would like to consider Hermitian and unitary
modules for $\bH.$ For this, we need a $*$-operation on $\bH.$ This is
given by \cite{BM3}:
\begin{align}\label{eq:3.2.1}
&t_{w}^*=t_{w^{-1}},\quad w\in W,
&\omega^*=-\overline\omega+\sum_{\al\in\Delta^+}\langle\omega,\al\rangle
t_{s_\al},\quad \omega\in\fh.
\end{align}

For every $\al\in\Pi,$ set
\begin{equation}\label{eq:3.2.2}
r_\al=t_{s_\al}\check\al-1
\end{equation}
(an element of $\bH$). Then, if $w\in W$ and $w=s_{\al_1}\cdot
s_{\al_2}\cdot\dots\cdot s_{\al_k}$ is a minimal decomposition of
$w$ in $W$, one can define
\begin{equation}\label{eq:3.2.3}
r_w=r_{\al_1}\cdot r_{\al_2} \cdot\,  \dots \, \cdot r_{\al_k}.
\end{equation}
Note that $r_w$ is independent of the choice   of the reduced
decomposition (\cf~\cite{BM3}),  and therefore is well-defined.
When $w=w_0$ (the long Weyl group element), we use $r_{w_0}$ to
define the {\it long intertwining
  operator}:
\begin{equation}\label{eq:3.2.4}
A(w_0,\chi)\colon  X(\chi)\longrightarrow X(w_0\chi),\quad
t_w\otimes \one_\chi\mapsto t_wr_{w_0}\otimes \one_{w_0\chi},\ w\in
W.
\end{equation}
We obtain the same formula if we replace the $r_{w_0}$ by the
element $r_{w_0}(\chi)\in\bC[W] $ defined as follows. If
\[
w_0=s_{\al_1}\cdot s_{\al_2} \cdot\,  \dots \, \cdot s_{\al_l}\] is
a minimal decomposition of $w_0$ in $W$ (with $\ell=\#\Delta^+$),
then
\begin{equation}
r_{w_0}(\chi)=r'_{\al_1}\cdot r'_{\al_2} \cdot\,  \dots \, \cdot
r'_{\al_\ell}
\end{equation}
with
\begin{equation}\label{eq:3.2.5}
r'_{\al_j}=-\langle\check\al_j,(s_{\al_{j+1}}s_{\al_{j+2}}\dots
s_{\al_\ell})\chi\rangle \,t_{s_{\al_j}}-1.
\end{equation}
%Moreover, the $\check\al_j$'s of the $r_{\al_j}$-factors in the
%decomposition $r_{w_0}=r_{\al_1}\cdot r_{\al_2}\cdot\dots\cdot
%r_{\al_\ell},$ where $\ell=|\Delta^+|,$ can be replaced by scalars
%\begin{equation}\label{eq:3.2.5}
%-\langle\check\al_j,(s_{\al_{j+1}}s_{\al_{j+2}}\dotsb s_{\al_\ell})\chi\rangle.
%\end{equation} We denote the resulting element by $r_{w_0}(\chi).$
%Note that $r_{w_0}(\chi)\in\bC[W].$

The following discussion is again the Hecke algebra analogue of a
classical result for groups, we refer the reader for example to
\cite{BM3} for more details.

\begin{lemma}\label{l:3.2} Let $A(w_0,\chi)$ be the operator defined in
  (\ref{eq:3.2.4}).
Then $A(w_0,\chi)$ is an intertwining operator, and  is polynomial
(hence entire) in $\chi.$
\end{lemma}

Let $(\psi,V_\psi)$ be a $W$-type. The operator $A(w_0,\chi)$
induces operators on $\Hom$ spaces
\begin{equation}\label{eq:3.2.7}
A_\psi(w_0,\chi) \colon \Hom_W(V_\psi,X(\chi))\longrightarrow
\Hom_W(V_\psi,X(w_0\chi)).
\end{equation}
Furthermore, by Frobenius reciprocity and (\ref{eq:3.1.6}), this
induces an operator
\begin{equation}\label{eq:3.2.8}
A_\psi(w_0,\chi)\colon V_\psi^*\longrightarrow V_\psi^*.
\end{equation}

The operator $A_{triv}(w_0,\chi)$ is the scalar
\begin{equation}\label{eq:3.2.9}
N(\chi)=\prod_{\al\in\Delta^+} (\langle\check\al,\chi\rangle+1).
\end{equation}
We normalize the operators by this scalar and call them
$\CA(w_0,\chi),$ respectively $\CA_\psi(w_0,\chi).$ Then
$\CA_{triv}(w_0,\chi)=1.$

\begin{proposition}\label{p:3.2} Let $\CA(w_0,\chi)$ be the operator
  defined above. Assume that $\chi$ is dominant.

\begin{enumerate}
\item The operator $\CA(w_0,\chi)$ has no poles.
\item The image of $\CA(w_0,\chi)$ is $L(\chi).$
\item The Hermitian dual of $L(\chi)$ is $L(-\chi)$. Therefore
  $L(\chi)$ is Hermitian if and only if $w_0\chi=-\chi.$ In this case, if
  $(\, ,\, )_h$ denotes the Hermitian pairing, then the Hermitian form
  on $L(\chi)$ is
\begin{equation}\label{eq:3.2.6}
(t_x\otimes \one_\chi,t_y\otimes \one_\chi):=(t_x\otimes \one_\chi,\frac 1{N(\chi)}
t_yr_{w_0}(\chi)\otimes \one_{-\chi})_h,\quad x,y\in W.
\end{equation}
\end{enumerate}

\end{proposition}

It is easy to see that $L(0)$ is (irreducible and) unitary. In
general, we have the following unitarity criterion.

\begin{corollary}\label{c:3.2} Assume $\chi$ is dominant and $w_0\chi=-\chi.$
The spherical irreducible module $L(\chi)$ is unitary if and only if
the operators $\CA_{\psi}(w_0,\chi)$ are positive semidefinite for
all $W$-types $\psi.$
\end{corollary}

Clearly, this shows that the computation of the spherical unitary
dual, for any given Hecke algebra $\bH$, is a finite problem.

\medskip

\noindent{\bf Example}.

(a) The operator on the sign $W$-type is easy to compute:
\begin{equation}\label{eq:3.2.10}
\CA_{sign}(w_0,\chi)=\prod_{\al\in\Delta^+}\frac{1-\langle\check\al,\chi\rangle}{1+\langle\check\al,\chi\rangle}.
\end{equation}
If $\chi$ is dominant, the module $X(\chi)$  is irreducible if and
only if $A_{sign}(\chi)\neq 0,$ that is, if and only if
\begin{equation}
\langle\check\al,\chi\rangle\neq 1, \text{ for all } \al\in\Delta^+.
\end{equation}
This is the Hecke algebra analogue of the statement that the
spherical quotient of the unramified principal series of a split (real
or adjoint $p$-adic) group admits Whittaker models if and only if it is the full principal series.

(b) Let us consider the simplest case, \ie the Hecke algebra of type
$A_1.$ Then there are two representations of $W$, the trivial and
the sign. There is a single positive root $\al.$ From corollary
\ref{c:3.2} and formula (\ref{eq:3.2.10}), it follows that $L(\chi)$
is unitary if and only if
\begin{equation}\label{eq:3.2.11}
-1\le\langle\check\al,\chi\rangle\le 1.
\end{equation}

\subsection{} Let $
w_0=s_{\al_1}\cdot s_{\al_2} \cdot\,  \dots \, \cdot s_{\al_l}$ be a
minimal decomposition of $w_0$ in $W$, and let
\begin{equation}\CA_\psi(w_0,\chi)=\CA_{\psi}(\al_1,w_1\chi)\dotsb
\CA_{\psi}(\al_\ell,w_\ell\chi)\end{equation} be the corresponding
decomposition of $\CA_\psi(w_0,\chi)$, with
$w_j=(s_{\al_{j+1}}s_{\al_{j+2}}\dotsb s_{\al_\ell})$. Then
\begin{equation}\label{3.3a.1}
\CA_{\psi}(\al_j,\nu)=\left\{\begin{matrix} 1  &\text{on the }
    (+1)\text{-eigenspace of } s_{\al_j} \text{ on } V_\psi^*\\
\frac
{1-\langle\check\al_j,\nu\rangle}{1+\langle\check\al_j,\nu\rangle}&\text{on
the }
    (-1)\text{-eigenspace of } s_{\al_j} \text{ on } V_\psi^*.\end{matrix}\right.
\end{equation}
If $\al$ is a simple root, we have the formula
\begin{equation}
t_{s_\al} r_w=r_wt_{s_{w^{-1}\al}}.
\end{equation}
>From this, since $s_{w^{-1}\al}=w^{-1}s_\al w$, it follows that
\begin{equation}\label{3.3a.2}
t_{w}r_{w}=r_{w}t_{w}, \text{ for all } w\in W.
\end{equation}

\begin{remark}\label{r:3.3a}
In particular, for $w=w_0$, we obtain  that every operator
$\CA_\psi(w_0,\chi)$ preserves the $(+1)$ and the $(-1)$ eigenspaces
of $w_0$ on $\psi^*.$
\end{remark}

\subsection{Relevant $W$-types}\label{sec:3.4} In light of   corollary \ref{c:3.2},
a spherical module $L(\chi)$ is unitary if and only if
$\CA_\psi(w_0,\chi)$ is positive semidefinite, for {\it all}
$W$-types $\psi\in \widehat{W}.$

In fact, one determines in \cite{Ba1},\cite{BC2},\cite{Ci} a strict
subset of $\widehat W$ which is sufficient to detect unitarity. We
call this  set  {\it the relevant $W$-types}.  Of course, relevant
$W$-types  are interesting if and only if they form a small set.
(This is the case, for example, for $W(E_8)$: while there are $112$
irreducible representations of $W(E_8)$, our set of relevant
$W(E_8)$-types contains only nine representations.)

 In
general, it is still unclear what the best way to define uniformly
this set would be. For example, one can notice, after the
calculations are done, that a possible set of relevant $W$-types
consists of representations which in the Springer's correspondence
(\cite{Sp}) are   attached to nilpotent orbits of level $4$ or less.
(One says that a nilpotent element $\check e$ has level $k$ if
$ad(\check e)^{k+1}=0,$ but $ad(\check e)^k\neq 0.$)

We now provide lists of relevant $W$-types. They will play  a role
in section \ref{sec:4.6}.

\begin{definition}\label{d:3.4} The following sets of $W$-types are called
  {\it relevant}:

\begin{enumerate}
\item[$A_{n-1}:$] $\{(n-m,m):\ 0\le m\le [n/2]\}$;
\item[$B_n,C_n:$] $\{(n-m,m)\times (0):\ 0\le m\le
  [n/2]\}\cup\{(n-m)\times (m): 0\le m\le n\}$;
\item [$D_n:$]  $\{(n-m,m)\times (0):\ 0\le m\le
  [n/2]\}\cup\{(n-m)\times (m): 0\le m\le [n/2]\}$;
\item [$G_2:$] $\{1_1,2_1,2_2\}$;
\item [$F_4:$] $\{1_1,4_2,2_3,8_1,9_1\}$;
\item[$E_6:$] $\{1_p,6_p,20_p,30_p,15_q\}$;
\item[$E_7:$] $\{1_a,7_a',27_a,56_a',21_b',35_b,105_b\}$;
\item[$E_8:$] $\{1_x,8_z,35_x,50_x,84_x,112_z,400_z,300_x,210_x\}$.
\end{enumerate}
The notation for  Weyl group representations  is as in \cite{Ca}.

\end{definition}

\subsection{Extended Hecke algebras}\label{sec:3.7}

We present some elements of the construction of  extended graded
Hecke algebras. This construction will be applied to the setting of
a Hecke algebra constructed from the set of good coroots and the
(dual of the) $R$-group, as in section \ref{sec:3.0a}.

Let $\bH$ denote the graded Hecke algebra associated to a root
system as in the previous section, and let $R$ be a subgroup of the
group of automorphisms of the root system for $\bH.$

\begin{definition}\label{d:3.7.1}
We define $\bH'$ to be the semidirect product
\begin{equation}
  \label{eq:3.7.1}
  \bH':=\bC[R]\ltimes \bH,
\end{equation}
where the action of $R$ on $\bH$ comes from the action of $R$ by outer
automorphisms on the root system of $\bH$.
\end{definition}

Set
\begin{equation}W':=R\ltimes W.
\end{equation}
In the same way as for usual graded Hecke algebras,  one  obtains
that:

\begin{lemma}\label{l:3.7.1}
The center of $\bH'$ is $\bA^{W'}.$
\end{lemma}

For every $\nu\in\fh^*,$ we fix the following notation:
\begin{align}\label{eq:3.7.2}
&R(\nu)=\text{the centralizer of }\nu\text{ in }R,\\\notag
&\bA'(\nu)=\bC[R(\nu)]\ltimes\bA,\\\notag
&\bH'(\nu)=\bC[R(\nu)]\ltimes\bH.
\end{align}

Let $triv_W$ be the trivial representation of $W$. This is
stabilized by $R$, so by Mackey induction, any representation
$\lambda$ of $R$  gives rise to a representation of $W',$ which we
denote by $\lambda\ltimes triv_W$; any $W'$-representation
containing $triv_W$ in its restriction to $W$ is obtained in this
way.

 Call an irreducible $W'$-representation  \emph{fine} if it is
of the form $\lambda\ltimes triv_W$ for some $\lambda\in \widehat
R.$ Note that if $\rho=\lambda\ltimes triv_W$ is a fine $W'$-type,
then $\Hom_W[\rho:triv_W]\cong \lambda^*.$

\begin{definition}\label{d:3.7.2}
We call a module $\pi$ of $\bH'$ {\it quasi-spherical} if
$\Hom_W[\pi:triv_W]\neq 0.$
\end{definition}
Clearly, any module $\pi$ containing a fine $W'$-type is
quasi-spherical.
Consider the principal series
\begin{equation}  \label{eq:3.7.3}
X'(\nu)=\bH'\otimes_{\bA'(\nu)}\bC_\nu=\bH'\otimes_{\bH'(\nu)}(\bH'(\nu)\otimes_{\bA'(\nu)}\bC_\nu).
\end{equation}

Every fine $W'$-type $\rho=\lambda\ltimes triv_W$ appears in
$X'(\nu)$ with multiplicity $\dim\lambda.$ Our main case of interest
is when $R$ is abelian, in this case all the multiplicities are one.

We can extend the definition of intertwining operators to this
setting. Assume $u w^0\in R\ltimes W.$ Then, similarly to section
\ref{sec:3.2}, we define a $\bH'$-operator
\begin{equation}\label{eq:3.7.5}
A'(u w^0, \nu)\colon X'(\nu)\to X'(u w^0\nu),\ x\otimes
\one_\nu\mapsto xu r_{w^0}\otimes \one_{u w^0\nu}.
\end{equation}
 We normalize the operator $A'(u w^0, \nu)$
 to be the identity on the trivial $W'$-type, and denote the
 resulting operator by $\CA'(u w^0, \nu)$.

For every $W'$-type $\psi'$, we obtain an operator
\begin{equation}\label{eq:3.7.6}
\CA'_{\psi'}(u w^0,\nu) \colon
\Hom_{W'}(\psi',X'(\nu))\longrightarrow\Hom_{W'}(\psi',X'(u
w^0\nu)).
\end{equation}

\begin{remark}\label{r:3.7}
When $w^0\nu=-\nu,$  the $\bH'$-operator $\CA'_{\psi'}(w^0,\nu)$ is
block diagonal, each block corresponding to a representation
$\lambda$ of $R(\nu)$, such that it is the restriction to
$\Hom_{R(\nu)}(\psi',\lambda)$ of $\CA_{\psi'|_W}(w^0,\nu).$
\end{remark}

% As hinted before, the construction of the
%previous section will be applied to the root system $\Delta_\delta$ of
%the good roots and the dual of $\widehat R_\delta$ of the R-group.

\section{Petite $K$-types for split real groups}\label{sec:3a}

In this section, we discuss the construction of petite $K$-types and
the relation between the real intertwining operators from section
\ref{sec:2}, and the graded Hecke algebra operators from section
\ref{sec:3}.

%\smallskip

\subsection{Operators for real split groups}\label{sec 4.1}
Assume that $G$ is a real split group, let $P=MAN$ be  a minimal
parabolic subgroup of $G$. Choose a (fine) $K$-type $\delta$ and a
weakly dominant character $\nu$ of $A$. By the results in section
\ref{sec:Hermitian}, a Langlands quotient $L(\delta,\nu)$ is
Hermitian if and only if there is a Weyl group element satisfying
\[ w\cdot Q = \bar{Q}, \quad w \cdot \delta \simeq \delta \quad \text{ and } \quad w \cdot \nu  = - \bar{\nu}\]
where $Q$ is the the parabolic defined by $\nu$, as in section
\ref{sec:2.9}. (If $\nu$ is strictly dominant, then $Q=P$ and
$w=w_0$, the long Weyl group element.) The intertwining operator
\[
\mathcal{A}_{\mu}(w,\delta,\nu)=\mu_{\delta}(w)
A'_{\mu}(w,\delta,\nu) \colon \Hom_{M}(\mu,\delta) \rightarrow
\Hom_{M}(\mu,\delta)
 \]
induces a nondegenerate invariant Hermitian  on $L(\delta,\nu)$.

The operator $\mathcal{A}_{\mu}(w,\delta,\nu)$ has been introduced
in sections \ref{sec:2.2} and \ref{sec:2.10}. We now describe its
properties in greater details.

\begin{theorem}\label{t:4.1}
Let $w=s_{\alpha_{r}}s_{\alpha_{r-1}}\ldots s_{\alpha_{1}}$ be a
minimal decomposition of $w$ as a product of simple reflections.
Then
\begin{equation}\label{eq:3a.1.1}
\mathcal{A}_{\mu}(w,\delta,\nu)=\prod\limits_{j=1}^{r}
\mathcal{A}_{\mu}(s_{\alpha_{j}},\rho_{j-1},\nu_{j-1})
\end{equation}
with
\begin{equation}\label{eq:3a.1.2}
\mathcal{A}_{\mu}(s_{\alpha_{j}},\rho_{j-1},\nu_{j-1})=
\mu_{\delta}(\sigma_{\alpha_j})A'_{\mu}(
s_{\alpha_{j}},\rho_{j-1},\nu_{j-1})
\end{equation}
 and $\rho_j$
equal to the unique copy of $\delta_j$ inside the fine $K$-type
$\mu_{\delta}$.
\end{theorem}

\begin{proof}
Recall that $A'_{\mu}(  w,\delta,\nu)$  is a normalization of the
standard intertwining operator $A_{\mu}(  w,\delta,\nu)$ introduced
in section \ref{sec:2.2}. As in (\ref{eq:2.2.4}), it has  a decomposition of the form:
 \begin{equation}\label{eq:3a.1.3}
A'_{\mu}(  w,\delta,\nu) =\prod\limits_{j=1}^{r}A'_{\mu}(
s_{\alpha_{j}},\delta_{j-1},\nu_{j-1})
 \end{equation}
with $x_0=1$, $\delta_{0}=\delta=x_0\cdot \delta,$ $\nu_{0}=\nu=x_0
\cdot \nu$, and
$ \delta_{j}=\underbrace{s_{\alpha_{j}}s_{\alpha_{j-1}}\ldots
s_{\alpha_{1}}}_{x_j}\cdot\, \delta= x_j\cdot\, \delta$,
$\nu_{j}=\underbrace{s_{\alpha_{j}}s_{\alpha_{j-1}}\ldots
s_{\alpha_{1}}}_{x_j}\cdot\, \nu= x_j\cdot\, \nu$,
for $j\geq 1$. The operator $\mathcal{A}_{\mu}(w,\delta,\nu)
=\mu_{\delta}(w)A'_{\mu}( w,\delta,\nu)$ inherits a similar
decomposition. We can write:
\begin{align*}
\mathcal{A}_{\mu}(w,\delta,\nu)  & = \mu_{\delta}(w) \left[
\prod\limits_{j=1}^{r}A'_{\mu}(
s_{\alpha_{j}},\delta_{j-1},\nu_{j-1})\right]
\\
& =\mu_{\delta}(x_r) \left[ \prod\limits_{j=1}^{r}A'_{\mu}(
s_{\alpha_{j}},x_{j-1}\cdot \delta,x_{j-1}\cdot\nu)\right]
\mu_{\delta}(x_0)^{-1}
\\
&=  \prod\limits_{j=1}^{r}\left[\mu_{\delta}(x_j)A'_{\mu}(
s_{\alpha_{j}},x_{j-1}\cdot \delta,x_{j-1}\cdot\nu)
\mu_{\delta}(x_{j-1})^{-1}\right]
\\
& =
\prod\limits_{j=1}^{r}\left[\mu_{\delta}(\sigma_{\alpha_j})A'_{\mu}(
s_{\alpha_{j}},\rho_{j-1},\nu_{j-1}) \right]\\
&= \prod\limits_{j=1}^{r}
\mathcal{A}_{\mu}(s_{\alpha_{j}},\rho_{j-1},\nu_{j-1}).
\end{align*}
Here $\rho_j$ denotes the unique copy of $\delta_j$ inside the
(fixed) fine $K$-type $\mu_{\delta}$.
We refer to the factorization
\begin{equation}\label{eq:3a.1.4}
\mathcal{A}_{\mu}(w,\delta,\nu)=\prod\limits_{j=1}^{r}
\mathcal{A}_{\mu}(s_{\alpha_{j}},\rho_{j-1},\nu_{j-1})
\end{equation}
as the Gindikin-Karpelevi\v c decomposition of the operator
$\mathcal{A}_{\mu}(w,\delta,\nu)$.
 Note that every factor
$$
\mathcal{A}_{\mu}(s_{\alpha_i},\rho_{i-1},\nu_{i-1}) \colon
\Hom_M(\mu,\rho_{i-1}) \longrightarrow \Hom_M(\mu,\rho_i)
$$
of this decomposition can be regarded as an operator on
$\Hom_M(\mu,\mu_{\delta})$. We summarize the factorization with the
following diagram. The horizontal lines on top represent the
decomposition (\ref{eq:3a.1.3}) of $A'_{\mu}( w,\delta,\nu)$. The
vertical lines are the operators $\mu_{\delta}(x_j)$'s and their
inverses; in particular, the first  and last vertical lines are the
identity and the operator $\mu_{\delta}(w)$ respectively. The
horizontal lines in the bottom represent the  factorization
(\ref{eq:3a.1.4}) of $\mathcal{A}_{\mu}( w,\delta,\nu).$\\

\begin{center}
\begin{picture}(300,10)(60,78)

\put(40,-24){\small{$\Hom_M(\mu,\delta)$}}
\put(40,87){\small{$\Hom_M(\mu,\delta)$}}

\put(56,77){$\shortparallel$} \put(56,-14){$\shortparallel$}

\put(37,-4){\small{$\Hom_M(\mu,\rho_0)$}}
\put(37,67){\small{$\Hom_M(\mu,\delta_0)$}}

\put(355,-24){\small{$\Hom_M(\mu,\delta)$}}
\put(345,87){\small{$\Hom_M(\mu,w\delta)$}}

\put(366,77){$\shortparallel$} \put(366,-14){$\shortparallel$}

\put(352,-4){\small{$\Hom_M(\mu,\rho_r)$}}
\put(352,67){\small{$\Hom_M(\mu,\delta_r)$}}

\put(26,33){\tiny{$\mu_{\delta}(x_0)^{-1}=1$}}
\put(346,33){\tiny{$\mu_{\delta}(x_r)$}} \put(369,33){\tiny{$=
\mu_{\delta}(w)$}}

\put(64,10){\vector(1,0){20}} \put(64,60){\vector(1,0){20}}
\put(88,10){\dots\dots \dots} \put(88,60){\dots\dots \dots}
\put(58,10){\circle*{4}} %\put(58,10){\circle{7}}
\put(58,60){\circle*{4}}  %\put(58,60){\circle{7}}

\put(58,15){\vector(0,1){40}}

\put(107,-4){\small{$\Hom_M(\mu,\rho_{j-1})$}}
\put(136,10){\circle*{4}} \put(145,10){\vector(1,0){60}}
 \put(224,10){\vector(1,0){60}}

\put(214,10){\circle*{4}}
\put(198,-4){\small{$\Hom_M(\mu,\rho_{j})$}}
\put(292,10){\circle*{4}}
\put(270,-4){\small{$\Hom_M(\mu,\rho_{j+1})$}} \put(301,10){\dots
\dots\dots} \put(341,10){\vector(1,0){20}} \put(368,10){\circle*{4}}

\put(107,67){\small{$\Hom_M(\mu,\delta_{j-1})$}}
\put(136,60){\circle*{4}} \put(145,60){\vector(1,0){60}}
 \put(224,60){\vector(1,0){60}}
\put(214,60){\circle*{4}}
\put(194,67){\small{$\Hom_M(\mu,\delta_{j})$}}
\put(292,60){\circle*{4}}
\put(266,67){\small{$\Hom_M(\mu,\delta_{j+1})$}} \put(301,60){\dots
\dots\dots} \put(341,60){\vector(1,0){20}} \put(368,60){\circle*{4}}
\put(368,55){\vector(0,-1){40}}

 \put(138,15){\vector(0,1){40}}
 \put(132,55){\vector(0,-1){40}}

 \put(211,55){\vector(0,-1){40}}
\put(217,15){\vector(0,1){40}}

  \put(288,55){\vector(0,-1){40}}
\put(295,15){\vector(0,1){40}}

\put(219,33){\tiny{$\mu_{\delta}(x_j)^{-1}$}}
 \put(189,33){\tiny{$\mu_{\delta}(x_j)$}}

\put(296,33){\tiny{$\mu_{\delta}(x_{j+1})^{-1}$}}
\put(258,33){\tiny{$\mu_{\delta}(x_{j+1})$}}

\put(140,33){\tiny{$\mu_{\delta}(x_{j-1})^{-1}$}}
\put(101,33){\tiny{$\mu_{\delta}(x_{j-1})$}}

\put(140,90){\tiny{$A'_{\mu}(s_{\alpha_j},\delta_{j-1},\nu_{j-1})$}}
\put(139,-22){\tiny{$\mathcal{A}_{\mu}(s_{\alpha_j},\rho_{j-1},\nu_{j-1})$}}

\put(225,90){\tiny{$A'_{\mu}(s_{\alpha_{j+1}},\delta_{j},\nu_{j})$}}
\put(232,-22){\tiny{$\mathcal{A}_{\mu}(s_{\alpha_{j+1}},\rho_{j},\nu_{j})$}}

\put(252,-12){$\Downarrow$}\put(253,-2){$\vdots$}\put(255,-2){$\vdots$}

\put(178,-12){$\Downarrow$}\put(179,-2){$\vdots$}\put(181,-2){$\vdots$}

\put(253,65){$\vdots$}\put(255,65){$\vdots$}
\put(252,78){$\Uparrow$}

\put(178,78){$\Uparrow$}\put(179,65){$\vdots$}\put(181,65){$\vdots$}

\end{picture}
\end{center}

\vskip1.4in
\end{proof}

Next, we show how to compute the various factors of the
decomposition (\ref{eq:3a.1.4}). Consider the action of $Z_{\alpha_j}^2$ on the space
$\Hom_M(\mu,\rho_{j-1})$ defined by
\begin{equation}
T \mapsto T \circ \text{d}\mu(Z_{\alpha_j})^2.
\end{equation}
 This action is well
defined because $Ad(m)Z_{\alpha_j}=\pm Z_{\alpha_j}$ for all $m\in
M$, so $Z_{\alpha_j}^2$ commutes with $M$. Every  eigenvalue of
$Z_{\alpha_j}^2$ on $\Hom_M(\mu,\rho_{j-1})$ is of the form
$(-l^2)$, with $l$ an integer. Precisely, $l$ is an even integer if
$\alpha_j$ is good for $\rho_{j-1}$ and an odd integer otherwise.
Write
\begin{equation}
 \Hom_M(\mu,\rho_{j-1}) := \bigoplus E^{\alpha_j}(-l^2)
 \end{equation}
  for
the decomposition of $\Hom_M(\mu,\rho_{j-1})$ in generalized
eigenspaces for $Z_{\alpha_j}^2$. \\ Then the operator
$\mathcal{A}_{\mu}(s_{\alpha_j},\rho_{j-1},\nu_{j-1})$ acts on
$E^{\alpha_j}(-l^2)$ by
\begin{equation}
T \mapsto c_{l}(\alpha_j, \nu_{j-1})
\,\mu_{\delta}(\sigma_{\alpha_j}) \, T \,
\mu(\sigma_{\alpha_j}^{-1}),
\end{equation}
where $\sigma_{\alpha_j}$ is a representative of $s_{\alpha_j}$ in
$K$. The constant $c_{l}(\alpha_j, \nu_{j-1})$ is equal  to 1 for
$l=0$ or $1$, and is given by \begin{align}
c_{2m+1}(\alpha_j,\nu_{j-1})&= (-1)^m \frac{(2-\lambda)(4-\lambda)
\cdots
(2m-\lambda)}{(2+\lambda)(4+\lambda) \cdots (2m+\lambda)} \\
 c_{2m}(\alpha_j,\nu_{j-1})&= (-1)^m \frac{(1-\lambda)(3-\lambda) \cdots
(2m-1-\lambda)}{(1+\lambda)(3+\lambda) \cdots (2m-1+\lambda)}
\end{align}
for $m>0$.  Here $\lambda=\langle \nu_{j-1},\check{\alpha}_j
\rangle.$ Notice the similarity with the formulas for the spherical
operators of $SL(2,\mathbb{R})$ in section  \ref{rank one}.
%This similarity should not be surprising:
Recall from section
\ref{sec:2.7.2} that
$\mathcal{A}_{\mu}(s_{\alpha_j},\rho_{j-1},\nu_{j-1})$ agrees with a
spherical  operator for the rank-one group $MG^{\alpha_j}$ on the
restriction of $\mu$ to $MK^{\alpha_j}$.

\begin{remark} For the purpose of computations, it is sometimes convenient to
construct simultaneously all the  operators
$\mathcal{A}_{\mu}(w,\delta',\nu)$ with $\delta'$ an $M$-type
occurring in $\mu_{\delta}$.

Let $\delta'$ be any $M$-type in the $W$-orbit of $\delta$. The
operator $\mathcal{A}_{\mu}(w,\delta' ,\nu)$, together with all its
factors
\[
\mathcal{A}_{\mu}(s_{\alpha_{j}},\rho'_{j-1},\nu_{j-1})\colon
\Hom_M(\mu,\rho_{j-1}') \longrightarrow \Hom_M(\mu,\rho'_{j}),
\]
can be regarded as   operators on $\Hom_M(\mu,\mu_{\delta})$. Let
\begin{equation}
\mathcal{A}_{\mu}(w,\nu)\colon \Hom_{M}( E_{\mu},E_{\mu_{\delta}})
\longrightarrow \Hom_{M}( E_{\mu},E_{\mu_{\delta}})
\end{equation}
be the direct sum (over all the $\delta'$ in the $W$-orbit of
$\delta$) of the operators $\mathcal{A}_{\mu}(w,\delta' ,\nu)$.

Via (\ref{eq:3a.1.4}), the operator $\mathcal{A}_{\mu}(w,\nu)$ has a decomposition:
\begin{equation}
\mathcal{A}_{\mu}(w,\nu)=\prod\limits_{j\geq 1}
\mathcal{A}_{\mu}(s_{\alpha_{j}},\nu_{j-1}).\end{equation} Every
factor is an operator    on $\Hom_{M}( E_{\mu},E_{\mu_{\delta}})$,
and is easy to compute. Write
\begin{equation}
\Hom_{M}( E_{\mu},E_{\mu_{\delta}})= \bigoplus\limits_{l \in
\mathbb{N}/2} E^{\alpha_j}(- l^2) \end{equation} for the
decomposition of $\Hom_{M}( E_{\mu},E_{\mu_{\delta}})$
 into eigenspaces for $ Z_{\alpha_j} ^2$.  Then   $\mathcal{A}_{\mu}(s_{\alpha_{j}},\nu_{j-1})$
acts on $E^{\alpha_j}(-l^2)$ by
\begin{equation}
T \mapsto c_l(\alpha_{j},\nu_{j-1})\,
\mu_{\delta}(\sigma_{\alpha_j}) T \mu(\sigma_{\alpha_j}^{-1}).
\end{equation}

\end{remark}

%\smallskip

\subsection{Operators on fine $K$-types}\label{sec 4.2}

Let $\mu_{\delta}$ be our (fixed) fine $K$-type containing $\delta$
(the same one used to define the operator $\mathcal{A}
(w,\delta,\nu)$). Our normalization imposes that
$\mathcal{A}_{\mu_{\delta}}(w,\delta,\nu)$  is  trivial. We will
show that, if $\mu$ is \emph{any fine $K$-type}  containing
$\delta$, then the operator
\[
\mathcal{A}_{\mu}(w,\delta,\nu)=\prod\limits_{j=1}^{r}
\mathcal{A}_{\mu}(s_{\alpha_{j}},\rho_{j-1},\nu_{j-1})
\]
is a (possibly nontrivial) scalar.

Choose $\mu$ in $A(\delta)$. To construct the $j^{th}$ factor of the
operator $\mathcal{A}_{\mu}(w,\delta,\nu)$, we look at the action of
$Z_{\alpha_j}^2$ on $\Hom_M(\mu,\rho_{j-1})$. Because $\mu$ is fine,
the only possible eigenvalues of $Z_{\alpha_j}^2$ on this space are
$0$ and $-1$, hence
\begin{equation}
\Hom_M(\mu,\rho_{j-1}) =\begin{cases} E_{\mu}^{\alpha_j}(0) & \text{
if $\alpha_j$ is good for $\rho_{j-1}$} \\
E_{\mu}^{\alpha_j}(-1) & \text{ otherwise} .\\
\end{cases}\end{equation}
In both cases,
 the operator
$\mathcal{A}_{\mu}(s_{\alpha_{j}},\rho_{j-1},\nu_{j-1})$ acts by $ T
\mapsto  \mu_{\delta}(\sigma_{\alpha_j}) \, T \,
\mu(\sigma_{\alpha_j}^{-1}). $ This follows immediately from  the
fact that $c_{0}(\alpha_j, \nu_{j-1})=c_{1}(\alpha_j, \nu_{j-1})
=1$.

 Then, the full intertwining operator
 \[
 \mathcal{A}_{\mu}(w,\delta,\nu) \colon \Hom_M(\mu,\delta)
 \longrightarrow  \Hom_M(\mu,\delta),\,T \mapsto \left[ \prod\limits_{j=1}^{r}
\mathcal{A}_{\mu}(s_{\alpha_{j}},\rho_{j-1},\nu_{j-1})
 \right]T
 \]
 acts by
 \begin{equation}
  T
\mapsto  \mu_{\delta}(\sigma_w) \, T \, \mu(\sigma_w),
\end{equation}
 where
$\sigma_w$ is a representative of $w$ in $K$.

Note that the space $\Hom_M(\mu,\delta)$ is one-dimensional, so
 $\mathcal{A}_{\mu}(w,\delta,\nu)$ is a scalar operator.
 We will now give a different interpretation of this scalar.

For  every (not necessarily fine) $K$-type $\pi$ containing
$\delta$, we   define a representation $\psi_{\pi}$ of $W_{\delta}$
on $\Hom_M(\pi,\delta)$, as follows.
 Since
$\delta$ appears in $\mu_{\delta}$ with multiplicity one, we can
 identify $\delta$ with its copy inside $\mu_{\delta}$,
and $V^{\delta}$ with $E_{\mu_{\delta}}(\delta)$ (the isotypic
component of $\delta$ in $\mu_{\delta}$). Let   $M'_{\delta}$  be
the preimage of $W_{\delta}$ in the normalizer $M'$ of $A$ in $K$.
The group  $M'_{\delta}$ acts on both $E_{\pi}(\delta)$ and
$E_{\mu_{\delta}}(\delta)$ by restriction of the appropriate action
of $K$,  hence it acts on
\[
\Hom_M(E_{\pi}, V^{\delta})=\Hom_M(E_{\pi},
E_{\mu_{\delta}}(\delta))=\Hom_M(E_{\pi}(\delta),
E_{\mu_{\delta}}(\delta))
\]
by
\begin{equation}
\sigma \cdot T(v) = \mu_{\delta}(\sigma) T (\pi(\sigma)^{-1}v) \quad
\forall \, v \in E_{\pi}. \end{equation}
Because $M$ acts trivially
on $\Hom_M(E_{\pi}, V^{\delta})$, we also get a representation
$\psi_{\pi}$ of $W_{\delta}$.
 The restriction of $\psi_{\pi}$ to $W_{\delta}^0$ will be
denoted by $\psi_{\pi}^0.$

Now assume that $\pi$ is  fine.  Let $(M'_{\delta})^0$ be the
subgroup of  $M'_{\delta}$ generated by the elements
$\sigma_{\alpha}$ for all the good real roots $\alpha$. Notice that
$(M'_{\delta})^0$ can be  smaller than the preimage of
$W_{\delta}^0$ in $M'$, because   it may not contain all of $M$.

The group $(M'_{\delta})^0$ acts trivially on the $\delta$-isotypic
component of every fine $K$-type $\pi$ containing $\delta$. To prove
this claim we observe that
 \[
\pi(\sigma_{\alpha}^2) v=  \pi(m_{\alpha}) v = \delta(m_{\alpha})v=v
\quad  \forall \, v \in E_{\pi}(\delta)
\]
for all $\sigma_{\alpha}\in (M'_{\delta})^0$ (recall that $\alpha$
is a good real root for $\delta$). So  $m_{\alpha}$  acts trivially
on $E_{\pi}(\delta)$. Because $\mu$ is fine  and has level at most
one, the $(+1)$-eigenspace of $\sigma_{\alpha}$ coincides with the
 $(+1)$-eigenspace
of $m_{\alpha}$  (and  the $(0)$-eigenspace to $Z_{\alpha}$). In
particular,
 $\sigma_{\alpha}$ acts trivially on $E_{\pi}(\delta).$

Let us go back to the space
\[ \Hom_M(\pi, \delta)= \Hom_M(E_{\pi},
V^{\delta})=\Hom_M(E_{\pi}(\delta), E_{\mu_{\delta}}(\delta)).
\]
Since $\pi$ and $\mu_{\delta}$ are both fine, the action of
$(M'_{\delta})^0$ on this space is obviously trivial. The group
$W_{\delta}^0 = [M(M'_{\delta})^0]/M$ also acts trivially on
$\Hom_M(\pi, \delta).$

Hence, for every fine $K$-type $\pi$ containing $\delta$, there is a
representation of the $R$ group $R_{\delta}=W_{\delta}/W_{\delta}^0$
on $\Hom_M(E_{\pi}, V^{\delta}).$ We denote this representation by
$\varrho_{\pi}$. If $[w]\in R_{\delta}$ and $\sigma$ is a
representative for $w$ in $M'_{\delta}$, then
\begin{equation}
(\varrho_{\pi}[w]T )(v) = \mu_{\delta}(\sigma) T (\pi(\sigma)^{-1}v)
\quad \forall \, v \in E_{\pi}.
\end{equation}

Notice that this is exactly the action of
$\mathcal{A}_{\pi}(w,\delta,\nu)$ on $\Hom_M(E_{\pi}, V^{\delta}).$

\begin{remark}
The representations $\psi_{\mu}$, $\psi_{\mu}^0$ and $\varrho_{\pi}$
depend on the choice of the (fixed) fine $K$-type $\mu_{\delta}\in
A(\delta)$. This is the same fine $K$-type used to define the
operator $\mathcal{A}(w,\delta,\nu)$.
\end{remark}

\begin{proposition}\label{prop:4.2} Choose a minimal parabolic subgroup
$P=MAN$   of $G$, a  representation $\delta$ of $M$ and a weakly
dominant character $\nu$ of $A$. Suppose that $X(\delta,\nu)$ is
Hermitian, and let $w\in W_{\delta}$ be such that $w\cdot Q =
\bar{Q} $ and $w\cdot \nu = - \bar{\nu}$. Denote by $\sigma$ a
representative for $w$ in $M'_{\delta}$ and by $[w]$ the equivalence
class of $w$ in $R_{\delta}$.

 Having  fixed a fine  $K$-type
$\mu_{\delta}\in A(\delta)$, we
 associate a
character $\varrho_{\pi}$ of   $R_{\delta}$ to every fine $K$-type
$\pi $ containing $\delta$,  as above. The intertwining operator
\[
 \mathcal{A}_{\pi}(w,\delta,\nu) = \mu_{\delta}(w) A'_{\pi}(w,\delta,\nu)
\]
acts
  on
$\Hom_M(\pi, \delta)$ by the scalar $\varrho_{\pi}[w].$

\end{proposition}

\begin{corollary}\label{c:4.2}
If $w\in W_{\delta}^0$, the operator
$\mathcal{A}_{\pi}(w,\delta,\nu)$ is trivial  for all $\pi \in
A(\delta).$
\end{corollary}

If  $w \not \in W_{\delta}^0$ (and $R_{\delta}\neq \{1\}$), then
different
 fine $K$-types may get
different sign.

%\smallskip

\subsection{A nonunitarity criterion} A first criterion for
nonunitarity is obtained by analyzing the signature of the
intertwining operator on the fine $K$-types.

\begin{proposition}\label{prop 4.3}
Let $G$ be a real split group. Choose a minimal parabolic subgroup
$P=MAN$   of $G$, a  representation $\delta$ of $M$ and a
\emph{strictly} dominant character $\nu$ of $A$. Let $L(\delta,\nu)$
be the corresponding irreducible representation of $G$. Suppose that
$w\cdot P = \bar{P} $ and $ w \cdot \nu  = - \bar{\nu}$ for some
$w\in W_{\delta}$, so that $L(\delta,\nu)$ is Hermitian, and assume
that there is no element $w^0\in W_{\delta}^0$ with the same
property. Then $L(\delta,\nu)$ is not unitary.
\end{proposition}

\begin{proof} Since $A(\delta)=\# R_{\delta} $ and $R_{\delta}$ is nontrivial,
 the principal series
$X(\delta,\nu)$ contains several fine $K$-types. Every   fine
$K$-type of $X(\delta,\nu)$  is contained in $L(\delta,\nu)$,
because $L(\delta,\nu)$ is the unique irreducible subquotient of
$X(\delta,\nu)$. So, to prove that $L(\delta,\nu)$ is nonunitary
(\ie $\mathcal{A}(w,\delta,\nu)$ is not positive semidefinite), it
is sufficient to show that there is a  fine  $K$-type $\pi$ in
$A(\delta)$ such that
$\mathcal{A}_{\pi}(w,\delta,\nu)=\varrho_{\pi}[w]< 0.$

Because $R_{\delta}$ is an abelian two-group, the dual
$\widehat{R_{\delta}}$ is a group isomorphic to $R_{\delta}$. In
particular, the
  number of distinct characters of $R_{\delta}$ equals the number of
fine $K$-types containing $\delta$. Also notice that,  with the
notation of the previous section, $\varrho_{\pi_1}=\varrho_{\pi_2}$
if and only if $\pi_1=\pi_2$. So every character of $R_{\delta}$ is
of the form  $\varrho_{\pi}$ for some $\pi \in A(\delta)$. We
conclude that there must be at least one fine $K$-type $\pi$ in
$A(\delta)$ such that $\varrho_{\pi}[w]=-1$, and we are done.
\end{proof}

Another explanation (suggested by D. Vogan) for why this result should
hold is as follows.
Choose $G_2\supset G$ as disconnected as possible (if $G$ is the
real points of an algebraic group, then $G_2$ is the preimage of the
real points of the adjoint group of $G$), and let $M_2A_2N_2$ be a
minimal parabolic subgroup of $G_2$ containing $MAN$. The group
$M_2$ is still abelian and contains $M$;  we can choose an extension
$\delta_2$ of $\delta$ to   $M_2$ such that $W_{\delta_2}=
W_{\delta}^0.$ Then the principal series $X(\delta_2,\nu)$ of $G_2$
is not Hermitian, and its restriction to $G$ cannot be unitary.
(Note that every fine $K$-type containing $\delta$ appears in the
restriction of the unique fine $K_2$-type containing $\delta_2$.)

\medskip

We now give some  applications of the proposition.
\begin{enumerate}

\item If $\nu=a>0$, the nonspherical Langlands quotient $L(sign,a)$ of $SL(2,\mathbb{R})$ is
Hermitian. The   element $w=s_{\alpha} \in W_{sign}$ carries $P$
into $\bar{P}$ and $\nu$ into $-\nu$. Because
 $W_{sign}^0=\{1\}$,  there is
no element in $W_{sign}^0$ with the same property.   It is not hard
to check that the operator $\mathcal{A}(s_{\alpha},sign,a)$ takes
opposite sign on the two fine $SO(2)$-types ($+1$ and $-1$)
containing $sign$. Both fine $SO(2)$-types are contained in
$L(sign,a)$, hence
 $L(sign,a)$ is not unitary. Notice that the corresponding representation of
 $G_2=SL^{\pm}(2,\mathbb{R})$ is not even Hermitian.

\item If $\nu$ is real and strictly dominant, the nonspherical Langlands quotient $L(\delta_p,\nu)$ of
$Sp(2n,\mathbb{R})$ is Hermitian (\cf~section \ref{sec:2.5}). The
element $w=-I \in W_{\delta_p}=C_{n-p}\times C_p $ carries $P$ into
$\bar{P}$ and $\nu$ into $-\nu$. Because
 $W_{\delta_p}^0=C_{n-p}\times D_p$, when $p$ is odd there is no element of
$W_{\delta_p}^0$  that can change sign to $\nu$ (only an even number
of sign changes can occur in the last $p$ entries of $\nu$).
 It
is not hard to check that the operator
$\mathcal{A}(-1,\delta_p,\nu)$ takes opposite sign on the two fine
$U(n)$-types ($\Lambda^p(\mathbb{C}_n)$ and
$\Lambda^p(\mathbb{C}_n)^\ast $) containing $\delta_p$. Both fine
$U(n)$-types are contained in $L(\delta_p,\nu)$, hence
$L(\delta_p,\nu)$ is nonunitary.

\end{enumerate}
In these examples    $R_{\delta}=\mathbb{Z}/{2\mathbb{Z}}$ and
$A(\delta)$ contains  two  fine $K$-types (which   correspond to the
trivial and the sign representation of
${\mathbb{Z}}/{2\mathbb{Z}}$ respectively). Because $[w] $ is a
generator for the $R_{\delta}$ group,  the operator
$\mathcal{A}(w,\delta,\nu)$ takes opposite sign on the two fine
$K$-types. Hence the Langlands quotient $L(\delta,\nu)$ (which
includes both  $K$-types) is nonunitary.

\begin{remark} In proposition \ref{prop 4.3}, $\nu$ is required to be dominant, so
$X(\delta,\nu)$ has a unique irreducible subquotient (which must of
course contain all the fine $K$-types for $\delta$). If $\nu$ is
only weakly dominant, the lowest $K$-types are distributed among the
various subquotients. If  one of these subquotients contains at
least two fine $K$-types, then a similar nonunitarity test may
apply.
\end{remark}

\begin{proposition}
Suppose that
\begin{enumerate}
\item  $R_{\delta} \neq \{1\} $, and
\item $R_{\delta}(\nu) \neq R_{\delta}$
\end{enumerate}
so that there is a subquotient of $X(\delta,\nu)$ containing (at
least) two fine $K$-types $\pi_1$ and $\pi_2$. Also assume that  $ w
\cdot \nu = - \bar{\nu}$ for some $w\in W_{\delta}$, and that there
is no element $w_0\in W_{\delta}^0$ with the same property.  If $
\varrho_{\pi_1}[w] \neq \varrho_{\pi_2}[w], $ then the Langlands
subquotient $L(\delta,\nu)(\pi_1)=L(\delta,\nu)(\pi_2)$ is not
unitary.
\end{proposition}

Petite $K$-types are meant to provide additional nonunitarity
certificates for  Langlands quotients of principal series.

\subsection{Definition of petite $K$-type}
The definition of petite $K$-types is rather technical. It is intended
to
provide a natural relation between intertwining operators for real
groups and intertwining operators for graded Hecke algebras. This
relation will be detailed in the next few sections.

\begin{comment}
Our definition of petite $K$-type will probably  appear too
technical and rather difficult
 to a first reader. It is, nonetheless, a ``natural''
definition to create a bridge between  intertwining operators for
real groups and intertwining operators for graded Hecke algebras.
This will become clear in the next few sections.
\end{comment}

\begin{definition}\label{d:4.4.1}
Let $\mu$ be a $K$-type containing $\delta$.
 We say that $\mu$ is
``petite for $\delta$'' if the following conditions hold.
\begin{enumerate}

\item If $\alpha$ is a good root for $\delta$,  and  $\alpha$ is
simple in  both $W_{\delta}^0$ and $W$, then the only possible
eigenvalues of $Z_{\alpha}^2$ on $\Hom_M(\mu,\delta)$ are $0$ and
$-4$. This gives: \[ \Hom_M(\mu,\delta) \equiv E^{\alpha}(0) \oplus
E^{\alpha}(-4).\]

\item  If $\alpha$ is a good root for $\delta$,   and  $\alpha$ is
simple in  $W_{\delta}^0$ but not in $W$,     choose
  a \emph{minimal}
decomposition of $s_{\alpha}$ in $W$ of the form \[s_{\alpha} =
(s_{\gamma_1}s_{\gamma_2}\cdots s_{\gamma_l} )s_{\xi} ( s_{\gamma_l}
\cdots s_{\gamma_2}  s_{\gamma_1}) \] with
\begin{itemize}

\item $\xi$ good for $ \tau=(s_{\gamma_{l}} \cdots
s_{\gamma_2} s_{\gamma_1}) \delta$, and

\item  $\gamma_j$   bad for both $\tau_j^-=
(s_{\gamma_{j-1}} \cdots s_{\gamma_2}  s_{\gamma_1}) \delta$ and
$\tau_j^+= (s_{\gamma_{j+1}} \cdots s_{\gamma_{l}} s_{\xi}
s_{\gamma_{l}} \cdots  s_{\gamma_1}) \delta$, for all  $ j=1\dots
l$.

\end{itemize}
Then
\begin{itemize}

\item[ (2a)] the only possible
eigenvalues of $Z_{\xi}^2$ on $\Hom_M(\mu,\tau)$ are $0$ and $-4$.

\item[ (2b)] the only possible
eigenvalue of $Z_{\gamma_j}^2$ on $\Hom_M(\mu,\tau_j^{\pm})$ is
$-1$.

\end{itemize}
Note that, in these cases,
\[ \Hom_M(\mu,\tau) \equiv E^{\xi}(0) \oplus E^{\xi}(-4)\]
and
\[ \Hom_M(\mu,\tau^{\pm}) \equiv E^{\gamma_j}(-1)\quad \forall \, j = 1 \dots l.\]

\item  If $\alpha$ is a bad root for  $\delta$ but $s_{\alpha}$ is in the stabilizer of $\delta$, then the only possible
eigenvalue of $Z_{\alpha}^2$ on $\Hom_M(\mu,\delta)$ is $-1$. This
gives:
\[ \Hom_M(\mu,\delta) \equiv E^{\alpha}(-1).\]
\end{enumerate}
\end{definition}

\begin{proposition}\label{p:4.4.1}
 $K$-types of level at
most 2 are petite.
\end{proposition}

\begin{definition}\label{d:4.4.2}
Let $\mu$ be a $K$-type containing $\delta$.
 We say that $\mu$ is
``quasi-petite for $\delta$'' if  conditions (1) and (2) hold.
\end{definition}

\begin{proposition}\label{p:4.4.2}
 $K$-types of level at most 3 are quasi-petite.
\end{proposition}

\subsection{Spherical Petite $K$-types}
Assume that $\delta$ is trivial and $\nu$ is weakly dominant. Choose
$w\in W$ such that $w\cdot Q=\bar{Q}$, $w\cdot \delta=\delta$ and
$w\cdot \nu = -\bar{\nu}$. Note that $W_{\delta}^0 = W_{\delta} = W$
in this case.

\begin{theorem}\label{t:4.5}
Let $\delta$ be the trivial character of $M$. For every spherical
$K$-type  $\mu$, let ${\psi_{\mu}^0}$ be the representation of
$W_{\delta}^0$ on $\Hom_M(\mu,\delta)$ introduced in section
\ref{sec 4.2}. If $\mu$ is petite, then
\[ \mathcal{A}_{\mu}(w,\delta,\nu) = \mathcal{A}_{({\psi_{\mu}^0})^\ast} ( w ,\nu)  \]
for all weakly dominant $\nu$.
\end{theorem}
The operator on the left is an operator for the  real group $G$; the
operator on the right is an operator for the affine graded Hecke
algebra corresponding to  $W$ (\cf section \ref{sec:3.1}).
\begin{proof}
Consider first the Hecke algebra operator
\begin{equation}
\mathcal{A}_{(\psi_{\mu}^0)^{\ast}}(w, \nu)\colon
\left(V_{(\psi_{\mu}^0)^{\ast}}\right)^{\ast}=\Hom_{M}(\mu, \delta)
\rightarrow
\left(V_{(\psi_{\mu}^0)^{\ast}}\right)^{\ast}=\Hom_{M}(\mu, \delta)
\end{equation}
(Note that domain and codomain of
$\mathcal{A}_{(\psi_{\mu}^0)^{\ast}}$ coincide with the ones of
 $\mathcal{A}_{\mu}(w,\delta,\nu)$). If
\begin{equation}
w = s_{\beta_r}\cdots s_{\beta_2} s_{\beta_1}
\end{equation}
is a minimal decomposition of $w$ in $W$, then
$\mathcal{A}_{(\psi_{\mu}^0)^{\ast}}$ has  a Gindikin-Karpelevi\v c
factorization of the form
\begin{equation}\label{eq:3a.5.3}
\mathcal{A}_{(\psi_{\mu}^0)^{\ast}}(w, \nu) = \prod\limits_{j=1}^r
\mathcal{A}_{(\psi_{\mu}^0)^{\ast}}(s_{\beta_j}, \gamma_{j-1})
\end{equation}
with $\gamma_{j-1} = s_{\beta_{j-1}}\cdots s_{\beta_2} s_{\beta_1}
\nu$ for all $j \geq 1$ and $\gamma_{0}=\nu$.  The $j^{th}$-factor
of the operator acts by
\begin{equation}
\mathcal{A}_{(\psi_{\mu}^0)^{\ast}}(s_{\beta_j},\gamma_{j-1})
:=\begin{cases} 1 & \, \text{on
the $(+1)$-eigenspace of ${\psi_{\mu}^0}(s_{\beta_j})$} \\
\frac{1-\langle \gamma_{j-1},\,\check{\beta_j} \rangle}{1+\langle
\gamma_{j-1},\,\check{\beta_j} \rangle} & \, \text{on the
$(-1)$-eigenspace of ${\psi_{\mu}^0}(s_{\beta_j})$}. \end{cases}
\end{equation}
\[\]
Here is a picture of the action of
$\mathcal{A}_{(\psi_{\mu}^0)^{\ast}}(s_{\beta_j},\gamma_{j-1})$:

\[
 \noindent  \hskip-4cm
 \begin{picture}(0,-40)(120,10)
\put(30,-7){$\mathcal{A}_{(\psi_{\mu}^0)^{\ast}}(s_{\beta_j},\gamma_{j-1}):$}
\put(140,10){\circle*{4}} \put(95,15){\tiny{$(+1)$-eigensp.~of
${\psi_{\mu}^0}(s_{\beta_j})$} } \put(140,10){\line(1,0){72}}
\put(210,10){\circle*{4}} \put(190,15){\tiny{$(-1)$-eigensp.~of
${\psi_{\mu}^0}(s_{\beta_j})$} } \put(140,3){\vector(0,-1){16}}
 \put(144,-7){\tiny{$  {\psi_{\mu}^0}(s_{\beta_j})$}}
\put(210,3){\vector(0,-1){16}} \put(213,-7){\tiny{$- \frac{1-
\langle \gamma_{j-1} ,\, \check{\beta_j} \rangle }{1+ \langle
\gamma_{j-1} ,\, \check{\beta_j}
\rangle}{\psi_{\mu}^0}(s_{\beta_j})$}}
 \put(140,-20){\circle*{4}}
 \put(95,-30){\tiny{$(+1)$-eigensp.~of
${\psi_{\mu}^0}(s_{\beta_j})$} } \put(140,-20){\line(1,0){72}}
\put(210,-20){\circle*{4}} \put(190,-30){\tiny{$(-1)$-eigensp.~of
${\psi_{\mu}^0}(s_{\beta_j})$} }
%\put(210,-20){\line(1,0){72}} \put(280,-20){\circle*{4}}
\end{picture}
\]
\[\]\[\]
\[\]
Notice that
\begin{equation}
\mathcal{A}_{(\psi_{\mu}^0)^{\ast}}(s_{\beta_j},\gamma_{j-1})\equiv\frac{1+\langle
\gamma_{j-1},\,\check{\beta_j} \rangle
{\psi_{\mu}^0}(s_{\beta_j})}{1+\langle
\gamma_{j-1},\,\check{\beta_j} \rangle}.
\end{equation}

Now look at the real operator. Corresponding to the same minimal
decompositions of $w$ in $W$, there is the
factorization (\ref{eq:3a.1.4}) for $\mathcal{A}_{\mu}(w,\delta,\nu)$:
\begin{equation}
\mathcal{A}_{\mu}(w,\delta,\nu) = \prod\limits_{j=1}^r
\mathcal{A}_{\mu}(s_{\beta_j},\rho_{j-1},\nu_{j-1}).
\end{equation}
Observe that $ \rho_{j-1} = \delta$ and $\nu_{j-1}= \gamma_{j-1}$,
for all $j=1\dots r$. Most importantly,
 both $\mathcal{A}_{(\psi_{\mu}^0)^{\ast}}(s_{\beta_j}, \gamma_{j-1})$ and
$\mathcal{A}_{\mu}(s_{\beta_j},\delta,\gamma_{j-1})$ are operators
on $\Hom_M(\mu,\delta)$. We will prove that the two factors match.

Recall from section \ref{sec 4.1} that the action of
$\mathcal{A}_{\mu}(s_{\beta_j},\delta,\gamma_{j-1})$ on
$\Hom_M(\mu,\delta)$ depends on the decomposition of
 $\Hom_M(\mu,\delta)$ into $Z_{\beta_j}^2$ eigenspaces. In particular,  $\mathcal{A}_{\mu}(s_{\beta_j},\delta,\gamma_{j-1})$
 acts on $E^{\beta_j}(-4m^2) $ by
\begin{equation} T \mapsto
c_{2m}(\beta_j, \gamma_{j-1}) \,\mu_{\delta}(\sigma_{\beta_j}) \, T
\, \mu(\sigma_{\beta_j}^{-1}).
\end{equation}
Because $\beta_j$ is a good root for $\delta$, we can re-write this
action as
\begin{equation}
T\mapsto c_{2m}(\beta_j, \gamma_{j-1})\,
{\psi_{\mu}^0}(s_{\beta_j})T.
\end{equation}
 Now, because $\mu$ is petite, and $\beta_j$ is a good
root for $\delta$ (simple in both $W$ and $W_{\delta}^0$), the only
eigenvalues of $ Z_{\beta_j}^2$ on $\Hom_M(\mu,\delta)$ are $0$ and
$-4$:
\[
\Hom_M(\mu,\delta) :=  E^{\beta_j}(0)\oplus E^{\beta_j}(-4).
\]
Hence only the constants $c_{0}(\beta_j,\gamma_{j-1})=1$ and $
c_{2}(\beta_j,\gamma_{j-1})=- \frac{1- \langle \gamma_{j-1} ,\,
\check{\beta_j} \rangle }{1+ \langle \gamma_{j-1} ,\,
\check{\beta_j} \rangle}$   appear:

 \noindent\hskip-.5cm \begin{picture}(0,30)(45,10)
 \put(140,10){\circle*{4}}  \put(130,16){\tiny{$E\left(0 \right)$}}
\put(140,10){\line(1,0){72}} \put(205,10){\circle*{4}}
\put(195,16){\tiny{$E\left(-4 \right)$}}
\put(210,10){\line(1,0){72}} \put(285,10){\circle*{4}}
\put(275,16){\tiny{$E\left(- 16 \right)$}}
\put(280,10){\line(1,0){72}} \put(350,10){\circle*{4}}
\put(355,10){\ldots\ldots \ldots\ldots }
%$\put(415,10){\circle*{4}}$
\put(345,16){\tiny{$E\left(-36 \right)$}} \put(355,-20){\ldots\ldots
\ldots\ldots }
%$\put(415,-20){\circle*{4}}$
 \put(345,-30){\tiny{$E\left(-36 \right)$}}
%\put(360,10){\ldots \ldots}%
\put(140,0){\vector(0,-1){12}}
 \put(143,-7){\tiny{$ {\psi_{\mu}^0}(s_{\beta_j})$}}
\put(205,0){\vector(0,-1){12}} \put(205,-7){\tiny{$-\frac{1- \langle
\gamma_{j-1}, \,\check{\beta_j} \rangle}{1+ \langle \gamma_{j-1},
\,\check{\beta_j} \rangle} {\psi_{\mu}^0}(s_{\beta_j})$}}

%\put(285,-7){\tiny{$c_{4}(\beta_j,\gamma_{j-1})
%{\psi_{\mu}^0}(s_{\beta_j})$}}
%\put(354,-7){\tiny{$c_{6}(\beta_j,\gamma_{j-1})
%{\psi_{\mu}^0}(s_{\beta_j})$}}

\put(285,0){\vector(0,-1){12}}
 \put(350,0){\vector(0,-1){12}}
 \put(140,-20){\circle*{4}}
\put(130,-30){\tiny{$E\left(0 \right)$}}
\put(140,-20){\line(1,0){72}} \put(205,-20){\circle*{4}}
\put(196,-30){\tiny{$E\left(-4 \right)$}}
\put(210,-20){\line(1,0){72}} \put(285,-20){\circle*{4}}
\put(275,-30){\tiny{$E\left(-16 \right)$}}
\put(280,-20){\line(1,0){72}} \put(350,-20){\circle*{4}}
  %\put(250,5){\line(1,1){20}}
 \put(275,13){\line(1,1){10}} \put(283,-2){\line(1,1){28}}
     \put(287,-25){\line(1,1){50}} \put(317,-25){\line(1,1){50}}
   \put(347,-25){\line(1,1){50}}   \put(377,-25){\line(1,1){30}}

\put(60,-7){$\mathcal{A}_{\mu}(s_{\beta_j},\delta,\gamma_{j-1}):$}

\end{picture}
\vskip.65in

Notice that \begin{align} E^{\beta_j}(0) &= (+1)\text{-eigenspace of
${\psi_{\mu}^0}(s_{\beta_j})$}
 \\
  E^{\beta_j}(-4)
& = (-1)\text{-eigenspace of ${\psi_{\mu}^0}(s_{\beta_j})$}
\end{align} so
\begin{equation}
\mathcal{A}_{\mu}(s_{\beta_j},\delta,\gamma_{j-1}) =
\mathcal{A}_{(\psi_{\mu}^0)^{\ast}}(s_{\beta_j}, \gamma_{j-1})\quad
\forall  j=1\dots r.
\end{equation}
We conclude  that
\begin{equation}
\mathcal{A}_{\mu}(w,\delta,\nu) =
\mathcal{A}_{(\psi_{\mu}^0)^{\ast}}(w, \nu).
\end{equation}
\end{proof}

\subsection{Nonspherical petite \boldmath{$K$}-types.}\label{sec 4.6} The matching
of operators in the \emph{nonspherical} case is much harder,
especially if the $K$-type is not pseudospherical. In this case,
$W_{\delta}^0$ is a proper subgroup of $W$ and the two intertwining
operators cannot be compared term by term. We point out some of the
difficulties:
\begin{itemize}

\item The element $w$ is in the stabilizer of $\delta$ in $W$, but does not
necessarily belong to $W_{\delta}^0$, so the operator
$\mathcal{A}_{(\psi_{\mu}^0)^{\ast}}(w, \nu)$ might be meaningless.

\item If $\alpha$ is a  bad root for $\delta$, and is simple in $W$, then the $\alpha$-factor of
$\mathcal{A}_{\mu}(w,\delta,\nu)$ does not have an immediate
correspondent in $\mathcal{A}_{(\psi_{\mu}^0)^{\ast}}(w,\nu)$.

\item If $\alpha$ is  a  good  root for $\delta$, and is simple   in $W_{\delta}^0$ but  not
in $W$, then   the $\alpha$-factor of
$\mathcal{A}_{(\psi_{\mu}^0)^{\ast}}(w,\nu)$ does not have an
immediate correspondent in $\mathcal{A}_{\mu}(w,\delta,\nu)$.

\end{itemize}

It is convenient to proceed by increasing level of difficulty.  At
this stage we will assume that  $w\in W_{\delta}^0$, so that we can
choose a minimal decomposition of $w$ in $W_{\delta}^0$ which is
``compatible'' with the one of $w$ in $W.$

\begin{lemma}\label{lem:4.6}
Let  $w = s_{\beta_r} \cdots s_{\beta_2} s_{\beta_1}$ be a minimal
decomposition of $w$ in $W_{\delta}^0$. Notice that the   roots
occurring in this factorization are obviously simple in
$W_{\delta}^0$, but need \emph{not} to be simple in $W$.  If
$\beta_i$ is not simple,  we can choose
  a \emph{minimal}
decomposition of $s_{\beta_i}$ in $W$ such that:
\begin{equation}\label{e:4.6.1}
s_{\beta_i} = (s_{\gamma_1}s_{\gamma_2}\cdots s_{\gamma_l}) s_{\xi}
( s_{\gamma_l} \cdots s_{\gamma_2}  s_{\gamma_1} )
\end{equation}
 with
\begin{itemize}
\item $\xi$ good for $ (s_{\gamma_{l}} \cdots
s_{\gamma_2} s_{\gamma_1} )\delta, and$

\item $\gamma_j$ bad for both $
(s_{\gamma_{j-1}} \cdots s_{\gamma_2} s_{\gamma_1}) \delta$ and $
(s_{\gamma_{j+1}} \cdots s_{\gamma_{l}} s_{\xi} s_{\gamma_{l}}
\cdots  s_{\gamma_1}) \delta.$

\end{itemize}
We can choose these minimal decompositions in a way that,  after
replacing every non simple reflection $ s_{\beta_i}$ by its
expression in \ref{e:4.6.1}, we obtain a minimal decomposition of
$w$ in $W$.

\end{lemma}

The following example will clarify the lemma. Let $G$ be the split
linear group of type $F_4$, with simple roots:
\begin{center}\end{center}
 {\footnotesize{
 \noindent\hskip-.5cm \begin{picture}(0,-40)(50,10)
 \put(140,10){\circle*{5}}  \put(90,20){$\alpha_1=\epsilon_1 - \epsilon_2 - \epsilon_3 - \epsilon_4$}
\put(140,10){\line(1,0){72}} \put(210,10){\circle*{5}}
\put(195,20){$\alpha_2=2\epsilon_4$} \put(211,12){\line(1,0){70}}
\put(211,8){\line(1,0){70}}
 \put(280,10){\circle*{5}}
\put(264,20){$\alpha_3=\epsilon_3-\epsilon_4$}
\put(280,10){\line(1,0){72}} \put(350,10){\circle*{5}}
\put(335,20){$\alpha_4=\epsilon_2 - \epsilon_3$}
\end{picture}}}
\begin{center}\end{center}
 There is a one-dimensional nongenuine character $\delta$ of $M$
(with a three-dimensional orbit $\{\delta,\delta',\delta''\}$ under
$W$) that admits
\begin{center}\end{center}
 {\footnotesize{
 \noindent\hskip-.5cm \begin{picture}(0,-40)(50,10)
 \put(140,10){\circle*{5}}
\put(210,10){\circle*{5}}
 \put(280,10){\circle*{5}}
 \put(350,10){\circle*{5}}
 \put(140,8){\line(1,0){72}}
  \put(140,12){\line(1,0){72}}
 \put(211,10){\line(1,0){70}}
 \put(280,10){\line(1,0){70}}
 \put(196,20){$\alpha_3=\epsilon_3-\epsilon_4$}
\put(125,20){$\alpha_2= 2\epsilon_4 $}
\put(265,20){$\alpha_4=\epsilon_2 - \epsilon_3$}
\put(339,20){$\beta=\epsilon_1 - \epsilon_2$}
\end{picture}}}
\[\]
as a basis for the good roots. Then $W_{\delta}^0$ is a Weyl group
of type $C_4.$ The long Weyl group element $w=-1$ has length 16 in
$W_{\delta}^0$, and length 24 in $W$. We choose
\begin{equation}\label{e:4.6.2}
 w =
s_{\alpha_2}(
s_{\alpha_3}s_{\alpha_2}s_{\alpha_3})(s_{\alpha_4}s_{\alpha_3}s_{\alpha_2}s_{\alpha_3}s_{\alpha_4})
(s_{\beta}s_{\alpha_4}s_{\alpha_3}s_{\alpha_2}s_{\alpha_3}s_{\alpha_4}s_{\beta})
\end{equation}
to be a  minimal decomposition of $w$ in $W_{\delta}^0$. Note that
$\beta$ is the only good root which is simple in $W_{\delta}^0$ but
not in $W$; we decompose it as:
\begin{equation}\label{e:4.6.3}
s_{\beta} = (s_{\alpha_1}
s_{\alpha_2})s_{\alpha_3}(s_{\alpha_2}s_{\alpha_1}).
\end{equation}
Then
\begin{itemize}
\item $\alpha_1$ is bad for  $\delta$. The reflection $s_{\alpha_1}$
 carries $\delta$ into $\delta'=s_{\alpha_1}\delta$  and  $\delta'$
 into $\delta$.

\item  $\alpha_2$ is bad for   $\delta'$. The reflection $s_{\alpha_2}$
 carries $\delta'$ into $\delta''=s_{\alpha_2}\delta'
 $  and $\delta''$ into $\delta'$.

\item $\alpha_3$ is good for  $\delta''$. . The reflection $s_{\alpha_3}$ obviously stabilizes
$\delta''$.

\item  The composition $s_{\alpha_1}
s_{\alpha_2}s_{\alpha_3}s_{\alpha_2}s_{\alpha_1}$ stabilizes
$\delta$.
\end{itemize}

 Replacing every occurrence  of  $s_{\beta}$ in $w$ by the
product $  s_{\alpha_1}
s_{\alpha_2}s_{\alpha_3}s_{\alpha_2}s_{\alpha_1}$, we obtain a
minimal decomposition of $w$ in $W$:

$ \quad w = s_{\alpha_2}(
s_{\alpha_3}s_{\alpha_2}s_{\alpha_3})(s_{\alpha_4}s_{\alpha_3}s_{\alpha_2}s_{\alpha_3}s_{\alpha_4})
\cdot $\\

$\quad  \quad  \quad (\underbrace{
s_{\alpha_1}s_{\alpha_2}s_{\alpha_3}s_{\alpha_2}s_{\alpha_1}}_{=s_{\beta}}s_{\alpha_4}s_{\alpha_3}s_{\alpha_2}s_{\alpha_3}s_{\alpha_4}
\underbrace{
s_{\alpha_1}s_{\alpha_2}s_{\alpha_3}s_{\alpha_2}s_{\alpha_1}}_{=s_{\beta}}).$

\begin{theorem}\label{t:4.6}
Let $\delta$ be a nontrivial representation of $M$ and  let $\nu$ be
a weakly dominant character of $A$.  Assume that  there exists $w\in
W_{\delta}^0$ such that $w\cdot Q=\bar{Q}$  and  $w\cdot \nu =
-\bar{\nu}$. For every
 $K$-type  $\mu$ containing $\delta$, let ${\psi_{\mu}^0}$ be
the representation of $W_{\delta}^0$ on $\Hom_M(\mu,\delta)$
introduced in section \ref{sec 4.2}.  If $\mu$ is petite, then
\[ \mathcal{A}_{\mu}(w,\delta,\nu) = \mathcal{A}_{({\psi_{\mu}^0})^\ast} ( w ,\nu).  \]
\end{theorem}
The operator on the left is an operator for the  real group $G$; the
operator on the right is an operator for the affine graded Hecke
algebra corresponding to $W_{\delta}^{0}$ (\cf section
\ref{sec:3.1}).

\begin{proof}
Suppose that \begin{equation}w = s_{\beta_r} \cdots s_{\beta_2}
s_{\beta_1} \end{equation} and
 \begin{equation}
w = s_{\beta_r}\cdots (\underbrace{ s_{\gamma_1}s_{\gamma_2}\cdots
s_{\gamma_l} s_{\xi}   s_{\gamma_l} \cdots s_{\gamma_2} s_{\gamma_1}
}_{s_{\beta_i}})s_{\beta_{i-1}} \cdots (\underbrace{
s_{\tau_1}s_{\tau_2}\cdots s_{\tau_l} s_{\zeta}   s_{\tau_l} \cdots
s_{\tau_2} s_{\tau_1} }_{s_{\beta_j}})\cdots s_{\beta_1}
\end{equation}
are minimal decompositions of $w$ in $W_{\delta}^0$ and $W$
respectively, that are compatible in the sense of the previous
lemma.

Recall that the Hecke algebra operator
\begin{equation}
\mathcal{A}_{(\psi_{\mu}^0)^{\ast}}(w, \nu)\colon
\left(V_{(\psi_{\mu}^0)^{\ast}}\right)^{\ast}=\Hom_{M}(\mu, \delta)
\rightarrow
\left(V_{(\psi_{\mu}^0)^{\ast}}\right)^{\ast}=\Hom_{M}(\mu, \delta)
\end{equation}
 has the factorization (\ref{eq:3a.5.3}):
%a Gindikin-Karpelevi\v c
%factorization, corresponding to the minimal decomposition of $w$ in
%$W_{\delta}^0$:
\begin{equation}
\mathcal{A}_{(\psi_{\mu}^0)^{\ast}}(w, \nu) = \prod\limits_{j=1}^r
\mathcal{A}_{(\psi_{\mu}^0)^{\ast}}(s_{\beta_j}, \gamma_{j-1}).
 \end{equation}
%Here $\gamma_{j-1} = s_{\beta_{j-1}}\cdots s_{\beta_2} s_{\beta_1}
%\nu$ for all $j \geq 1$ and $\gamma_{0}=\nu$.\\
  The $j^{th}$-factor
of $\mathcal{A}_{(\psi_{\mu}^0)^{\ast}}(w, \nu)$ is again an
operator on $\Hom_{M}(\mu, \delta)$, and acts by \begin{equation}
\mathcal{A}_{(\psi_{\mu}^0)^{\ast}}(s_{\beta_j},\gamma_{j-1})
:=\begin{cases} 1 & \, \text{on
the $(+1)$-eigenspace of ${\psi_{\mu}^0}(s_{\beta_j})$} \\
\frac{1-\langle \gamma_{j-1},\,\check{\beta_j} \rangle}{1+\langle
\gamma_{j-1},\,\check{\beta_j} \rangle} & \, \text{on the
$(-1)$-eigenspace of ${\psi_{\mu}^0}(s_{\beta_j})$}. \end{cases}
\end{equation}

The real operator $\mathcal{A}_{\mu}(w,\delta,\nu)$, on the other
hand, has the factorization (\ref{eq:3a.5.3}).
 This factorization
involves  more factors; every factor corresponds to a simple root in
$W$ (which might, of course, be either good or bad for $\delta$). We
will prove that:
\begin{itemize}
\item[(a)] If $\beta$
is a good root  for $\delta$, and $\beta$ is simple in both
$W_{\delta}^0$ and $W$, then the $\beta$-factor of  $
\mathcal{A}_{(\psi_{\mu}^0)^{\ast}}(w, \nu)$ matches the
corresponding $\beta$-factor of $\mathcal{A}_{\mu}(w,\delta,\nu)$.

\item[(b)]  If $\beta$
is a good root   for $\delta$, and $\beta$ is simple   in
$W_{\delta}^0$ but not in $W$, write
\[
s_{\beta} = (s_{\gamma_1}s_{\gamma_2}\cdots s_{\gamma_l}) s_{\xi} (
s_{\gamma_l} \cdots s_{\gamma_2}  s_{\gamma_1})
\]
for a minimal decomposition of $s_{\beta}$ in $W$ as in lemma
\ref{lem:4.6}. Then the $\beta$-factor of $
\mathcal{A}_{(\psi_{\mu}^0)^{\ast}}(w, \nu)$ matches  the product of
the all the factors of $\mathcal{A}_{\mu}(w,\delta,\nu)$ coming from
$(s_{\gamma_1}s_{\gamma_2}\cdots s_{\gamma_l}) s_{\xi} (s_{\gamma_l}
\cdots s_{\gamma_2}  s_{\gamma_1})$.

\end{itemize}
In  our $F_4$ example, the matchings are as follows:
\begin{itemize}
\item[(a')]

$\mathcal{A}_{(\psi_{\mu}^0)^{\ast}}(\alpha_i,\gamma) =
\mathcal{A}_{\mu}(\alpha_i,\delta,\gamma) $ for all $i=2,3,4$, and all $\gamma$.\\

\item[(b')]

$\mathcal{A}_{(\psi_{\mu}^0)^{\ast}}(\beta,\gamma)=\mathcal{A}_{\mu}(\alpha_1,\delta',s_{\alpha_2}s_{\alpha_3}s_{\alpha_2}s_{\alpha_1}\gamma)
\circ
\mathcal{A}_{\mu}(\alpha_2,\delta'',s_{\alpha_3}s_{\alpha_2}s_{\alpha_1}\gamma)
\circ $ \\

$\qquad \qquad \circ
\mathcal{A}_{\mu}(\alpha_3,\delta'',s_{\alpha_2}s_{\alpha_1}\gamma)\circ
\mathcal{A}_{\mu}(\alpha_2,\delta',s_{\alpha_1}\gamma)\circ
\mathcal{A}_{\mu}(\alpha_1,\delta,\gamma) $\\
for all $\gamma.$
\end{itemize}
Instead of proving  claims $(a)$ and  $(b)$, we will prove $(a')$
and  $(b')$ instead. The general idea will emerge from this simpler
case.

 Condition $(a')$ is easy to
prove. Let $\beta=\alpha_i$, for $i=2,3,4$. Then  $\beta$ is a good
root for $\delta$, and is simple in both $W$ and $W_{\delta}^0$.
Because $\mu$ is petite,
 $\Hom_M(\mu,\delta)= E^{\beta}(0) \oplus E^{\beta}(-4)$.   The
same argument used in the spherical case shows that the
$\beta$-factor of the Hecke algebra operator matches the
$\beta$-factor  of the real operator:
\begin{equation}
\mathcal{A}_{(\psi_{\mu}^0)^{\ast}}(\beta,\gamma)=
\mathcal{A}_{\mu}(\beta,\delta,\gamma) \end{equation} for all
$\gamma$. Notice that both factors act on $\Hom_M(\mu,\delta)$.

 Condition $(b')$ is more delicate. Choose $\beta= \alpha_1$ and
 $s_{\beta} =
 (s_{\alpha_1}s_{\alpha_2}s_{\alpha_3}s_{\alpha_2}s_{\alpha_1})$.
The Hecke algebra operator
$\mathcal{A}_{(\psi_{\mu}^0)^{\ast}}(\beta,\gamma)$ acts on
$\Hom_M(\mu,\delta)$ in the usual way. The composition of the
factors of the real operator corresponding to $
(s_{\alpha_1}s_{\alpha_2}s_{\alpha_3}s_{\alpha_2}s_{\alpha_1})$ also
acts on $\Hom_M(\mu,\delta)$; the single factors, however, do not.
Set $\delta'=s_{\alpha_1}\delta$ and
$\delta''=s_{\alpha_2}s_{\alpha_1}\delta.$ Then
\begin{align}
 \mathcal{A}_{\mu}(s_{\alpha_1},\delta,\gamma)&\colon \Hom_M(\mu,\delta) \,\,\,\longrightarrow  \Hom_M(\mu,\delta')\\
\mathcal{A}_{\mu}(s_{\alpha_2},\delta',s_{\alpha_1}\gamma)&\colon
\Hom_M(\mu,\delta')
\,\longrightarrow  \Hom_M(\mu,\delta'')\\
\mathcal{A}_{\mu}(s_{\alpha_3},\delta'',s_{\alpha_2}s_{\alpha_1}\gamma)&\colon
\Hom_M(\mu,\delta'') \longrightarrow  \Hom_M(\mu,\delta'')\\
\mathcal{A}_{\mu}(s_{\alpha_2},\delta'',s_{\alpha_3}s_{\alpha_2}s_{\alpha_1}\gamma)&\colon
\Hom_M(\mu,\delta'') \longrightarrow  \Hom_M(\mu,\delta')\\
 \mathcal{A}_{\mu}(s_{\alpha_1},\delta',s_{\alpha_2}s_{\alpha_3}s_{\alpha_2}s_{\alpha_1}\gamma)&\colon \Hom_M(\mu,\delta')  \,\,\longrightarrow
 \Hom_M(\mu,\delta).
\end{align}
 Because $\mu$ is petite and $\alpha_1$ is bad for both
$\delta$ and $\delta'$, serious restrictions are imposed on the
eigenvalues of $\text{d}\mu(Z_{\alpha_1})^2$ on both
$\Hom_M(\mu,\delta)$ and
$\Hom_M(\mu,\delta')$:\\
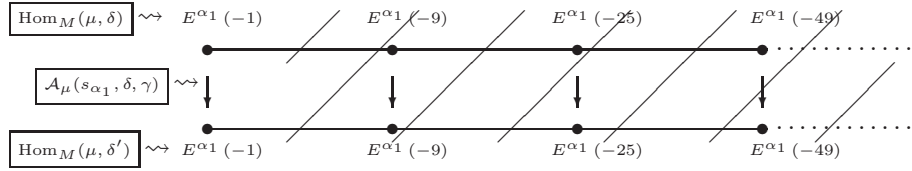
\[\]
 \noindent\hskip-.5cm
 \begin{picture}(0,-40)(50,10)

 \put(65,20){$\tiny{\boxed{\Hom_M(\mu,\delta)}}$\,$\rightsquigarrow $}
\put(65,-30){$\tiny{\boxed{\Hom_M(\mu,\delta')}}$\,$\rightsquigarrow
$}
\put(75,-5){$\tiny{\boxed{\mathcal{A}_{\mu}(s_{\alpha_1},\delta,\gamma)}}$\,$\rightsquigarrow
$}

 \put(140,10){\circle*{4}}  \put(130,20){\tiny{$E^{\alpha_1}\left(-1\right)$}}
\put(140,10){\line(1,0){72}} \put(210,10){\circle*{4}}
\put(200,20){\tiny{$E^{\alpha_1}\left(-9\right)$}}
\put(210,10){\line(1,0){72}} \put(280,10){\circle*{4}}
\put(270,20){\tiny{$E^{\alpha_1}\left(-25\right)$}}
\put(280,10){\line(1,0){72}} \put(350,10){\circle*{4}}
\put(355,10){\ldots\ldots \ldots\ldots }
%$\put(415,10){\circle*{4}}$
\put(345,20){\tiny{$E^{\alpha_1}\left(-49\right)$}}
\put(355,-20){\ldots\ldots \ldots\ldots }
%$\put(415,-20){\circle*{4}}$
 \put(345,-30){\tiny{$E^{\alpha_1}\left(-49\right)$}}
%\put(360,10){\ldots \ldots}%
\put(140,0){\vector(0,-1){12}}
 %\put(145,-7){\tiny{$1\cdot s_{\alpha_i}$}}
\put(210,0){\vector(0,-1){12}} \put(280,0){\vector(0,-1){12}}
 \put(350,0){\vector(0,-1){12}}
 \put(140,-20){\circle*{4}}
\put(130,-30){\tiny{$E^{\alpha_1}\left(-1\right)$}}
\put(140,-20){\line(1,0){72}} \put(210,-20){\circle*{4}}
\put(200,-30){\tiny{$E^{\alpha_1}\left(-9\right)$}}
\put(210,-20){\line(1,0){72}} \put(280,-20){\circle*{4}}
\put(270,-30){\tiny{$E^{\alpha_1}\left(-25\right)$}}
\put(280,-20){\line(1,0){72}} \put(350,-20){\circle*{4}}
  \put(170,5){\line(1,1){20}}  \put(170,-25){\line(1,1){50}}
     \put(210,-25){\line(1,1){50}}  \put(250,-25){\line(1,1){50}}  \put(290,-25){\line(1,1){50}}
   \put(330,-25){\line(1,1){50}}  \put(370,-25){\line(1,1){30}}
\end{picture}
\vskip.75in

 Then, for all $\gamma$, we have:
\begin{equation}
\mathcal{A}_{\mu}(s_{\alpha_1},\delta,\gamma)\colon
\Hom_M(\mu,\delta) \longrightarrow  \Hom_M(\mu,\delta'),\, T \mapsto
\mu_{\delta}(\sigma_{\alpha_1}) T \mu(\sigma_{\alpha_1})^{-1}
\end{equation}
and
\[ \mathcal{A}_{\mu}(s_{\alpha_1},\delta',s_{\alpha_2}s_{\alpha_3}s_{\alpha_2}s_{\alpha_1}\gamma)\colon  \Hom_M(\mu,\delta') \longrightarrow  \Hom_M(\mu,\delta),\, T \mapsto
\mu_{\delta}(\sigma_{\alpha_1}) T \mu(\sigma_{\alpha_1})^{-1}.
\]
Note that, because $T$ is  $M$-invariant and $\sigma_{\alpha_1}^2\in
M$, we can also write
\begin{equation}
 \mathcal{A}_{\mu}(s_{\alpha_1},\delta',s_{\alpha_2}s_{\alpha_3}s_{\alpha_2}s_{\alpha_1}\gamma)T=\mu_{\delta}(\sigma_{\alpha_1})^{-1} T
 \mu(\sigma_{\alpha_1}).
\end{equation}
There are similar restrictions on the eigenvalues of
$\text{d}\mu(Z_{\alpha_2})^2$  on $\Hom_M(\mu,\delta')$ and
$\Hom_M(\mu,\delta'')$. Hence
\begin{equation} \mathcal{A}_{\mu}(s_{\alpha_2},\delta',s_{\alpha_1}\gamma)\colon  \Hom_M(\mu,\delta') \longrightarrow  \Hom_M(\mu,\delta''),\, T \mapsto
\mu_{\delta}(\sigma_{\alpha_2}) T \mu(\sigma_{\alpha_2})^{-1}
\end{equation}
and
\begin{equation} \mathcal{A}_{\mu}(s_{\alpha_2},\delta'',s_{\alpha_3}s_{\alpha_2}s_{\alpha_1}\gamma)\colon  \Hom_M(\mu,\delta'') \longrightarrow  \Hom_M(\mu,\delta'),\, T \mapsto
\mu_{\delta}(\sigma_{\alpha_2})^{-1} T \mu(\sigma_{\alpha_2}).
\end{equation}
 We are only missing the central factor, $\mathcal{A}_{\mu}(s_{\alpha_3},\delta'',s_{\alpha_2}s_{\alpha_1}\gamma)$.
 Because $\mu$ is petite, and $\alpha_3$ is  good for $\delta''$,
 $Z_{\alpha_3}^2$ acts on $\Hom_M(\mu,\delta'')$ with eigenvalues $0$ and $(-4)$.
 Then
\begin{equation}
\mathcal{A}_{\mu}(s_{\alpha_3},\delta'',s_{\alpha_2}s_{\alpha_1}\gamma)\colon
\Hom_M(\mu,\delta'')\longrightarrow \Hom_M(\mu,\delta'')
\end{equation}
acts by
\begin{equation}
 T \mapsto
\frac{T + \langle
s_{\alpha_2}s_{\alpha_1}\gamma,\check{\alpha_3}\rangle\mu_{\delta}(\sigma_{\alpha_3})
T \mu(\sigma_{\alpha_3})^{-1}   }{1+ \langle
s_{\alpha_2}s_{\alpha_1}\gamma,\check{\alpha_3}\rangle}.
\end{equation}
Finally, we look at the composition of the 5 factors:\\

$\mathcal{A}_{\mu}(\alpha_1,\delta',s_{\alpha_2}s_{\alpha_3}s_{\alpha_2}s_{\alpha_1}\gamma)
\circ
\mathcal{A}_{\mu}(\alpha_2,\delta'',s_{\alpha_3}s_{\alpha_2}s_{\alpha_1}\gamma)
\circ $ \\

$\qquad \qquad \circ
\mathcal{A}_{\mu}(\alpha_3,\delta'',s_{\alpha_2}s_{\alpha_1}\gamma)\circ
\mathcal{A}_{\mu}(\alpha_2,\delta',s_{\alpha_1}\gamma)\circ
\mathcal{A}_{\mu}(\alpha_1,\delta,\gamma), $\\
\\
which acts by
\begin{equation}T\mapsto \frac{T + \langle
s_{\alpha_2}s_{\alpha_1}\gamma,\check{\alpha_3}\rangle\,\mu_{\delta}(\sigma_{\beta})
T \mu(\sigma_{\beta})^{-1}   }{1+ \langle s_{\alpha_2}s_{\alpha_1}
\gamma,\check{\alpha_3}\rangle}.
\end{equation} Because $\beta$ is a
good root for $\delta$ and
\[
 \langle s_{\alpha_2}s_{\alpha_1} \gamma,\check{\alpha_3}\rangle =
\langle \gamma,s_{\alpha_1}s_{\alpha_2}\check{\alpha_3}\rangle =
\langle \gamma, \check{\beta}\rangle,
\]
we can re-write this action as \[ T \mapsto \frac{T + \langle
\gamma, \check{\beta}\rangle\,{\psi_{\mu}^0}(s_{\beta}) T }{1+
\langle \gamma, \check{\beta}\rangle}. \] Hence  the product of
operators  behaves exactly like
$\mathcal{A}_{(\psi_{\mu}^0)^{\ast}}(\beta,\gamma)$. This concludes
the proof of the theorem.
\end{proof}

We now go one step forward, and discuss the case $w \in W^{\delta}$
(not necessarily in $W_{\delta}^0$).
 Recall that $W$ is the semidirect product of
$W_{\delta}^0$ with
\[
R_{\delta}^c=  \{w \in W \colon w({}^{\vee}\Delta_{\delta}^{+})=
{}^{\vee}\Delta_{\delta}^{+}  \}. \]
 $R_{\delta}^{c} $ is an abelian two
group, and is generated by reflections through strongly orthogonal
bad roots perpendicular to $\rho(\Delta_{\delta})$.

\begin{theorem}\label{t:4.7}
Let $\delta$ be a (nontrivial) representation of $M$ and  let $\nu$
be a weakly dominant character of $A$.  Assume  that there exists
$w\in W_{\delta}$ such that $w\cdot Q=\bar{Q}$ and  $w\cdot \nu =
-\bar{\nu}$. Write $w=u w^0$, with $u\in R_{\delta}^c$ and $w^0 \in
W_{\delta}^0.$
 For every
 $K$-type  $\mu$ containing $\delta$, let ${\psi_{\mu}}$ be
the representation of $W_{\delta}$ on $\Hom_M(\mu,\delta)$
introduced in section \ref{sec 4.2}, and let  ${\psi_{\mu}^0}$ be
its restriction to  $W_{\delta}^0$. If $\mu$ is petite, then
\[ \mathcal{A}_{\mu}(w,\delta,\nu) =\psi_{\mu}(u)\mathcal{A}_{(\psi_{\mu}^0)\ast}(w^0,\nu).  \]
\end{theorem}
The operator on the left is an operator for the  real group $G$; the
operator on the right is an operator for  the extended  Hecke
algebra corresponding to $W_{\delta}$ (\cf~section \ref{sec:3.7}).

\begin{proof}
Choose a minimal decomposition of $w=uw^0$ in $W$ of the form
\begin{equation}
w= \underbrace{s_{\zeta_t}\cdots  s_{\zeta_1}}_{u}
\underbrace{s_{\alpha_s}\cdots s_{\alpha_1}}_{w^0}
 \end{equation}
with $s_{\zeta_i}\in R_{\delta}^c$ and $(s_{\alpha_s}\cdots
s_{\alpha_1})$ a minimal decomposition of $w^0$ in $W$. Notice that
$\zeta_1\dots \zeta_t$ are bad roots for $\delta$, but the
corresponding reflections stabilize $\delta.$

The intertwining operator $\mathcal{A}_{\mu}(w,\delta,\nu)$ factors:
\begin{equation}
\mathcal{A}_{\mu}(w,\delta,\nu) =
\mathcal{A}_{\mu}(s_{\zeta_t},\delta,\nu_{s+t-1}) \cdots
\mathcal{A}_{\mu}(s_{\zeta_1},\delta,\nu_{s})
\left[\prod\limits_{j=1}^s
\mathcal{A}_{\mu}(s_{\alpha_j},\rho_{j-1},\nu_{j-1})\right].
\end{equation}

If $\mu$ is petite, then
\begin{equation}
\left[\prod\limits_{j=1}^s
\mathcal{A}_{\mu}(s_{\alpha_j},\rho_{j-1},\nu_{j-1})\right]
=\mathcal{A}_{(\psi_{\mu}^0)\ast}(w^0,\nu)
\end{equation}
by   theorem \ref{t:4.6}.

 Let us look at the remaining factors of the operator
 $\mathcal{A}_{\mu}(w,\delta,\nu)$.
Since $s_{\zeta_i}$ stabilizes $\delta$, each
$\mathcal{A}_{\mu}(s_{\zeta_i},\delta,\nu_{r+i-1})$ is an operator
on $\Hom_M(\mu,\delta)$. Because $\mu$ is petite,  $Z_{\gamma_i}^2$
acts on $\Hom_M(\mu,\delta)$ with eigenvalue $(-1)$.
 Then
\begin{equation}
 \mathcal{A}_{\mu}(s_{\zeta_i},\delta,\nu_{s+i-1}) \colon \Hom_M(\mu,\delta) \longrightarrow  \Hom_M(\mu,\delta),\, T \mapsto
\mu_{\delta}(\sigma_{\zeta_i}) T \mu(\sigma_{\zeta_i})^{-1}
\end{equation}
for all $i=1\dots t.$  We can re-write this action in terms of the
$W_{\delta}$-representation $\psi_{\mu}$ introduced in section
\ref{sec 4.2}:
\begin{equation}
\mathcal{A}_{\mu}(s_{\zeta_i},\delta,\nu_{s+i-1}) \colon
\Hom_M(\mu,\delta) \longrightarrow  \Hom_M(\mu,\delta),\, T \mapsto
\psi_{\mu}(s_{\zeta_i}) T .
\end{equation}
Composing the various factors of the operator we find:
\begin{equation}
\mathcal{A}_{\mu}(w,\delta,\nu) =\underbrace{
\psi_{\mu}(s_{\zeta_{s+t-1}}) \cdots
\psi_{\mu}(s_{\zeta_1})}_{\psi_{\mu}(u)}
\mathcal{A}_{(\psi_{\mu}^0)\ast}(w^0,\nu) =
\psi_{\mu}(u)\mathcal{A}_{(\psi_{\mu}^0)\ast}(w^0,\nu) .
\end{equation}
This concludes the proof of the theorem.
\end{proof}

\begin{remark}
If $\mu$ is fine, then \[\mathcal{A}_{\mu}(w,\delta,\nu) =
\psi_{\mu}(u)\mathcal{A}_{(\psi_{\mu}^0)\ast}(w,\nu) =
\psi_{\mu}(u)\psi_{\mu}^0(w^0)=\psi_{\mu}(w)=\varrho_{\mu}[w]\] (in
accordance with proposition \ref{prop:4.2}).
\end{remark}

\smallbreak

\subsection{Matching of unitary duals}

Let $G$ be a split real group.  Suppose, for simplicity, that the
principal series $X^G(\delta,\nu)$ has a unique irreducible
Langlands subquotient $L^G(\delta,\nu)$. Let $w\in W$ be such that
\[
w\cdot \delta \simeq \delta\quad w\cdot \nu = -\bar{\nu} \quad
\text{and} \quad w\cdot Q = \bar{Q}.
\]
The unitarity of the Langlands quotient $L^G(\delta,\nu)$ depends on
the signature of the Hermitian intertwining operator
$\mathcal{A}(w,\delta,\nu)$. More precisely, $L^G(\delta,\nu)$ is
unitary if and only if $\mathcal{A}_{\mu}(w,\delta,\nu)$ is
positive-semidefinite for all $\mu \in \widehat{K}$ containing
$\delta$.

Because the Weyl group element $w$ stabilizes $\delta$,   we can
write $w=uw^0$, with $u\in R^c_{\delta}$ and $w^0\in W_{\delta}^0$.
 If  $\mu$  is a \emph{petite} $K$-type containing $\delta$, then
\begin{equation}
 \mathcal{A}_{\mu}(w,\delta,\nu)= \psi_{\mu}(u) \mathcal{A}_{(\psi_{\mu}^0)*}
 (w^0,\nu).
\end{equation} Here $\psi_{\mu}$ is the $W_{\delta}$-representation
on $\Hom_{M}(\mu,\delta)$ and $\psi_{\mu}^0$ is its restriction to
$W_{\delta}^0$. The operator
$\mathcal{A}_{(\psi_{\mu}^0)*}(w^0,\nu)$ is a spherical intertwining
operator for a graded Hecke algebra $\mathbb{H}^0$ defined as in
section \ref{sec:3.1}, with $W$=$W_{\delta}^0$ and $\Pi$ the simple
roots of $\Delta_{\delta}$.

 Recall from section \ref{sec:3.4} that the
unitarity  of the irreducible spherical $\mathbb{H}^0$-module
$L^{\mathbb{H}^0}(\nu)$ is detected by a finite number of
\emph{relevant} $W_{\delta}^0$-types.

If $A(\delta)$ has cardinality one, then $W_{\delta}^0=W_{\delta}$,
$R_{\delta}$ is trivial and $u=1$.  Hence \begin{equation}
 \mathcal{A}_{\mu}(w,\delta,\nu)= \mathcal{A}_{(\psi_{\mu}^0)*}
 (w,\nu)\end{equation}
  for all $\mu$ petite.

Suppose that every relevant $W_{\delta}^0$ type comes from a petite $K$-type via
the correspondence $ \mu \rightarrow \psi_{\mu}^0$. In this case,
the unitarity of the $G$-module $L^G(\delta,\nu)$ implies the
unitarity of Hecke algebra-module $L^{\mathbb{H}^0}(\nu)$.
%\[\]\[\]

\[
\text{\doublebox{$L^G(\delta,\nu)$ is unitary}}
\,{\color{black}{\,\;=====================>\;\,}}
\text{\doublebox{$L^{\mathbb{H}^0}( \nu)$ is unitary}}
\]
\[
 \quad \Updownarrow \qquad \qquad \qquad \qquad\qquad \qquad \qquad \qquad
\qquad \qquad \qquad  \qquad \qquad \quad \quad \Updownarrow
\]
\[
\boxed{\begin{array}{c} \mathcal{A}_{\mu}(w,\delta,\nu) \geq 0\\
\\
\forall \, \mu \in \widehat{K}
\end{array}}
\Rightarrow
\boxed{\begin{array}{c} \mathcal{A}_{\mu}(w,\delta,\nu) \geq 0\\
\\
\forall \, \mu \text{ petite}
\end{array}}
\Rightarrow \boxed{\begin{array}{c}\mathcal{A}_{\tau}
 (w,\nu)\geq 0\\
\\
\forall \, \tau \text{ relevant}
\end{array}}
\Leftrightarrow \boxed{\begin{array}{c} \mathcal{A}_{\tau}
 (w ,\nu) \geq 0\\
\\
\forall \, \tau \in \widehat{W_{\delta}^0}
\end{array}}
\]
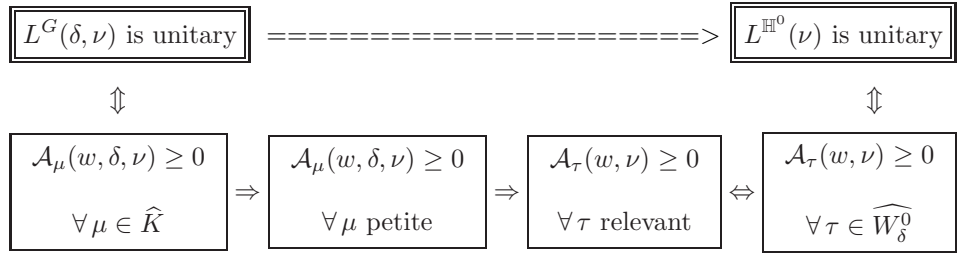
\[\]
It follows that the portion of the unitary dual of $G$ induced by
$\delta$ is embedded in the spherical unitary dual of
$\mathbb{H}^0$.
\begin{theorem}[\cite{Ba1},\cite{Ba2}]
Let $G$ be any real split group and let  $\delta$ be the trivial
representation of $M$. Then  $R_{\delta}=\{1\}$,
$W_{\delta}^0=W_{\delta} = W$ and every relevant $W$-type comes from
a petite $K$-type (via
 $ \mu \rightarrow \psi_{\mu}$).
%\begin{enumerate}

%\item
As a consequence, the spherical unitary dual of $G$  is embedded in the spherical unitary dual of
the corresponding Hecke algebra $\mathbb{H}$.

%\item If $G$ is a \emph{classical} group,  the inclusion becomes an equality.
%\end{enumerate}
\end{theorem}

Suppose that $A(\delta)$ has cardinality bigger than one, so that
 $R_{\delta}$ is nontrivial. Because we are
 assuming the existence of a unique irreducible subquotient,
  every fine $K$-type in $A(\delta)$ is contained in $L^G(\delta,\nu)$. We
 distinguish two cases:
\begin{itemize}
\item[(a)] $w \in W_{\delta}^0$
\item[(b)] $w \in W_{\delta} \diagdown W_{\delta}^0$, and there is
no  $w^0 \in W_{\delta}^0$ satisfying
\[
 w^0\cdot \nu = -\bar{\nu} \quad
\text{and} \quad w^0\cdot Q = \bar{Q}.
\]
\end{itemize}

If  $w \in W_{\delta}^0$, the operator $ \mathcal{A}(w,\delta,\nu)$
acts trivially on every fine $K$-type contained in
$L^G(\delta,\nu)$, and acts by
\[
 \mathcal{A}_{(\psi_{\mu}^0)*}
 (w ,\nu)
 \]
 on every petite $K$-type containing $\delta$. Assume that
every relevant $W_{\delta}^0$ type can be matched with a petite
$K$-type. Then the same analysis
 performed above shows that  $L^G(\delta,\nu)$ is unitary only if
  $L^{\mathbb{H}^0}(\delta,\nu)$ is unitary. Hence we obtain an
 embedding of unitary duals.

If $w \in W_{\delta} \diagdown W_{\delta}^0$, and there is no  $w^0
\in W_{\delta}^0$ satisfying
\[
w^0\cdot \nu = -\bar{\nu} \quad \text{and} \quad w^0\cdot Q =
\bar{Q},
\]
then there are at least two fine $K$-types $\pi_1$ and $\pi_2$ such
that \begin{equation} \psi_{\pi_1}(w) = \varrho_{\pi_1}[w] = -
\varrho_{\pi_2}[w] = - \psi_{2}(w).\end{equation} Because the
operators
 $ \mathcal{A}_{\pi_1}(w,\delta,\nu)$ and  $
 \mathcal{A}_{\pi_1}(w,\delta,\nu)$ have opposite sign, the Langlands quotient
$L^G(\delta,\nu)$ is not unitary.

\subsubsection{Example: $Sp(2n,\bR)$}\label{sec:3a.7.1}
Let $G$ be the real split group $Sp(2n,\bR)$. Then $K=U(n)$ and
$M=(\mathbb{Z}_2)^n.$ There are $(n+1)$ $W$-conjugacy classes of
$M$-types:
 the spherical $M$-type
\[
\delta_0= trivial
\]
and the nonspherical $M$-types
\[
\delta_p=
(\underbrace{+,+,\dots,+}_{n-p},\underbrace{-,-,\dots,-}_p) \qquad
p=1\dots n.
\]
 Recall from section \ref{sec:2.5} that $W_{\delta_0}^0=W_{\delta_0}=W$,  while $
W_{\delta_p}^0=C_{n-p}\times D_p $ and $ W_{\delta_p}=C_{n-p}\times
C_p $ for $p\geq 1$.  The $R_{\delta_p}$-group is
$\frac{\mathbb{Z}}{2\mathbb{Z}},$ and indeed $\delta_p$ is contained
into two fine $K$-types:
\[
\mu_p(+) = (\underbrace{1,\dots ,1}_{p},\underbrace{0,\dots
,0}_{n-p}) = \Lambda^p(\mathbb{C}^n)
\]
and
\[
\mu_p(-) = (\underbrace{0,\dots ,0}_{n-p},\underbrace{-1,\dots
,-1}_{p}) = [\Lambda^p(\mathbb{C}^n)]^\ast.
\]
If $p=1\dots n$, we will assume that the last entry of
\[\nu=(a_1,a_2,\dots,a_n)\]
is nonzero and that $p$ is even. Under these assumptions, the
principal series $X(\delta_p,\nu)$ has a unique irreducible
subquotient and there exists $w\in W_{\delta_p}^0 $ mapping $P$ into
$\bar{P}$ and $\nu$ into $-\bar{\nu}$.

The case of $Sp(2n,\bR)$ is very special, because  every relevant
$W_{\delta}$-type can be matched with a petite $K$-type (via
$\psi_{\mu} \leftrightarrow \mu$). The matching is as follows:

\[
  \begin{tabular}{|c|c|}
\hline
&\\
$\qquad \quad \text{relevant $W(C_n)$-type $\psi$}$ \qquad \quad &$  \text{spherical petite $K$-type such that $\psi_{\mu}=\psi$}$ \\
 &\\
\hline \hline & \\
$ (n-k)\times (k) $&$(\unb{k}{2,\dots
,2},0,\dots ,0)$\\
% &\\
\hline & \\
$(n-k,k)\times (0) $&$(\unb{k}{1,\dots
,1},0,\dots ,0,\unb{k}{-1,\dots ,-1})$ \\
%& \\
\hline \hline
  \end{tabular}
\]

\[\]

\[
\begin{tabular}{|c|c|}
\hline
&\\
$\text{relevant $W(C_{n-p}\times C_{p})$-type $\psi$}$&$  \text{nonspherical petite $K$-type such that $\psi_{\mu}=\psi$}$ \\
 &\\
\hline \hline & \\
 $(triv)\otimes [(a,p-a)\times (0)]$
&$(\unb{a}{1,\dots,1},\unb{n-p}{1,\dots,1},0,\dots,0,\unb{a}{-1,\dots,-1})$\\
%&\\
\hline & \\
$(triv)\otimes [(p-c)\times (c)] $&$(\unb{c}{2,\dots
,2},\unb{n-p}{1,\dots,1},0,\dots ,0)$\\
% &\\
\hline & \\
$[(n-p-c)\times (c)]\otimes (triv) $&$(\unb{n-p-c}{1,\dots
,1},0,\dots ,0,\unb{c}{-1,\dots ,-1})$ \\
%& \\
 \hline & \\
$[(c,n-p-c)\times (0)]\otimes (triv) $&$(\unb{c}{2,\dots
,2},\unb{n-p-2c}{1,\dots ,1},\unb{p+c}{0,\dots,0})$ \\
%&\\
\hline \hline
  \end{tabular}
\]
\[\]
To prove that the $K$-types in these tables are actually petite, we
observe that they all appear in the
 tensor products
\begin{equation*}
  \begin{aligned}
&\mu_+(p)\otimes\mu_-(p)=
\sum_{2a+b=2p} (\underbrace{1,\dots ,1}_{a},\underbrace{0,\dots ,0}_{b},\underbrace{-1,\dots ,-1}_{n-a-b}),\\
&\mu_+(p)\otimes\mu_+(k)= \sum_{2a+b=2p} (\underbrace{2,\dots
,2}_{a},\underbrace{1,\dots ,1}_{b},\underbrace{0,\dots
,0}_{n-a-b}).
\end{aligned}
\end{equation*}
Because $\mu_{\pm}(p)$ has level 1, every summand of these
 decompositions has level at most two and is automatically petite.
We get another set of petite $K$-types by changing all the signs to
minuses, \ie by passing to the dual.

The existence of this matching between relevant $W_{\delta}$-types
and petite $K$-types containing $ \delta$  allows us to draw the
following conclusions.
\begin{itemize}

\item If $\delta=\delta_0$, let $\mathbb{H}$ be the affine graded Hecke algebra corresponding to $W$
(defined as in section \ref{sec:3.1}). Then the spherical Langlands
quotient $L^G(\delta_0,\nu)$ is unitary only if the spherical Hecke
algebra-module $L^{\mathbb{H}}(\nu)$ is unitary.

\item If $\delta=\delta_p$, assume that $a_n \neq 0$ and $p$ is even. Let $\mathbb{H}^0$ be the affine graded Hecke algebra
corresponding to $W_{\delta}^0$ (defined as in section
\ref{sec:3.1}). Then the nonspherical Langlands quotient
$L^G(\delta_p,\nu)$ is unitary only if the spherical Hecke
algebra-module $L^{\mathbb{H}^0}(\nu)$ is unitary.

\end{itemize}

\subsubsection{Remarks} So far we assumed that the principal series
$X(\delta,\nu)$ had a unique irreducible subquotient. We conclude
this section with some brief remarks on the general case.

The dual $R$-group  $\widehat{R_{\delta}}$ acts simply transitively
on $A(\delta)$, so the number of fine $K$-types containing $\delta$
is equal to the cardinality of $R_{\delta}$ (recall that
$R_{\delta}$ is an abelian group, isomorphic to its dual). Two fine
$K$-types occur in the same irreducible subquotient if and only if
they lie in the same orbit of $ R_{\delta}^{\perp}(\nu)$ (the
annihilator of $R_{\delta}(\nu)$ inside $\widehat{R_{\delta}}$).
Every subquotient contains the same number of  fine $K$-types, equal
to the cardinality of $ R_{\delta}^{\perp}(\nu)$. We conclude that
there are exactly $ \frac{\# \widehat{R_{\delta}}}{\#
R_{\delta}^{\perp}(\nu)} $ irreducible subquotients. It is easy to
see that $R_{\delta}^{\perp}(\nu)$ is the Kernel of the restriction
map
\[ \widehat{R_{\delta}} \rightarrow
\widehat{R_{\delta}(\nu)}, \, \chi \mapsto \chi|_{R_{\delta}(\nu)}.
\]
(Since the groups are abelian, all irreducible representations are
one-dimensional, so this is map is well defined and surjective.)
Hence
\[
\# \text{ irreducible subquotients } = \#
\frac{\widehat{R_{\delta}}}{R_{\delta}^{\perp}(\nu)}=\#\widehat{R_{\delta}(\nu)}=\#
R_{\delta}(\nu).
\]

Now suppose that $R_{\delta}(\nu)$ is not trivial. Then the
principal series $X(\delta,\nu)$ contains several irreducible
subquotients, and the intertwining operator
$\mathcal{A}_{\mu}(w,\delta,\nu)$ has a block diagonal structure
(with one block per subquotient). The blocks are in one-to-one
correspondence
 the set of orbits of $ R_{\delta}^{\perp}(\nu)$ in $A(\delta)$.

If $\mu$ is petite, the operator $\mathcal{A}_{\mu}(w,\delta,\nu)$
matches the operator $\mathcal{A}'_{\psi_{\mu}^\ast}(w,\nu)$ for the
quasi-spherical module $X'(\nu)$ for the extended graded Hecke
algebra $\mathbb{H}'$ corresponding to $W_{\delta}$ (\cf~section
\ref{sec:3.7}). The operator $\mathcal{A}'_{\psi_{\mu}^\ast}(w,\nu)$
also has a   block diagonal structure, with one block per character
of $R_{\delta}(\nu)$.

We notice  that there is a one-to-one correspondence between blocks
of $\mathcal{A}_{\mu}(w,\delta,\nu)$ and blocks of
$\mathcal{A}'_{\psi_{\mu}^\ast}(w,\nu)$. Equivalently, there is a
bijection between
 orbits  of $ R_{\delta}^{\perp}(\nu)$ in $A(\delta)$ and characters of
 $R_{\delta}(\nu)$:
fix a base point $\mu_{\delta}$ in  $A(\delta)$ (\ie a fine $K$-type
containing $\delta$) and identify $A(\delta)$ with
$\widehat{R_{\delta}}$ via $\mu \mapsto \psi_{\mu}$. The set of
orbits is then identified with $
\frac{\widehat{R_{\delta}}}{R_{\delta}^{\perp}(\nu)}$, which we know
is isomorphic to $\widehat{R_{\delta}(\nu)}$.

Then, to check the signature of the Hermitian form on a petite
$K$-type $\mu$ occurring in the subquotient
$L(\delta,\nu)(\mu_{\delta})$, one can look at the the signature of
the appropriate block of the Hecke algebra operator
$\mathcal{A}'_{\psi_{\mu}^\ast}(w,\nu)$. Similarly for the other
subquotients (but an issue with the normalization may arise).

\section{Spherical unitary dual}\label{sec:4}

In this section we describe the spherical unitary modules, with real
infinitesimal (central) character. To distinguish
between the real and the $p$-adic case, we will denote by
$L^{\bR}(\chi)$ and $L^{\bH}(\chi)$ the corresponding spherical
modules for $G(\bR),$ respectively Hecke algebra $\bH,$ and
similarly for all related notation.

\begin{comment}
We will call a
parameter $\chi\in \check\fh$ {\it unitary} if the corresponding
spherical module $L(\chi)$ is unitary.
\end{comment}

\subsection{Parameters}\label{sec:4.1} To every $\chi\in\check\fh_\bR$ we attach uniquely a
nilpotent $\check G$-orbit $\ch\CO(\chi)$ in $\check\fg$ as
follows. Consider
\begin{equation}\label{eq:4.1.1}
\check\fg_{1,\chi}=\{x\in \check\fg: ad(\chi)(x)=x\},\quad \check
G_{0,\chi}=\{g\in \check G: Ad(g)(\chi)=\chi\}.
\end{equation}
It is known that $\check G_{0,\chi}$ acts with finitely many orbits on
  $\check\fg_{1,\chi}$, and as a consequence, there is a unique open
  orbit in $\fg_{1,\chi}.$ Let $\check\CO(\chi)$ denote the $\check
  G$-saturation of this open orbit. For the classification and
  relevant facts about nilpotent orbits in complex Lie algebras, we
  refer the reader to \cite{Ca} and \cite{CM}.

Let $\check\CO$ be a nilpotent $\check G$-orbit in $\check\fg.$
\begin{definition}[1]\label{d:4.1} The {\it $\check\CO$-complementary
    series} are the
  sets
\begin{align}\label{eq:4.1.2}
&CS^{\bH}(\check\CO)=\{\chi: L^\bH(\chi) \text{ is unitary and }
\check\CO(\chi)=\check\CO\},\text{ in the $p$-adic case},\\
&CS^{\bR}(\check\CO)=\{\chi: L^\bR(\chi) \text{ is unitary and }
\check\CO(\chi)=\check\CO\},\text{ in the real case}.
\end{align}
When we want to refer to both sets simultaneously, we will just use
the notation $CS(\check\CO).$

\end{definition}
The spherical unitary parameters are the disjoint union of the
corresponding complementary series $\sqcup_{\check\CO}
CS(\check\CO).$

Fix a Lie triple $\{\check e,\check h,\check f\}$ in $\check\CO.$ Let
$\fz(\check e,\check h,\check f)$ denote the centralizer in
$\check\fg$ of $\check e,\check h,\check f.$ Then the orbit
$\check\CO(\chi)$ can be described differently. By \cite{BM1}, the
orbit $\check\CO(\chi)$ is the unique one satisfying the conditions:
\begin{align}\label{eq:4.1.3}
&\text{(1) $w\chi=\check h/2+\nu,$ for some $\nu\in\fz(\check e,\check
  h,\check f),$ and}\\\notag
&\text{(2) $\check\CO(\chi)$ is maximal with respect to condition (1).}
\end{align}
Clearly $\check\CO(\check h/2)=\check\CO.$ In fact these parameters
are special, they are instances of {\it Arthur parameters}.

\begin{definition}[2] The modules $L^\bH(\check h/2)$ and
  $L^\bR(\check h/2)$ are called {\it
    special unipotent}.
\end{definition}

The conjectures in \cite{A} suggest that the special unipotent parameters
should be unitary. In the Hecke algebra case, there exists the
Iwahori-Matsumoto involution $IM$ which preserves unitarity and takes
the special unipotent modules to tempered modules. The unitarity of
the corresponding group representations when $G$ is $p$-adic is then
implied by \cite{BM1}.

In the real case, the unitarity of special unipotent modules is proved
in \cite{Ba1} for split classical groups. This is beyond the scope of
this exposition, and we refer the reader to
section 9 in \cite{Ba1} for details.

\medskip

Returning to the $\check\CO$-complementary series, we see that
$CS^{\bH}(\check\CO)$ contains at least one element, $\check h/2$,
corresponding to the special unipotent. In fact, when $\check\CO$ is a
distinguished orbit, the conditions (\ref{eq:4.1.3}) imply that
$\check h/2$ is the unique element of $CS^{\bH}(\check\CO).$

\subsection{$0$-complementary series}\label{sec:4.2} The basic case one needs to
compute is when $\check\CO$ is the trivial nilpotent orbit. The
parameters $\chi$, which we will assume dominant, such that $\check\CO(\chi)=0$ are precisely those such
that
\begin{equation}\label{eq:4.2.1}
\langle\check\al,\chi\rangle\neq 1,\text{ for any }\al\in\Delta^+.
\end{equation}
In the (adjoint) $p$-adic case, these parameters correspond to the
modules which are both spherical and admit Whittaker models, in other
words, to the irreducible principal series $X(\chi).$

It is clear that the operators $\CA_\mu(w_0,\chi)$ in section
\ref{sec:3.2} are isomorphisms in any open region of the complement
of the hyperplane arrangement given by (\ref{eq:4.2.1}) in the
dominant Weyl chamber in $\check\fh_\bR.$ Due to the Hermitian
condition, we may only consider $\chi$ lying in the
$(-1)$-eigenspace of $w_0$:
\begin{equation}\label{eq:4.2.1a}
E_0=\{\chi\in\check\fh_\bR: w_0\chi=-\chi\}.
\end{equation}

Consequently,
these operators have constant signature inside any such open regions,
and we see that $CS(0)$ is a union of open regions intersected with $E_0.$

Let $\C C_0$ denote the {\it fundamental alcove}:
\begin{equation}\label{eq:4.2.2}
\C C_0=\{\chi\in\check\fh_\bR: 0\le\langle\check\al,\chi\rangle<1,\text{ for all }\al\in\Pi\}.
\end{equation}
Any open region conjugate to $\C
C_0$ by the affine Weyl group is called an alcove.

The following result is well-known. We will regard $CS(0)$ in this
section as a subset of the fundamental Weyl chamber, \ie we will only
consider dominant $\chi.$

\begin{lemma}\label{l:4.2}

\begin{enumerate}
\item $\C C_0\cap E_0\subset CS(0)$.
\item Every open region contributing to $CS(0)$ is bounded.
\end{enumerate}

\end{lemma}

%\begin{proof}
Part (1) is implied by the fact that $X(0)$ is irreducible and
unitary. Part (2) is a signature calculation. In the Hecke algebra
setting, one can easily show that $\CA_{\mu}(w_0,\nu)$, with $\mu$
the reflection $W$-representation, can only be positive definite in
the bounded regions.
%\end{proof}

\medskip

The first result is the determination of the $0$-complementary
series. The proof for classical groups is in \cite{BM3} and
\cite{Ba1}, with a different proof in \cite{BC3}, while for exceptional groups it is in \cite{Ci} for $F_4$
and \cite{BC} for $E_6,E_7,E_8$. Earlier, for type $A$, this result was obtained in
\cite{V2} in the real case and \cite{Ta} in the $p$-adic case, while
for $G_2$, it is part of \cite{V3} for the real case, and \cite{Mu}
for the $p$-adic.

\begin{theorem} The $0$-complementary series $CS(0)$ is formed of:

\begin{enumerate}
\item[A.] $\C C_0\cap E_0$;
\item[B.] $\C C_0$;
\item[C,D.] $2^{[n/2]}$ alcoves intersected with $E_0$, where $n$ is the rank of the group;
\item[G2.] $2$ alcoves;
\item[F4.] $2$ alcoves;
\item[E.] $2$ alcoves intersected with $E_0$ for $E_6,$ $2^3$ alcoves
  for $E_7$, and $2^4$ alcoves for $E_8.$
\end{enumerate}

 In explicit coordinates, in type $C_n$ or $D_n$, if $\chi=(\nu_1,\dots,\nu_n)$,
  with $w_0\chi=-\chi$, the set
  $CS(0)$ is formed of the parameters satisfying the following
  condition: there exists an index $i$ such that
  \begin{align}\label{eq:4.2.3}\notag
&0\le\nu_1\le\dotsb\le\nu_i<1-\nu_{i-1}<\nu_{i+1}<\dotsb<\nu_n<1,
\text{ and between any $\nu_j<\nu_{j+1}$,} \\
&i\le j<n,
\text{ there is an odd
  number of $(1-\nu_l)$, $1\le l<i.$}
\end{align}
The explicit description for the exceptional types is in section \ref{sec:zero}.
\end{theorem}

%\begin{proof}
In the $p$-adic case, the proof has two components:

\begin{enumerate}
\item in any region which has a wall of the form
  $\langle\check\al,\chi\rangle=0,$ for some simple coroot
  $\check\al$, one knows by unitary induction the signature of the
  operators $\CA^\bH_\psi(w_0,\nu)$ on every $W$-type $\psi.$
\item any open region for which all the walls are of the form
  $\langle\check\al,\chi\rangle=1,$ for $\al\in\Delta^+,$ cannot be
  unitary. One can use a signature argument here, for example in the
  classical cases one can show that the signature of the operator
  $\CA^\bH_\psi(w_0,\chi),$ where $\psi=(\text{refl})+Sym^2(\text{refl})$ is indefinite (see \cite{BC3}). For
  $E_7$ and $E_8,$ a different argument is used in \cite{BC}, namely
  it is shown that any such region must have a codimension two edge
  given by two coroots which form an $A_2$, and by a simple signature
  argument, the parameter cannot be unitary in such a
  case.

Explicit calculations showed that in fact in all cases,
  classical and exceptional, a parameter $\chi$ is in $CS(0)$ if and
  only if the operator $\CA^\bH_\psi(w_0,\chi)$ is positive definite for
  $\psi=(\text{refl})+Sym^2(\text{refl}),$ but we do not have a conceptual proof of
  this fact for $E_7$ and $E_8.$ However the outline above gives conceptual proofs
  in all cases of the fact that $\chi$ is in $CS(0)$ if and only if
\begin{equation}\label{eq:4.2.4}
\CA^\bH_\psi(w_0,\chi) \text{ is positive definite for all relevant
  $W$-types }\psi.
\end{equation}

\end{enumerate}

In the real case, the same result holds by the following extra argument:

\begin{enumerate}
\item[(3)] by (\ref{eq:4.2.4}), the unitary spherical and generic
  parameters for $G(\bR)$ are a subset of $CS(0)$;
\item[(4)] any $\chi\in CS(0)$ has the property that
  $\langle\check\al,\chi\rangle\neq m,$ for all roots
  $\al\in\Delta^+$, and all $m\in\bZ_{>0}.$
  Therefore the spherical minimal principal series $X(\chi)$ in
  the real case is irreducible as well for all $\chi\in CS(0)$;
\item [(5)] also in the real case, any $\chi\in CS(0)$ can be proven
  unitary by irreducible deformations and unitary induction.

\end{enumerate}

%\end{proof}

\subsection{}\label{sec:zero}
 We record next the precise description of the $0$-complementary series
(that is, the generic spherical unitary
parameters) for simple exceptional split groups. We use the Bourbaki simple
roots in type $E$. To simplify the notation, we will write $\check\al<1$ instead
of $\langle\check\al,\chi\rangle<1$ in the description of the unitary
regions.

\subsubsection{$G_2$}\label{sec:genspherG2} The parameter is
$\chi=(\nu_1,\nu_1+\nu_2,-2\nu_1-\nu_2),$ with $\nu_1\ge 0,$ $\nu_2\ge
0.$ The $0$-complementary series is:

\begin{itemize}
\item[(1)] $\check\al_6<1,$ and $\check\al_1,\check\al_2\ge 0.$
\item[(2)] $\check\al_4<1,$ $\check\al_5>1$, and $\check\al_1\ge 0.$
\end{itemize}

We use the coroots $\check\al_1=(2/3,-1/3,-1/3)$,
$\check\al_2=(-1,1,0)$, $\check\al_6=(0,1,-1)$,
$\check\al_5=(1,0,-1),$ and $\check\al_4=(1/3,1/3,-2/3).$

\subsubsection{$F_4$}\label{sec:genspherF4} The parameter is
$\chi=(\nu_1,\nu_2,\nu_3,\nu_4)$ assumed dominant. The
$0$-complementary series is:

\begin{itemize}
\item[(1)] $\check\al_{24}<1,$ and
  $\check\al_1,\check\al_2,\check\al_3,\check\al_4\ge 0.$
\item[(2)] $\check\al_{22}<1,$ $\check\al_{23}>1$, and
  $\check\al_1,\check\al_2,\check\al_4\ge 0.$
\end{itemize}

We use the coroots $\check\al_1=(1,-1,-1,-1),$
$\check\al_2=(0,0,0,2),$ $\check\al_3=(0,0,1,-1)$,
$\check\al_4=(0,1,-1,0),$ and $\check\al_{24}=(2,0,0,0),$
$\check\al_{23}=(1,1,1,1)$, $\check\al_{22}=(1,1,1,-1).$

\subsubsection{$E_6$}\label{sec:genspherE6} In $W(E_6)$, the longest Weyl
group element $w_0$ does not act by $-1.$ Therefore, we only consider
dominant parameters $\chi$ such that $w_0\chi=-\chi.$
\noindent In coordinates,
\begin{small}\begin{align}\notag
\chi=&\left(\frac {\nu_1-\nu_2}2-\nu_3,\frac {\nu_1-\nu_2}2-\nu_4,\frac {\nu_1-\nu_2}2+\nu_4,\frac {\nu_1-\nu_2}2+\nu_3,\frac {\nu_1+\nu_2}2,\notag\right.\\
&\left. -\frac {\nu_1+\nu_2}2,-\frac {\nu_1+\nu_2}2,\frac {\nu_1+\nu_2}2\right).
\end{align}\end{small}

The $0$-complementary series is:

\begin{itemize}
\item[(1)] $\check\al_{36}<1,$ and $\check\al_1,\check\al_2,\check\al_3,\check\al_4,\check\al_5,\check\al_6\ge 0.$
\item[(2)] $\check\al_{34}<1,$ $\check\al_{35}>1,$ and $\check\al_1,\check\al_2,\check\al_3,\check\al_5,\check\al_6\ge 0.$
\end{itemize}

\subsubsection{$E_7$}\label{sec:genspherE7}

The  parameters are
$\chi=(\nu_1,\nu_2,\nu_3,\nu_4,\nu_5,\nu_6,-\nu_7,\nu_7),$ assumed
dominant.

The $0$-complementary series is:

\medskip

\begin{enumerate}

\item $\check\alpha_{63}<1,$ and $\check\al_1,\check\al_2,\check\al_3,\check\al_4,\check\al_5,\check\al_6,\check\al_7\ge 0.$

\item $\check\alpha_{61}<1,\ \ \check\alpha_{62}>1$ and
  $\check\al_1,\check\al_2,\check\al_4,\check\al_5,\check\al_6,\check\al_7\ge 0.$

\item $\check\alpha_{58}<1,\ \check\alpha_{59}<1,\ \ \check\alpha_{60}>1$ and
  $\check\al_1,\check\al_3,\check\al_4,\check\al_6,\check\al_7\ge 0.$

\item $\check\alpha_{53}<1,\ \check\alpha_{54}<1,\ \check\alpha_{55}<1,\ \
  \check\alpha_{56}>1,\ \check\alpha_{57}>1$ and $\check\al_1,\check\al_3,\check\al_5\ge 0.$

\item $\check\alpha_{46}<1,\ \check\alpha_{47}<1,\ \check\alpha_{48}<1,\
  \check\alpha_{49}<1,\ \  \check\alpha_{50}>1,\   \check\alpha_{51}>1,\  \check\alpha_{52}>1$
  and $\check\al_2\ge 0.$

\item $\check\alpha_{53}<1,\ \check\alpha_{59}<1,\ \ \check\alpha_{56}>1$ and
  $\check\al_1,\check\al_3,\check\al_4,\check\al_5,\check\al_6\ge 0.$

\item $\check\alpha_{49}<1,\ \check\alpha_{53}<1,\ \check\alpha_{54}<1,\ \
  \check\alpha_{52}>1,\ \check\alpha_{56}>1$ and $\check\al_3,\check\al_4,\check\al_5\ge 0.$

\item $\check\alpha_{47}<1,\ \check\alpha_{48}<1,\ \check\alpha_{49}<1,\ \check\alpha_{53}<1,\
  \ \check\alpha_{51}>1,\ \check\alpha_{52}>1$ and $\check\al_2,\check\al_4\ge 0.$

\end{enumerate}

\subsubsection{$E_8$}\label{sec:genspherE8}

The parameters are
$\chi=(\nu_1,\nu_2,\nu_3,\nu_4,\nu_5,\nu_6,\nu_7,\nu_8),$  assumed
dominant.

The $0$-complementary series is:

\medskip

\begin{enumerate}

\item $\check\alpha_{120}<1$ and $\check\alpha_1,\check\alpha_2,\check\alpha_3,\check\alpha_4,\check\alpha_5,\check\alpha_6,\check\alpha_7,\check\alpha_8\ge0$.

\item $\check\alpha_{113}<1,\check\alpha_{114}<1$; $\check\alpha_{115}>1$ and $\check\alpha_1,\check\alpha_4,\check\alpha_5,\check\alpha_6,\check\alpha_7,\check\alpha_8\ge 0$.

\item $\check\alpha_{109}<1,\check\alpha_{110}<1$; $\check\alpha_{111}>1,\check\alpha_{112}>1$ and $\check\alpha_3,\check\alpha_5,\check\alpha_6,\check\alpha_7,\check\alpha_8\ge 0$.

\item $\check\alpha_{91}<1,\check\alpha_{92}<1,\check\alpha_{97}<1,\check\alpha_{98}<1$; $\check\alpha_{95}>1,\check\alpha_{96}>1,\check\alpha_{101}>1$ and $\check\alpha_3,\check\alpha_4\ge 0$.

\item $\check\alpha_{90}<1,\check\alpha_{91}<1,\check\alpha_{92}<1,\check\alpha_{97}<1$; $\check\alpha_{94}>1,\check\alpha_{95}>1,\check\alpha_{96}>1$  and $\check\alpha_1,\check\alpha_3\ge 0$.

\item $\check\alpha_{89}<1,\check\alpha_{90}<1,\check\alpha_{91}<1,\check\alpha_{92}<1$; $\check\alpha_{93}>1,\check\alpha_{94}>1,\check\alpha_{95}>1,$ $\check\alpha_{96}>1$ and $\check\alpha_1\ge 0$.

\item $\check\alpha_{104}<1,\check\alpha_{110}<1$; $\check\alpha_{107}>1,\check\alpha_{112}>1$ and $\check\alpha_3,\check\alpha_4,\check\alpha_5,\check\alpha_7,\check\alpha_8\ge 0$.

\item $\check\alpha_{104}<1,\check\alpha_{105}<1,\check\alpha_{106}<1$; $\check\alpha_{107}>1,\check\alpha_{108}>1$ and $\check\alpha_2,\check\alpha_4,\check\alpha_7,\check\alpha_8\ge 0$.

\item $\check\alpha_{118}<1$; $\check\alpha_{119}>1$ and $\check\alpha_1,\check\alpha_2,\check\alpha_3,\check\alpha_4,\check\alpha_5,\check\alpha_6,\check\alpha_8\ge 0$.

\item $\check\alpha_{97}<1,\check\alpha_{110}<1$; $\check\alpha_{101}>1,\check\alpha_{112}>1$ and $\check\alpha_3,\check\alpha_4,\check\alpha_5,\check\alpha_6,\check\alpha_7\ge 0$.

\item $\check\alpha_{97}<1,\check\alpha_{105}<1,\check\alpha_{106}<1$; $\check\alpha_{101}>1,\check\alpha_{108}>1$ and $\check\alpha_2,\check\alpha_4,\check\alpha_6,\check\alpha_7\ge 0$.

\item $\check\alpha_{116}<1$; $\check\alpha_{117}>1$ and $\check\alpha_1,\check\alpha_2,\check\alpha_3,\check\alpha_4,\check\alpha_6,\check\alpha_7,\check\alpha_8\ge 0$.

\item $\check\alpha_{97}<1,\check\alpha_{98}<1,\check\alpha_{106}<1$; $\check\alpha_{101}>1,\check\alpha_{102}>1$ and $\check\alpha_2,\check\alpha_4,\check\alpha_5,\check\alpha_6\ge 0$.

\item $\check\alpha_{97}<1,\check\alpha_{98}<1,\check\alpha_{99}<1$; $\check\alpha_{96}>1,\check\alpha_{101}>1,\check\alpha_{102}>1$ and $\check\alpha_2,\check\alpha_4,\check\alpha_5\ge 0$.

\item $\check\alpha_{97}<1,\check\alpha_{98}<1,\check\alpha_{99}<1,\check\alpha_{100}<1$; $\check\alpha_{101}>1,\check\alpha_{102}>1,\check\alpha_{103}>1$ and $\check\alpha_2,\check\alpha_5\ge 0$.

\item $\check\alpha_{114}<1$; $\check\alpha_{112}>1$ and $\check\alpha_1,\check\alpha_3,\check\alpha_4,\check\alpha_5,\check\alpha_6,\check\alpha_7,\check\alpha_8\ge 0$.

\end{enumerate}

\subsubsection{Roots for type $E$}\label{roots}
The notation for the positive coroots which appear in the
$0$-complementary series is as follows.

\medskip
\noindent{$E_6$}.

\begin{longtable}{ll}
$\check\al_{34}=\frac 12(-1,1,-1,1,1,-1,-1,1)$ &$\check\al_{35}=\frac
  12(-1,-1,1,1,1,-1,-1,1)$\\
$\check\al_{36}=\frac 12(1,1,1,1,1,-1,-1,1)$
\end{longtable}
\medskip

\noindent{$E_7$}.

\begin{longtable}{ll}
$\check\al_{46}=\frac 12(-1,1,-1,1,1,-1,-1,1)$ &$\check\al_{47}=\frac
  12(-1,1,1,-1,-1,1,-1,1)$\\
$\check\al_{48}=\frac 12(1,-1,-1,1,-1,1,-1,1)$ &$\check\al_{49}=\ep_5+\ep_6$\\

$\check\al_{50}=\frac 12(-1,-1,1,1,1,-1,-1,1)$ &$\check\al_{51}=\frac
  12(-1,1,-1,1,-1,1,-1,1)$\\
$\check\al_{52}=\frac 12(1,-1,-1,-1,1,1,-1,1)$ &$\check\al_{53}=\frac 12(1,1,1,1,1,-1,-1,1)$\\
$\check\al_{54}=\frac 12(-1,-1,1,1,-1,1,-1,1)$ &$\check\al_{55}=\frac
  12(-1,1,-1,-1,1,1,-1,1)$\\

$\check\al_{56}=\frac 12(1,1,1,1,-1,1,-1,1)$ &$\check\al_{57}=\frac
  12(-1,-1,1,-1,1,1,-1,1)$\\

$\check\al_{58}=\frac 12(1,1,1,-1,1,1,-1,1)$ &$\check\al_{59}=\frac
  12(-1,-1,-1,1,1,1,-1,1)$\\

$\check\al_{60}=\frac 12(1,1,-1,1,1,1,-1,1)$ &$\check\al_{61}=\frac
  12(1,-1,1,1,1,1,-1,1)$\\

$\check\al_{62}=\frac 12(-1,1,1,1,1,1,-1,1)$ &$\check\al_{63}=-\ep_7+\ep_8$\\

\end{longtable}

\medskip

\noindent{$E_8$}.

\begin{longtable}{ll}
$\check\al_{89}=\frac 12(1,-1,1,1,1,1,-1,1)$ &$\check\al_{90}=\frac
  12(1,1,-1,1,1,-1,1,1)$\\
$\check\al_{91}=\frac 12(1,1,1,-1,-1,1,1,1)$ &$\check\al_{92}=\frac
  12(-1,-1,-1,1,-1,1,1,1)$\\

$\check\al_{93}=\frac 12(-1,1,1,1,1,1,-1,1)$ &$\check\al_{94}=\frac
  12(1,-1,1,1,1,-1,1,1)$\\
$\check\al_{95}=\frac 12(1,1,-1,1,-1,1,1,1)$ &$\check\al_{96}=\frac
  12(-1,-1,-1,-1,1,1,1,1)$\\

$\check\al_{97}=-\ep_7+\ep_8$ &$\check\al_{98}=\frac 12(-1,1,1,1,1,-1,1,1)$\\
$\check\al_{99}=\frac 12(1,-1,1,1,-1,1,1,1)$ &$\check\al_{100}=\frac
  12(1,1,-1,-1,1,1,1,1)$\\
$\check\al_{101}=-\ep_6+\ep_8$ &$\check\al_{102}=\frac 12(-1,1,1,1,-1,1,1,1)$\\
$\check\al_{103}=\frac 12(1,-1,1,-1,1,1,1,1)$ &$\check\al_{104}=-\ep_5+\ep_8$\\
$\check\al_{105}=\frac 12(-1,1,1,-1,1,1,1,1)$ &$\check\al_{106}=\frac
  12(1,-1,-1,1,1,1,1,1)$\\

$\check\al_{107}=-\ep_4+\ep_8$ &$\check\al_{108}=\frac 12(-1,1,-1,1,1,1,1,1)$\\

$\check\al_{109}=-\ep_3+\ep_8$ &$\check\al_{110}=\frac 12(-1,-1,1,1,1,1,1,1)$\\

$\check\al_{111}=-\ep_2+\ep_8$ &$\check\al_{112}=\frac 12(1,1,1,1,1,1,1,1)$\\

$\check\al_{113}=\ep_1+\ep_8$ &$\check\al_{114}=-\ep_1+\ep_8$\\

$\check\al_{115}=\ep_2+\ep_8$ &$\check\al_{116}=\ep_3+\ep_8$\\

$\check\al_{117}=\ep_4+\ep_8$ &$\check\al_{118}=\ep_5+\ep_8$\\

$\check\al_{119}=\ep_6+\ep_8$ &$\check\al_{120}=\ep_7+\ep_8$\\

\end{longtable}

\subsection{$\CO$-complementary series}\label{sec:4.3} In this section we state the
main theorems as they follow from \cite{Ba1}, \cite{BC}.

\subsubsection{}
\begin{theorem}[\cite{Ba1},\cite{BC},\cite{Ci}]\label{t:4.3.1} Fix $\check\CO$ a
  nilpotent $\check G$-orbit in $\check\fg$ and a Lie triple $\{\check
  e,\check h,\check f\}.$ Let $\chi$ be a (hyperbolic) semisimple
  element such that $\check\CO(\chi)=\check\CO$, and which we write as $\chi=\check h/2+\nu,$ with
  $\nu\in\fz(\check e,\check h,\check f).$ Then
\begin{equation}\label{eq:4.3.1}
\chi \text{ is in }\check CS^\bH(\check\CO)\text{ if and only if } \nu
\text{ is in }  CS^{\bH(\fz(\check e,\check h,\check f))}(0),
\end{equation}
unless $\check\CO$ is one of the following orbits:
\begin{enumerate}
\item[$\bullet$] $A_1+\wti A_1$ in $F_4$,
\item[$\bullet$] $A_2+3A_1$ in $E_7$,
\item[$\bullet$] $A_4+A_2+A_1$, $A_4+A_2$, $D_4(a_1)+A_2$, $A_3+2A_1$,
  $A_2+2A_1$, and $4A_1$ in $E_8.$
\end{enumerate}
The explicit description of the complementary series is in section
\ref{sec:5}.
\end{theorem}

In the case of the exceptions, unless the orbit is $4A_1$
in $E_8,$ the complementary series is smaller than the one for the
centralizer, and for $4A_1$, it is larger.

\subsubsection{}

\begin{theorem}[\cite{Ba1},\cite{BC2}]\label{t:4.3.2} A spherical module
  $L^\bH(\chi)$ is unitary if and only if the operators
  $\CA^\bH_\psi(w_0,\chi)$ are positive semidefinite for all relevant
  $W$-types $\psi.$
\end{theorem}

\subsubsection{} We record the results for the spherical principal series
of a split real group.

\begin{theorem}[\cite{Ba1},\cite{Ba2}]\label{t:4.3.3} Every relevant $W$-type $\psi$
  appears as a $W$-subrepresentation of the $(V_\mu^M)^*$ space of a petite
  $K$-type $(\mu,V_\mu).$
\end{theorem}
The construction of petite $K$-types was explained in section
\ref{sec:3a}.

\begin{corollary}\label{c:4.3.3}
For every nilpotent orbit $\check\CO$, one has
\begin{equation}\label{eq:4.3.3}
CS^\bR(\check\CO)\subseteq CS^\bH(\check\CO).
\end{equation}
\end{corollary}

\subsubsection{}We have already seen in section \ref{sec:4.2} that in fact $CS^\bR(0)=CS^\bH(0).$

\begin{theorem}[\cite{Ba1}]\label{t:4.3.4} If $G(\bR)$ is split classical, then
\begin{equation}\label{eq:4.3.4}
CS^\bR(\check\CO)=CS^\bH(\check\CO),
\end{equation}
for every nilpotent orbit $\check\CO.$
\end{theorem}
In addition to the unitarity of the special unipotent representations
already mentioned in section \ref{sec:4.1}, in order to establish this
theorem, \cite{Ba1} needs to analyze the irreducibility of
parabolically induced representations (see section 10 in
\cite{Ba1}). One of the ideas which makes this tractable is a
combinatorial parameterization of the spherical representations, which
we recall next.

\subsection{Strings}\label{sec:4.4} In this section, we give a description of the
parameterization of the spherical dual as in \cite{Ba1}, section 2, by means of
a generalization of the {\it multisegments} or {\it strings}
introduced for $GL(n)$ in \cite{Ze}.

To every parameter $\chi$, one
attaches a multisegment, so that the orbit $\check\CO(\chi)$ and the
decomposition $\chi=\check h/2+\nu$ used in the previous section can
be read off easily. These multisegments arise naturally in the
setting of the geometry for the Hecke algebra. More precisely, they
parameterize the unique open $\check G_{0,\chi}$-orbit in
$\check\fg_{1,\chi}$ (see section \ref{sec:4.1}.)

Let $[a]=(a_1,\dotsc,a_k)$ be a set of numbers.
\begin{definition}\label{d:4.4} We call $[a]$ an {\it increasing
    (respectively decreasing)
    string} if $-a_{i-1}+a_i=1$ (respectively $a_{i-1}-a_i=1$) for all $i$.
\end{definition}

We will explain next how one builds from $\chi$: the multisegment, the
orbit $\check\CO(\chi)$, the middle element $\check h$ and the
parameter $\nu$ in the centralizer $\fz(\check e,\check h,\check f).$
Recall that for simple classical types the complex nilpotent orbits
$\check\CO$ are parameterized by partitions as follows:
\begin{enumerate}
\item[$\bullet$] partitions of $n$, when $\check\fg=sl(n,\bC)$;
\item[$\bullet$] partitions of $2n,$ with odd parts occurring with even
  multiplicity, when $\check\fg=sp(2n,\bC)$;
\item[$\bullet$] partitions of $2n+1$, with even parts occurring with
  even multiplicity, when $\check\fg=so(2n+1,\bC)$;
\item[$\bullet$] partitions of $2n$, with even parts occurring with
  even multiplicity, when $\check\fg=so(2n,\bC)$. In this case,
  when the partition is very even, \ie all parts are even, there are
  two distinct orbits corresponding to it.
\end{enumerate}

\subsubsection{$G$ of type $B_n$} We partition the entries of the
character $\chi=(\nu_1,\dotsc,\nu_n)$ into subsets $A_{\tau}$ where
$0\le\tau\le 1/2$ and
\begin{equation}\label{eq:4.4.1}
A_\tau=\{\nu_i: \nu_i\text{ or }-\nu_i\equiv 1/2+\tau (mod\ \bZ)\}.
\end{equation}
There are two cases $0\le\tau<1/2$ and $\tau=1/2.$

\smallskip

When $0\le\tau<1/2$, form $A'_\tau=A_\tau\sqcup(-A_\tau)$. We partition
$A'_\tau$ into a disjoint union of increasing strings
$M^+_{\tau,1},\dotsc,M^+_{\tau,\ell}$ and decreasing strings
$M^-_{\tau,1},\dotsc, M^-_{\tau,\ell}$, where $M^+_{\tau,i}=-M^-_{\tau,i}$ as
follows. Remove the smallest entry, say $-a$ in $A'_\tau$ and place it in
$M^+_{1,\tau}$, and the largest entry $a$ and place it in
$M^-_{1,\tau}.$ If $a-1$ appears in $A'_\tau$, remove it from
$A'_\tau$ and place it $M^-_{1,\tau}$
and similarly, remove $1-a$ and place it in $M^-_{\tau,1}$. Continue
with $a-2,a-3,\dotsc$ until this is not possible. This completes the
construction of $M^\pm_{\tau,1}$. Then repeat the
process with the remaining entries in $A'_\tau$ to construct
$M^\pm_{\tau,2}$, $M^\pm_{\tau,3}$, etc. Once this is finish, every
pair $(M^+_{\tau,i},M^-_{\tau,i})$ adds:

\begin{enumerate}
\item  a pair $(l_i,l_i)$, where $l_i=length(M^+_{\tau,i})=length(M^-_{\tau,i}))$ to the partition of
$\check\CO(\chi)$;
\item the entries
  $[l_i]=(-(l_i-1),-(l_i-1)+2,\dotsc,(l_i-1))$ to $\check h$ (this is
  the middle element of the principal orbit in $gl(l_i)$);
\item the entry $|\nu_i|$, where $\nu_i(1,\dotsc,1)=M^+_{\tau,i}-1/2[l_i]$
  to $\nu.$

\end{enumerate}

We give two examples of this process. For example, if $\tau=0$ and
$A_0=(0,0,1,1,1,1,2,3,3,4,5),$  then $M^+_{0,1}=(-5,-4,-3,-2,-1,0,1),$
$M^+_{0,2}=(-4,-3)$, and $M^+_{0,3}=(-1,0,1).$ Of course, always,
$M^-_{\tau,i}=-M^+_{\tau,i}.$ This means we add to the partition of
$\check\CO(\chi)$ the entries $(2,2,3,3,7,7).$ To $\check h$ we add
$(-1,1)$, $(-2,0,2)$ and $(-6,-4,-2,0,2,4,6),$ and to $\nu$ we add
$7/2,$ $0$, and $2.$

If $\tau=1/4,$ and $A_{1/4}=(1/4,1/4,3/4,5/4,5/4)$, then
$M^+_{\tau,1}=(-5/4,-1/4,3/4)$ and $M^+_{\tau,2}=(-5/4,-1/4).$ This
adds the entries $(2,2,3,3)$ to the partition of $\check\CO(\chi).$ To
$\check h$ we add $(-1,1)$ and $(-2,0,2)$, and to $\nu$ we add $3/4$
and $1/4.$

\smallskip

Now assume $\tau=1/2.$ Form again $A'_{1/2}=A_{1/2}\sqcup
(-A_{1/2})$. We only construct increasing strings $M^+_{1/2,i}$,
$i=1,\ell$ in this case. Remove the smallest entry $b$ in $A'_{1/2}$
and place it in $M^+_{1/2,1}$. Continue with $b+1,b+2,\dotsc$ until
this is not possible. This concludes the construction of
$M^+_{1/2,1}$. We repeat the process with the remaining entries in
$A'_{1/2},$ until we remove all of them. Any two strings
$M^+_{1/2,i}$ and $M^+_{1/2,j}$ such that $M^+_{1/2,i}=-M^+_{1/2,j}$
contribute:

\begin{enumerate}
\item a pair $(k_i,k_i)$, where
$k_i=length(M^+_{1/2,i})=length(M^+_{1/2,j})$ to $\check\CO(\chi)$;
\item the entries $[k_i]=(-(k_i-1)-(k_i-1)+2,\dotsc,(k_i-1))$ to
  $\check h$;
\item the entry $|\nu_i|,$ where
  $\nu_i(1,\dotsc,1)=M^+_{1/2,i}-1/2[k_i].$

\end{enumerate}

Remove them from the list of strings and repeat. If the number of
strings was odd, there is one remaining string at the end, say
$M^+_{1/2,k}.$ We call this string {\it distinguished}. The motivation
is that the positive part of this string is $1/2$ the middle
element of a distinguished nilpotent orbit in a symplectic complex
algebra as in \cite{CM}. Add the
partition corresponding to that orbit to $\check\CO(\chi).$ The
contribution to $\check h$ is twice the positive part of
$M^+_{1/2,k},$ while there is no contribution to $\nu.$

For example, if $A_{1/2}=(1/2,1/2,1/2,3/2,3/2,3/2,5/2,5/2,5/2,7/2),$
then we extract five strings:
$M^+_{1/2,1}=(-7/2,-5/2,-3/2,-1/2,1/2,3/2,5/2,7/2)$,
$M^+_{1/2,2}=M^+_{1/2,3}=(-5/2,-3/2,-1/2,1/2,3/2,5/2),$
$M^+_{1/2,4}=(5/2)$, and $M^+_{1/2,5}=(-5/2).$ Then the distinguished
string is $M^+_{1/2,1}.$ Its positive part is $(1/2,3/2,5/2,7/2)$,
which is $1/2$ the middle element of the principal nilpotent orbit in
$sp(8).$ It adds $(8)$ to $\check\CO(\chi),$ $(1,3,5,7)$ to $\check
h$, and nothing to $\nu.$ The other four strings add
$(6,6)$ and $(1,1)$ respectively to $\check\CO(\chi).$ Their
contribution to $\check h$ are $(-5,-3,-1,1,3,5)$, and $(0)$, and to
$\nu$, they contribute $0$, respectively $5/2.$

\smallskip

In conclusion, if our $\chi$ where the disjoint union
\begin{equation}\label{eq:4.4.2}
\chi=A_{0}\sqcup A_{1/4}\sqcup A_{1/2},
\end{equation}
with $A_\tau$ as above, then we rearrange the entries of the
nilpotent orbit increasingly, \eg
\begin{equation}\label{eq:4.4.3}
\check\CO(\chi)=(1,1;2,2,2,2;3,3,3,3;6,6;7,7;;8),
\end{equation}
 and $\check h$ and
$\nu$ are permuted accordingly. We separate the distinguished part by
$;;$ and the groups of identical entries by $;$. After this arrangement
\begin{equation}\label{eq:4.4.4}
\nu=(5/2;3/4,7/2;0,1/4;0;2;;\ ).
\end{equation}

\subsubsection{$G$ of type $C_n$} The algorithm of forming strings is
the same as the one for $B_n$ except that $A'_\tau=A_\tau\sqcup
(-A_\tau)\sqcup \{0\}),$ and the special case is $\tau=0$. For
$\tau=0$ we apply the algorithm as in the case $\tau=1/2$ for $B_n,$
and the distinguished string corresponds to a distinguished nilpotent
orbit in $so(2k+1,\bC).$ In the case $0<\tau\le 1/2$ the algorithm is
identical to $0\le\tau<1/2$ in $B_n.$

\subsubsection{$G$ of type $D_n$} The algorithm of forming strings is
the same as the one for $B_n$ except that  the special case is $\tau=0$. For
$\tau=0$ we apply the algorithm as in the case $\tau=1/2$ for $B_n,$
and the distinguished string corresponds to a distinguished nilpotent
orbit in $so(2k,\bC).$ In the case $0<\tau\le 1/2$ the algorithm is
identical to $0\le\tau<1/2$ in $B_n.$ There is a minor complication
when the parameter belongs to a very even nilpotent orbit which we ignore here. We refer the
reader to section 2.7 in \cite{Ba1} for the details of this case.

\subsection{Testing unitarity}\label{sec:4.5} Once the decomposition of a character
$\chi$ into strings is completed as in section \ref{sec:4.4}, testing
unitarity by theorem \ref{t:4.3.1} is easy.

Let assume that from a character (spherical parameter) $\chi$ we obtained in section
\ref{sec:4.4} the strings and the parameter $\nu$ for the centralizer
$\fz(\check e,\check h,\check f).$ Let us assume
\begin{align}\label{eq:4.5.1}
&\check\CO(\chi)=(\underbrace{1,\dotsc,1}_{2\ell_1};\underbrace{2,\dotsc,2}_{2\ell_2};\dots;\underbrace{k\dotsc k}_{2\ell_k};;[\lambda]),\\\notag
&\nu=(\nu_{1,1},\dotsc,\nu_{1,\ell_1};\nu_{2,1},\dotsc,\nu_{2,\ell_2};\dots;\nu_{k,1},\dotsc,\nu_{k,\ell_k};;\ ),
\end{align}
where $\ell_i\ge 0,$ $i=1,k$, and the conventions for notation are as
in (\ref{eq:4.4.3}) and (\ref{eq:4.4.4}). Moreover, we may permute the
entries in $\nu$ between any two consecutive $;$'s to be increasing.

The type of the centralizer $\fz(\check e,\check h,\check f)$ is well-known (see
\cite{CM} or \cite{Ca}). It is a product of types $B_{\ell_i}$,
$D_{\ell_i}$ or $C_{\ell_i},$ depending on the type of $G$. Then one
checks if the corresponding entries
$(\nu_{i,1},\dotsc,\nu_{i,\ell_i})$ in $\nu$ satisfy the conditions in
(\ref{eq:4.2.3}) for the corresponding (dual) type.

\smallskip

\noindent{\bf Example.} If $(\check\CO(\chi),\nu)$ are as in (\ref{eq:4.4.3}) and
$(\ref{eq:4.4.4})$, then $\check\fg=sp(32,\bC),$ and $\fz(\check
e,\check h,\check f)$ has type $C_1\times D_2\times C_2\times
D_1\times C_1$, in the same order as $\check\CO(\chi)$ is written, and
where by $D_1$ we mean a one-dimensional torus. We test the
corresponding $\nu_i$'s against (\ref{eq:4.2.3}) (for a torus $D_1$, the
only unitary parameter is $0$):

\begin{center}\begin{tabular}{|c|c|c|}
\hline
$\nu$ &centralizer &unitary?\\
\hline
$(5/2)$ &$C_1$      &no\\
$(3/4,7/2)$ &$D_2$  &no\\
$(0,1/4)$  &$C_2$   &yes\\
$(0)$      &$D_1$   &yes\\
$(2)$      &$C_1$   &no\\
\hline
\end{tabular}\end{center}

\noindent In conclusion, $\chi$ is not unitary.

\subsection{Maximal parabolic cases}\label{sec:4.6} In the next
sections we give a sketch of some of the ideas involved in the proofs
of the results stated in section \ref{sec:4.3}. We will be concerned  with theorems \ref{t:4.3.1}
and \ref{t:4.3.2}, which are proved in
\cite{Ba1} for classical split groups,\cite{Ci} for $G_2$ and $F_4$, and
\cite{BC2} for types $E$. The method for proving theorem
\ref{t:4.3.1}, but not the statement about relevant $W$-types, was
used for the first time in \cite{BM3} for classical split
groups.

Recall that we are in the setting of the affine graded Hecke
algebra $\bH=\bH_G.$ The method
consists in a double induction:

\begin{enumerate}
\item[(a)] an upward induction, by the rank of the group $G$, and
\item[(b)] a downward induction, in the closure ordering for nilpotent
  orbits $\check\CO$ in $\check\fg$.
\end{enumerate}

In this scheme, one determines the $0$-complementary series last. We
have seen in section \ref{sec:4.2} that there is an alternate method
for finding the $0$-complementary series directly.

As remarked in section \ref{sec:4.1}, there is nothing to do for
distinguished orbits $\check\CO.$ Therefore, the first basic cases to
do are when the nilpotent $\check\CO$ is parameterized in the
Bala-Carter classification by the Levi component of a maximal
parabolic subalgebra. We call them {\it maximal parabolic cases}.

Let us assume therefore that the Lie triple $\{\check e,\check
h,\check f\}$ of $\check\CO$ is contained in the Levi subalgebra
$\check\fm$ of a maximal parabolic subalgebra
$\check\fp=\check\fm+\check\fn\subset \check\fg.$ The graded Hecke
algebra corresponding to $\check\fm$, which we denote $\bH_M$, is
naturally a subalgebra of $\bH.$ Let $\chi=\check
h/2+\nu$ be the character as before. Note that $\check h/2$ defines a
special unipotent representation $L_M(\check h/2)$ for $M$.

\begin{lemma}\label{l:4.6}
With the notation as above, assume that $\nu$ is dominant with respect
to the roots in $\fn.$ Then $L(\chi)$ is the unique irreducible
quotient of
\begin{equation}\label{eq:4.6.1}
X(M,\check h/2,\nu):=\Ind_{\bH_M}^{\bH_G}(L_M(\check h/2)\otimes\bC_\nu).
\end{equation}
\end{lemma}

%\begin{proof}
We refer the reader to \cite{Ba1} (also \cite{BM3}) for the proof of
this result. It is proved using the Iwahori-Matsumoto involution and
the geometric classification of $\bH$-modules as it follows from
\cite{KL} and \cite{L1}.
%\end{proof}

One knows that the module $L(\chi)$ is Hermitian if and only if
there exits $w\in W$ such that
\begin{equation}\label{eq:4.6.2}
wM=M,\quad w\check h=\check h,\quad \text{and } w\nu=-\nu.
\end{equation}
If there exists such a $w$, we choose $w_m$ to be a minimal element in
the double coset $W(M)wW(M).$ Then one can define a generalization of
the intertwining operators from section \ref{sec:3.2}:
\begin{align}\label{eq:4.6.3}
&A(w_m,\nu):X(M,\check h/2,\nu)\longrightarrow X(M,\check
h/2,-\nu),\quad\text{ and}\\\notag
&A_\psi(w_m,\nu):\Hom_{W(M)}[V_\psi:L_M(\check h/2)]\longrightarrow \Hom_{W(M)}[V_\psi:L_M(\check h/2)],
\end{align}
for every $W$-type $(\psi,V_\psi)$.

In section 6 of \cite{Ba1}, these operators are computed explicitly for all
relevant $W$-types in the classical groups as in definition
\ref{d:3.4}. One can reduce the calculation, so that the only cases
that one considers there are for $\check\fh$ the middle element of the
principal nilpotent orbit on $\check\fm.$ In those cases, the
dimension of the $\Hom$ spaces in equation (\ref{eq:4.6.3}) is always
$1$. That makes $A_\psi(w_m,\nu)$ a scalar, which is normalized so
that it is
$+1$ when $\psi=triv.$

\smallskip

\noindent{\bf Example 4.6.1.} Let us assume that $G$ is of type $B_{n+k}.$ Then
$\check\fg=sp(2n+2k,\bC),$ and we consider $\check\CO=(k,k,2n)$,
sothat $\chi=(\frac 12,\frac 32,\dots,n-\frac
12,-\frac{k-1}2+\nu,\dots,\frac{k-1}2+\nu),$ and  $\check\fm=sp(2n,\bC)\times
gl(k,\bC).$ The calculation in \cite{Ba1} gives:

\begin{center}
\begin{tabular}{|c|c|}
\hline
relevant $W$-type &$A_\psi(w_m,\nu)$\\
\hline
$(n+k-m)\times (m)$ &$\displaystyle{\prod_{0\le j\le
    m-1}\frac{n+\frac {k}2-j-\nu}{n+\frac{k}2-j+\nu}}$\\
\hline
$(m,n+k-m)\times (0)$ &$\displaystyle{\prod_{0\le j\le m-1}
  \frac{n+\frac{k}2-j-\nu}{n+\frac{k}2-j+\nu}\cdot \frac{\frac{k}2-j-\nu} {\frac{k}2-j+\nu}
         }$\\
\hline
\end{tabular},
\end{center}
where $0\le m\le k.$

For example, when $n=1,$ $k=1,$ we are in the
case $sp(2)\times gl(1)\subset sp(4).$ The only $W$-types in the
$\Ind_{\bH_{C_1}}^{\bH_{C_2}}(triv\otimes\bC_\nu)$ are $2\times 0,$
$1\times 1$ and $11\times 0.$ Then $\chi=(\nu,\frac 12),$ $\nu>0.$ The
operators are $+1$, $\frac {3/2-\nu}{3/2+\nu}$, and $\frac
{(1/2-\nu)(3/2-\nu)}{(1/2+\nu)(3/2+\nu)}$ respectively. Since
  $X(M,\check h/2,\nu)$ is irreducible at $\nu=0$ in this case, we
  deduce that, in $sp(4,\bC)$, $CS^\bH((1,1,2))=\{\chi=(\nu,1/2):
  0\le\nu<1/2\}$, which is identical with the $0$-complementary series
  of the centralizer of type $A_1$ of the nilpotent $(1,1,2).$

On the other hand, if one considers $n=0,$ $k=2,$ we are in the case
$gl(2)\subset sp(4).$ The only $W$-types in the
$\Ind_{\bH_{A_1}}^{\bH_{C_2}}(triv\otimes\bC_\nu)$ are $2\times 0,$
$1\times 1$ and $0\times 2.$ Then $\chi=(-\frac 12+\nu,\frac 12+\nu),$ $\nu>0.$ The
operators are $+1$, $\frac {1-\nu}{1+\nu}$, and $-\frac
{(1-\nu)}{(1+\nu)}$ respectively. This shows $X(M,\check h/2,\nu)$,
$\nu>0$ if irreducible, it must be nonunitary. In this case, we know that
  $X(M,\check h/2,\nu)$ is reducible at $\nu=0$. We
  deduce that, in $sp(4,\bC)$, $CS^\bH((2,2))=\{\chi=(-1/2+\nu,1/2+\nu):
  \nu=0\}$, which is the same as the unitary set for the centralizer,
  which is of type $D_1,$ \ie a one-dimensional torus.

\smallskip

For exceptional groups, it is not true anymore that the relevant
$W$-types appear with multiplicity $1$, even if we induce from the
trivial representation of $\bH_M.$ Explicit calculations of the
operators $A_\psi(w_m,\nu)$ are done in \cite{BC2}.

\smallskip

\noindent{\bf Example 4.6.2.} Let us consider $G$ of type $E_8$,
$\check\CO=E_6+A_1$, so  that $\check\fm=E_6+A_1.$ The character is
$\chi=(0,1,2,3,4,-\frac 92,-\frac 72,4)+\nu\omega_7$, $\nu>0,$ where
$\omega_7$ is the fundamental weight corresponding to the coroot
$\check\al_7.$ (Recall that we use Bourbaki's notation.) In this
example, computing the determinants of $A_\psi(w_m,\nu)$ turns out to
be sufficient. we give below the table with the relevant $W$-types
which appear in this case.

\begin{center}
\begin{tabular}{|c|c|c|}
\hline
$W$-type &$\dim \Hom$-space &Determinant of $A_\psi(w_m,\nu)$\\
\hline

$8_z$ &$1$ &\begin{tiny}$\displaystyle{\frac {\frac {19}2-\nu}{\frac {19}2+\nu}}$\end{tiny}\\
\hline
$35_x$ &$3$  &\begin{tiny}$\displaystyle{\frac {(\frac {19}2-\nu)^2(\frac
    {17}2-\nu)(\frac {11}2-\nu)(\frac 92-\nu)(\frac 52-\nu)}{(\frac
    {19}2+\nu)^2(\frac {17}2+\nu)(\frac {11}2+\nu)(\frac 92+\nu)(\frac
    52+\nu)}}$\end{tiny}\\
\hline
$112_z$ &$2$ &\begin{tiny}$\displaystyle{\frac {(\frac {19}2-\nu)^2(\frac
    {17}2-\nu)(\frac {11}2-\nu)^2(\frac 92-\nu)(\frac 52-\nu)(\frac
    32-\nu)}{(\frac {19}2+\nu)^2(\frac {17}2+\nu)(\frac
    {11}2+\nu)^2(\frac 92+\nu)(\frac 52+\nu)(\frac 32+\nu)}}$\end{tiny}\\
\hline
$84_x$ &$2$  &\begin{tiny}$\displaystyle{\frac {(\frac {19}2-\nu)^2(\frac
    {17}2-\nu)(\frac {11}2-\nu)^2(\frac 92-\nu)(\frac 52-\nu)(\frac
    32-\nu)(\frac 12-\nu)}{(\frac {19}2+\nu)^2(\frac {17}2+\nu)(\frac
    {11}2+\nu)^2(\frac 92+\nu)(\frac 52+\nu)(\frac 32+\nu)(\frac
    12+\nu)}}$\end{tiny}\\
\hline
\end{tabular}
\end{center}

\smallskip

We plot the signatures of these determinants in a table:

\smallskip

\begin{center}
\begin{tabular}{|l|lllllllllllllll|}
\hline
$\nu$ &   &$\frac 12$ & &$\frac 32$  & &$\frac {5}2$  &  &$\frac {9}2$
& &$\frac {11}2$  & &$\frac {17}2$ & &$\frac {19}2$ &  \\
\hline
$8_z$ &$+$ &$+$ &$+$ &$+$ &$+$  &$+$ &$+$ &$+$  &$+$ &$+$ &$+$  &$+$  &$+$ &$0$ &$-$ \\
\hline
$35_x$ &$+$ &$+$ &$+$  &$+$ &$+$  &$0$  &$-$ &$0$ &$+$ &$0$ &$-$ &$0$ &$+$ &$0$ &$+$ \\
\hline
$84_x$ &$+$ &$0$ &$-$  &$0$ &$+$  &$0$ &$-$ &$0$ &$+$ &$0$ &$+$ &$0$ &$-$ &$0$ &$-$ \\
\hline
$112_z$ &$+$ &$+$ &$+$  &$0$ &$-$  &$0$ &$+$ &$0$ &$-$ &$0$ &$-$ &$0$ &$+$ &$0$ &$+$ \\
\hline
\end{tabular}
\end{center}

\smallskip

Since $X(M,\check h/2,\nu)$ is irreducible at $\nu=0,$ we conclude
that, in $E_8,$ $CS(E_6+A_1)=\{\chi: 0\le\nu<1/2\},$ which is
identical with the $0$-complementary series of the centralizer which
is type $A_1.$

\smallskip

An example of a maximal parabolic nilpotent orbit with torus centralizer in $E_8$, and
therefore with the unitary set formed only of $\nu=0$, is
$\check\CO=D_5+A_2.$ Similar tables as above are available in that
case, but we skip the details here.

\medskip

\noindent{\bf Remarks.} There are three important remarks to be made
which were assumed implicitly in the examples above:

\begin{enumerate}
\item The reducibility points of the induced module $X(M,\check
  h/2,\nu)$ for $\nu>0$ are known a priori, by different
  methods (see \cite{BC} and the references therein). Therefore, we
  can compare with the explicit calculations of operators
  $A_\psi(w_m,\nu)$, and see that the relevant $W$-types do not miss
  any reducibility points.
\item The reducibility of $X(M,\check h/2,\nu),$ at $\nu=0$ is known  by geometric considerents (\cite{KL}).
\item Whenever $X(M,\check h/2,\nu)$ is reducible at $\nu>0,$ the
  spherical quotient is parameterized by a larger nilpotent
  $\check\CO'\supsetneq\check\CO.$ Therefore, whenever one is
  concerned with the orbit $\check\CO,$ these reducibility points do
  not need to be considered. The unitarity of the corresponding
  modules was {\it already} checked in our inductive procedure.
\end{enumerate}

We summarize now the main consequence of these types of calculations.

\begin{proposition}[\cite{Ba1},\cite{BC2}] Assume $\check\CO$ is a
  maximal parabolic nilpotent orbit. Then a spherical module
  parameterized by $\chi$ such that $\check\CO(\chi)=\check\CO$ is
  unitary if and only if the intertwining operators are positive
  semidefinite on the relevant $W$-types.

\end{proposition}

We remark that one can obtain the unitarity results for maximal
parabolic cases more easily if the constraint of working with relevant
$W$-types only is removed. This is the approach of \cite{BM3} and
\cite{BC}. In there, one obtains the complementary series by checking
the signature of two $W$-types which are Springer representations for
$\check\CO$ and for an orbit $\check\CO'$ consecutive to $\check\CO$
in the closure ordering.

\subsection{Induction}\label{sec:4.7} In this section, we exemplify
the inductive step in the proof. We are still in the setting of the
affine graded Hecke algebra. The main idea is as follows: assume we
start with a parameter $\chi=\check h/2+\nu$ associated to $\check\CO$,
where $\check\CO$ is not a maximal parabolic nilpotent. Let $M$ and
$X(M,\check h/2,\nu)$ be as in lemma \ref{l:4.6}. One knows explicitly
for which values of $\nu$ the module $X(M,\check h/2,\nu)$ is
reducible. If $\nu=\nu_1$ is such a value (and $\nu_1$ is dominant with
respect to the roots in $\check\fn$), then the spherical module
$L(\chi_1),$ where $\chi_1=\check h/2+\nu_1$, is attached to a
nilpotent orbit $\CO_1$, with the property that
\begin{equation}\label{sec:4.7.1}
\check\CO_1\supset\check\CO\text{ and }\check\CO_1\neq\check\CO.
\end{equation}
By induction, we already know if $L(\chi_1)$ is unitary or
not.

\begin{lemma}[1]\label{l:4.7.1} Let $\chi(t)=\check h/2+\nu(t):[0,1]\to \check\fh$ be a continuous
  function, such that $\chi(t)=\chi_0$ and $\chi(1)=\chi_1.$ Assume
  that $X(M,\check h/2,\nu(t))$ is irreducible for $0\le t<1$ and
  $X(M,\check h/2,\nu_1)$ is reducible.

If
  $L(\chi_1)$ is not unitary, then $L(\chi_0)=X(M,\check h/2,\nu_0)$
  is not unitary as well.
\end{lemma}

This well-known criterion for nonunitarity can be applied therefore,
and it is the main tool for ruling out nonunitary parameters. Of
course, in general, there are many delicate combinatorial issues that
arise; in the classical cases, one needs to choose carefully how to
deform parameters $\chi=\check h/2+\nu.$ In the exceptional groups,
another complication arises: there are cases (\eg
in $\check\CO=2A_2\subset E_8$) of parameters $\chi=\chi_0$, which
turn out to be nonunitary, but in all possible deformations as in
lemma \ref{l:4.7.1}(1), the modules $L(\chi_1)$ are unitary. In those
cases we apply some ad-hoc signature arguments.

\medskip

The second (again well-known) criterion which we employ is the
following {\it complementary series} method.

\begin{lemma}[2]\label{l:4.7.2}Let $\chi(t)=\check h/2+\nu(t):[0,1]\to \check\fh$ be a continuous
  function, such that $\chi(t)=\chi_0$ and $\chi(1)=\chi_1.$ Assume
  that

\begin{enumerate}
\item[(a)]$X(M,\check h/2,\nu(t))$ is Hermitian and irreducible for $0\le t\le 1$, and
\item
  $X(M,\check h/2,\nu_1)$ is unitarily induced, \ie $X(M,\check h/2,\nu_1)=\Ind_{\bH_{M'}}^{\bH_{G}}(V_{M'}),$ where
  $M'$ is a Levi component,  $V_{M'}$ is a Hermitian spherical
  module for $\bH_{M'}$.
\end{enumerate}

Then
  $L(\chi_0)=X(M,\check h/2,\nu_0)$ is unitary in $\bH_G$ if and only
  if $V_{M'}$ is unitary in $\bH_{M'}.$
\end{lemma}

Since, by induction the unitarity of $\bH_{M'}$-modules is known, we
can apply this criterion. This is our main criterion for proving the
unitarity of modules $L(\chi).$

\medskip

\begin{figure}[h]
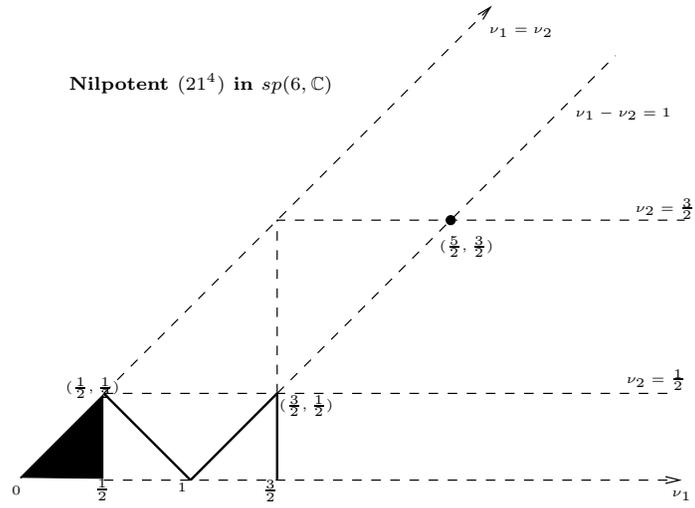
\label{fig:C3}
\input C3_21111.pstex_t
\caption{Example: $\check\CO=(21^4)$ for the Hecke algebra of type $C_3$}
\end{figure}

\noindent{\bf Example.} To illustrate this discussion, we conclude the section with a simple
example. Consider $\check\CO=(1^4,2)$ in $sp(6,\bC),$ so
$\check\fm=gl(1)^2\times sp(2).$ We write
$\chi=(\nu_1,\nu_2,1/2),$ where by conjugation with the Weyl group, we
can assume that $0\le\nu_2\le\nu_1.$ The centralizer has type $C_2$
in this case. The lines where $X(M,\check h/2,(\nu_1,\nu_2))$ becomes
reducible are drawn in figure \ref{fig:C3}.
In this picture, the solid lines denote unitary spherical parameters,
while the dotted lines denote nonunitary parameters. The reducibility
lines (in the Hecke algebra case) are:

\begin{enumerate}

\item[$\bullet$] $\nu_1=1/2,\ \nu_2=1/2$: $\CO(\chi)=(2211)$;
\item[$\bullet$] $\nu_1=3/2,\ \nu_2=3/2$: $\CO(\chi)=(411)$;
\item[$\bullet$] $\nu_1\pm\nu_2=1$: $\CO(\chi)=(222)$.
\end{enumerate}

On the lines $\nu_2=0$, and $\nu_1=\nu_2$, the module is unitarily
induced from a spherical module parameterized by $(211)$ in
$sp(4,\bC)$, respectively $(11)$ in $gl(2,\bC).$

\medskip

\newpage

\section{Lists of unitary spherical parameters}\label{sec:5}

In this section, we give tables with the explicit description of the
spherical unitary parameters for split $p$-adic groups as in
theorem \ref{t:4.3.1}. By theorems \ref{t:4.3.3} and \ref{t:4.3.4},
these lists also describe the spherical unitary dual for split real
classical groups, and, by \cite{V3}, also for $G_2$. For real split
$F_4,E_6,E_7,E_8$, the tables represent the only spherical parameters
which {could be} unitary, and it is natural to expect that,
in fact, this is the spherical unitary dual in these cases as well.
% The references for the tables of this section are
%\cite{Ba1}, \cite{Ci}, and \cite{BC}.

\subsection{Examples of classical groups}

\subsubsection{$A_4$}\label{sec:A4}

\begin{small}
\begin{center}

\begin{longtable}{|c|c|c|}
\caption{Table of spherical unitary parameters for $A_4$}\label{table:A4}\\
\hline
\multicolumn{1}{|c|}{$\mathbf{\check\CO}$} &
\multicolumn{1}{c|}{$\mathbf\chi$}
& \multicolumn{1}{c|}{unitary $\nu$} \\ \hline
\endfirsthead

\multicolumn{3}{c}%
{{  \tablename\ \thetable{} -- continued from previous page}} \\
\hline
\multicolumn{1}{|c|}{$\mathbf{\check\CO}$} &
\multicolumn{1}{c|}{$\mathbf\chi$}
& \multicolumn{1}{c|}{unitary $\nu$}
 \\ \hline
\endhead

\hline \hline
\endlastfoot

$(5)$ &$(-2,-1,0,1,2)$  &\\
\hline
$(41)$ &$(-3/2,-1/2,0,1/2,3/2)$ &\\
\hline
$(32)$ &$(-1,-1/2,0,1/2,1)$   &\\
\hline
$(311)$ &$(-1,0,1,-\nu,\nu)$ &$0\le\nu<1/2$\\
\hline
$(221)$ &$(-1/2-\nu,1/2-\nu,0,-1/2+\nu,1/2+\nu)$ &$0\le\nu<1/2$\\
\hline
$(21^3)$ &$(-1/2,0,1/2,-\nu,\nu)$ &$0\le\nu<1/2$\\
\hline
$(1^5)$ &$(-\nu_2,-\nu_1,0,\nu_1,\nu_2)$ &$0\le\nu_2\le\nu_2<1/2$\\

\end{longtable}
\end{center}
\end{small}

%\newpage

\subsubsection{$B_4$}\label{sec:B4}

\begin{small}
\begin{center}
\begin{longtable}{|c|c|c|}
\caption{Table of spherical unitary parameters for $B_4$}\label{table:B4}\\
\hline
\multicolumn{1}{|c|}{$\mathbf{\check\CO}$} &
\multicolumn{1}{c|}{$\mathbf\chi$}
& \multicolumn{1}{c|}{unitary $\nu$} \\ \hline
\endfirsthead

\multicolumn{3}{c}%
{{  \tablename\ \thetable{} -- continued from previous page}} \\
\hline
\multicolumn{1}{|c|}{$\mathbf{\check\CO}$} &
\multicolumn{1}{c|}{$\mathbf\chi$}
& \multicolumn{1}{c|}{unitary $\nu$}
 \\ \hline
\endhead

\hline \hline
\endlastfoot

$(8)$ &$(1/2,3/2,5/2,7/2)$  &\\
\hline
$(62)$  &$(1/2,1/2,3/2,5/2)$   &\\*
\hline
$(611)$ &$(\nu,1/2,3/2,5/2)$  &$0\le\nu<1/2$\\
\hline
$(44)$  &$(-3/2+\nu,-1/2+\nu,1/2+\nu,3/2+\nu)$ &$\nu=0$\\
\hline

$(422)$ &$(-1/2+\nu,1/2+\nu,1/2,3/2)$ &$0\le\nu<1$\\
\hline

$(4211)$ &$(\nu,1/2,1/2,3/2)$  &$0\le\nu<1/2$\\
\hline

$(41^4)$ &$(\nu_1,\nu_2,1/2,3/2)$ &$0\le\nu_1\le\nu_2<1/2$\\
\hline

$(332)$ &$(-1+\nu,\nu,1+\nu,1/2)$ &$0\le\nu<1/2$\\
\hline

$(3311)$ &$(-1+\nu_1,\nu_1,1+\nu_1,\nu_2)$ &$0\le\nu_1<1/2,0\le\nu_2<1/2$\\
\hline

$(2^4)$ &$(-1/2+\nu_1,1/2+\nu_1,-1/2+\nu_2,1/2+\nu_2)$ &$0\le\nu_1\le\nu_2<1/2$\\

\hline
$(2^311)$ &$(-1/2+\nu_1,1/2+\nu_1,\nu_2,1/2)$ &$0\le\nu_1<1/2,0\le\nu_2<1/2$\\
\hline
$(221^4)$ &$(\nu_1,\nu_2,-1/2+\nu_3,1/2+\nu_3)$ &$\nu_3=0,~0\le\nu_1\le\nu_2<1/2$\\
\hline
$(21^6)$ &$(\nu_1,\nu_2,\nu_3,1/2)$  &$0\le\nu_1\le\nu_2\le\nu_3<1/2$\\
\hline
$(1^8)$ &$(\nu_1,\nu_2,\nu_3,\nu_4)$
&$0\le\nu_1\le\nu_2\le\nu_3\le\nu_4<1/2$\\

\end{longtable}
\end{center}
\end{small}

\subsubsection{$C_4$}\label{sec:C4} 

\begin{small}
\begin{center}
\begin{longtable}{|c|c|c|}
\caption{Table of spherical unitary parameters for $C_4$}\label{table:C4}\\
\hline
\multicolumn{1}{|c|}{$\mathbf{\check\CO}$} &
\multicolumn{1}{c|}{$\mathbf\chi$}
& \multicolumn{1}{c|}{unitary $\nu$} \\ \hline
\endfirsthead

\multicolumn{3}{c}%
{{  \tablename\ \thetable{} -- continued from previous page}} \\
\hline
\multicolumn{1}{|c|}{$\mathbf{\check\CO}$} &
\multicolumn{1}{c|}{$\mathbf\chi$}
& \multicolumn{1}{c|}{unitary $\nu$}
 \\ \hline
\endhead

\hline \hline
\endlastfoot

$(9)$ &$(1,2,3,4)$ &\\

\hline
$(711)$ &$(\nu,1,2,3)$ &$\nu=0$\\

\hline
$(531)$ &$(0,1,1,2)$  &\\

\hline
$(522)$ &$(-1/2+\nu,1/2+\nu,1,2)$ &$0\le\nu<1/2$\\

\hline
$(441)$ &$(-3/2+\nu,-1/2+\nu,1/2+\nu,3/2+\nu)$ &$0\le\nu<1/2$\\

\hline
$(51^4)$ &$(\nu_1,\nu_2,1,2)$ &$0\le\nu_1\le\nu_2<1-\nu_1$   \\

\hline
$(333)$ &$(-1+\nu,\nu,1+\nu,1)$ &$0\le\nu<1$ \\

\hline
$(331^3)$ &$(-1+\nu_1,\nu_1,1+\nu_1,\nu_2)$  &$\nu_1=0,~0\le\nu_2<1$  \\

\hline
$(32211)$ &$(-1/2+\nu_1,1/2+\nu_1,\nu_2,1)$  &$\nu_2=0,~0\le\nu_1<1/2$  \\

\hline
$(31^6)$  &$(\nu_1,\nu_2,\nu_3,1)$           &$\nu_3=0,~0\le\nu_1\le\nu_2<1-\nu_1$   \\

\hline
$(2^41)$  &$(-1/2+\nu_1,1/2+\nu_1,-1/2+\nu_2,1/2+\nu_2)$ &$0\le\nu_1\le\nu_2<1/2$  \\

\hline
$(221^5)$ &$(-1/2+\nu_1,1/2+\nu_1,\nu_2,\nu_3)$  &$0\le\nu_1<1/2,~0\le\nu_2<\nu_3<1-\nu_2$  \\

\hline
$(1^9)$   &$(\nu_1,\nu_2,\nu_3,\nu_4)$
&$0\le\nu_1\le\nu_2\le\nu_3\le\nu_4\le 1-\nu_3$  \\
          &
&$0\le\nu_1\le\nu_2\le\nu_3<1-\nu_2<\nu_4<1-\nu_1$\\

\end{longtable}
\end{center}
\end{small}

%\newpage

\subsubsection{$D_4$}\label{sec:D4} 

\begin{small}
\begin{center}
\begin{longtable}{|c|c|c|}
\caption{Table of spherical unitary parameters for $D_4$}\label{table:D4}\\
\hline
\multicolumn{1}{|c|}{$\mathbf{\check\CO}$} &
\multicolumn{1}{c|}{$\mathbf\chi$}
& \multicolumn{1}{c|}{unitary $\nu$} \\ \hline
\endfirsthead

\multicolumn{3}{c}%
{{  \tablename\ \thetable{} -- continued from previous page}} \\
\hline
\multicolumn{1}{|c|}{$\mathbf{\check\CO}$} &
\multicolumn{1}{c|}{$\mathbf\chi$}
& \multicolumn{1}{c|}{unitary $\nu$}
 \\ \hline
\endhead

\hline \hline
\endlastfoot

$(71)$ &$(0,1,2,3)$ & \\
\hline
$(53)$ &$(0,1,1,2)$ &  \\
\hline
$(51^3)$ &$(\nu,0,1,2)$ &$0\le\nu<1$\\
\hline
$(44)_+$ &$(-3/2+\nu,-1/2+\nu,1/2+\nu,3/2+\nu)$ &$0\le\nu<1/2$\\
\hline
$(44)_-$ &$(3/2-\nu,-1/2+\nu,1/2+\nu,3/2+\nu)$ &$0\le\nu<1/2$\\
\hline
$(3311)$ &$(\nu_1,-1+\nu_2,\nu_2,1+\nu_2)$ &$\nu_1=0,~\nu_2=0$ \\
\hline
$(31^5)$ &$(\nu_1,\nu_2,0,1)$ &$0\le\nu_1\le\nu_2<1-\nu_1$  \\
\hline
$(2^4)_+$ &$(-1/2+\nu_1,1/2+\nu_1,-1/2+\nu_2,1/2+\nu_2)$ &$0\le\nu_1,\nu_2<1/2$\\
\hline
$(2^4)_-$ &$(1/2-\nu_1,1/2+\nu_1,-1/2+\nu_2,1/2+\nu_2)$  &$0\le\nu_1,\nu_2<1/2$\\
\hline
$(221^4)$ &$(-1/2+\nu_1,1/2+\nu_2,\nu_3,\nu_4)$ &$0\le\nu_1<1/2,~0\le|\nu_2|\le\nu_3<1-|\nu_2|$\\
\hline
$(1^8)$  &$(\nu_1,\nu_2,\nu_3,\nu_4)$ &$0\le|\nu_1|\le\nu_2\le\nu_3\le\nu_4<1-\nu_3$\\
         &
&$0\le|\nu_1|\le\nu_2\le\nu_3<1-\nu_2<\nu_4<1-|\nu_1|$\\

\end{longtable}
\end{center}
\end{small}

\subsection{Exceptional groups}\label{tables}
\begin{comment}
We record  the explicit description of the unitary spherical
dual for split exceptional groups over $p$-adic fields (\cf theorem
\ref{t:4.3.1}). For the corresponding split real groups, these lists represent
the only spherical parameters which {\it could} be
unitary. Conjecturally, all of them are unitary for the real group as
well, but this is only known for $G_2$, \cf \cite{V3}.
\end{comment}
The tables contain the
nilpotent orbits $\check\CO\neq 0$ in the notation of \cite{Ca} (the
case $\check\CO=0$ was recorded in section \ref{sec:zero}), the Hermitian infinitesimal
character, and the coordinates and type of the centralizer.

The nilpotent orbits which are exceptions are marked with $*$ in the
tables. The complementary series for them are listed after the tables.
For the rest of the nilpotent orbits, an \ic $\chi$ is in the
complementary series if and only if the corresponding parameter $\nu$
is in the $0$-complementary series for $\fz(\check\CO).$ The parameter $\nu$
is given by a string $(\nu_1,\dots,\nu_k)$, and the order agrees with
the way the centralizer $\fz(\check\CO)$ is written in the tables. The parts of
$\nu$ corresponding to a torus $T_1$ or $T_2$  in $\fz(\check\CO)$ must be
$0$, in order for $\chi$ to be unitary. In addition, if $\nu$
corresponds to $A_1$, the complementary series is $0\le\nu<\frac 12,$
while the notation $A_1^\ell$ means that it is $0\le\nu<1.$ If a string
$(\nu_1,\dots,\nu_k)$ of $\nu$ corresponds to type $A_k,$ the last
$k-[\frac k2]$ entries must be $0$ in order for $\chi$ to be
unitary. For example, in the table for $E_8$, for the nilpotent
$A_4+A_1$, the $\nu$-string is $(\nu_1,\nu_2,\nu_3)$ and the
centralizer is $A_2+T_1.$ This means that the unitary parameters are
those for which $\nu_3=0$ (this is the $T_1$-piece), $\nu_2=0$ and
$0\le\nu_1<\frac 12$ (this is the $0$-complementary series of $A_2$).

There is one difference in $E_6$ due to the fact that we only
consider Hermitian spherical infinitesimal characters $\chi.$ In
this table, the $\nu$-string already refers to the semisimple and
Hermitian spherical parameter of the centralizer. For example, the
nilpotent $A_2+A_1$ in $E_6$ has centralizer $A_2+T_1,$ and the
corresponding $\chi$ has a single $\nu$. This $\nu$ corresponds to
the Hermitian parameter in the $A_2$ part of $\fz(\check\CO)$, so it must
satisfy $0\le\nu<\frac 12.$

\subsubsection{$G_2$}
\begin{small}
\begin{center}
\begin{longtable}{|c|c|c|}
\caption{Table of parameters $(\check\CO,\nu)$ for $G_2$}\label{table:G2}\\
\hline
\multicolumn{1}{|c|}{$\mathbf{\check\CO}$} &
\multicolumn{1}{c|}{$\mathbf\chi$}
& \multicolumn{1}{c|}{$\mathbf{\fz(\check\CO)}$} \\ \hline
\endfirsthead

\multicolumn{3}{c}%
{{  \tablename\ \thetable{} -- continued from previous page}} \\
\hline
\multicolumn{1}{|c|}{$\mathbf{\check\CO}$} &
\multicolumn{1}{c|}{$\mathbf\chi$}
& \multicolumn{1}{c|}{$\mathbf{\fz(\check\CO)}$}
 \\ \hline
\endhead

\hline \hline
\endlastfoot

$G_2$ &$(1,2,-3)$ &$1$\\
\hline
$G_2(a_1)$ &$(0,1,-1)$&$1$\\
\hline
$\wti A_1$ &$(1,-\frac 12+\nu,-\frac 12-\nu)$ &$A_1$\\
\hline
$A_1$  &$(-\frac 12+\nu,\frac 12+\nu,-2\nu)$ &$A_1$\\
\hline
\end{longtable}
\end{center}
\end{small}

\subsubsection{$F_4$}

\begin{small}
\begin{center}
\begin{longtable}{|c|c|c|}
\caption{Table of parameters $(\check\CO,\nu)$ for $F_4$}\label{table:F4}\\
\hline
\multicolumn{1}{|c|}{$\mathbf{\check\CO}$} &
\multicolumn{1}{c|}{$\mathbf\chi$}
& \multicolumn{1}{c|}{$\mathbf{\fz(\check\CO)}$} \\ \hline
\endfirsthead

\multicolumn{3}{c}%
{{  \tablename\ \thetable{} -- continued from previous page}} \\
\hline
\multicolumn{1}{|c|}{$\mathbf{\check\CO}$} &
\multicolumn{1}{c|}{$\mathbf\chi$}
& \multicolumn{1}{c|}{$\mathbf{\fz(\check\CO)}$}
 \\ \hline
\endhead

%\hline \multicolumn{4}{|r|}{{Continued on next page}} \\ \hline
%\endfoot

\hline \hline
\endlastfoot

${F_4}$ &$(\frac {11}2,\frac 52,\frac 32,\frac 12)$ &$1$\\
\hline
${F_4(a_1)}$ &$(\frac 72,\frac 32,\frac 12,\frac 12)$ &$1$\\
\hline
${F_4(a_2)}$ &$(\frac 52,\frac 32,\frac 12,\frac 12)$ &$1$\\
\hline
${C_3}$ &$(\nu,\frac 52,\frac 32,\frac 12)$ &$A_1$ \\
\hline
${B_3}$ &$(\frac 32+\nu,-\frac 32+\nu,\frac 32,\frac 12)$ &$A_1^\ell$\\
\hline
${F_4(a_3)}$ &$(\frac 32,\frac 12,\frac 12,\frac 12)$ &$1$\\
\hline
${C_3(a_1)}$ &$(\nu,\frac 32,\frac 12,\frac 12)$ &$A_1$ \\
\hline
${A_1+\wti A_2}$ &$(\frac 14+\frac {3\nu}2,\frac 34+\frac {\nu}2,-\frac 14+\frac {\nu}2,-\frac 54+\frac {\nu}2)$ &$A_1$\\
\hline
${B_2}$ &$(\nu_1,\nu_2,\frac 32,\frac 12)$ &$2A_1$\\
\hline
${\wti A_1+A_2}$ &$(\frac 12+2\nu,\nu,-1+\nu,\frac 12)$ &$A_1$\\
\hline
${\wti A_2}$ &$(\nu_2+\frac{3\nu_1}2,1+\frac{\nu_1}2,\frac {\nu_1}2,-1+\frac {\nu_1}2)$ &$G_2$\\
\hline
${A_2}$ &$(\frac 12+\nu_1+\nu_2,-\frac 12+\nu_1,-\frac 12+\nu_2,\frac 12)$ &$A_1+T_1$\\
\hline
$*{A_1+\wti A_1}$ &$(\nu_1,\frac 12+\nu_2,-\frac 12+\nu_2,\frac 12)$ &$A_1+A_1^\ell$\\
\hline
${\wti A_1}$ &$(\nu_1+\nu_2,\nu_1-\nu_2,\frac 12+\nu_3,-\frac 12+\nu_3)$ &$B_2+T_1$ \\
\hline
${A_1}$ &$(\nu_1,\nu_2,\nu_3,\frac 12)$ &$C_3$\\

\end{longtable}
\end{center}
\end{small}

\smallskip

\noindent{\bf $\mathbf{F_4}$ exception}:

\noindent $\mathbf{A_1+\wti A_1}$: $\{\nu_1+2\nu_2<\frac
32,\ \nu_1<\frac 12\}\cup\{2\nu_2-\nu_1>\frac 32,\nu_2<1\}.$

\subsubsection{$E_6$}

\begin{small}
\begin{center}
\begin{longtable}{|c|c|c|}
\caption{Table of Hermitian parameters $(\check\CO,\nu)$ for $E_6$}\label{table:E6}\\
\hline
\multicolumn{1}{|c|}{$\mathbf{\check\CO}$} &
\multicolumn{1}{c|}{$\mathbf\chi$}
& \multicolumn{1}{c|}{$\mathbf{\fz(\check\CO)}$} \\ \hline
\endfirsthead

\multicolumn{3}{c}%
{{  \tablename\ \thetable{} -- continued from previous page}} \\
\hline
\multicolumn{1}{|c|}{$\mathbf{\check\CO}$} &
\multicolumn{1}{c|}{$\mathbf\chi$}
& \multicolumn{1}{c|}{$\mathbf{\fz(\check\CO)}$}
 \\ \hline
\endhead

%\hline \multicolumn{4}{|r|}{{Continued on next page}} \\ \hline
%\endfoot

\hline \hline
\endlastfoot

%\begin{small}
%\begin{table}\label{table:E6}
%\caption{Table of Hermitian spherical parameters for $E_6$}
%\noindent\begin{longtable}{|c|c|c|c|}
%\hline
%$\mathbf\check\CO$ &$\mathbf\chi$ &$\mathbf\nu$ &$\mathbf{\fz(\check\CO)}$\\
%\hline
%\hline

$E_6$ &$(0,1,2,3,4,-4,-4,4)$  &$1$\\
\hline
$E_6(a_1)$ &$(0,1,1,2,3,-3,-3,3)$  &$1$\\
\hline
$D_5$ &$(\frac 12,\frac 12,\frac 32,\frac 32,\frac 52,-\frac 52,-\frac
  52,\frac 52)$  &$T_1$\\
\hline
$E_6(a_3)$ &$(0,0,1,1,2,-2,-2,2)$  &$1$\\
\hline
$D_5(a_1)$ &$(\frac 14,\frac 34,\frac 34,\frac 54,-\frac 74,-\frac
  74,\frac 74)$
  &$T_1$\\
\hline
$A_5$ &$(-\frac {11}4,-\frac 74,-\frac 34,\frac 14,\frac 54,-\frac
  54,-\frac 54,\frac 54)+\nu(\frac 12,\frac 12,\frac 12,\frac 12,\frac 12,-\frac 12,-\frac
  12,\frac 12)$  &$A_1$\\

\hline
$A_4+A_1$ &$(0,\frac 12,\frac 12,1,\frac 32,-\frac 32,-\frac 32,\frac
  32)$  &$T_1$\\
\hline
$D_4$ &$(0,1,2,3,\nu,-\nu,-\nu,\nu)$  &$A_2$\\
\hline
$A_4$ &$(-2,-1,0,1,2,0,0,0)+\nu(\frac 12,\frac 12,\frac 12,\frac 12,\frac 12,-\frac
  12,-\frac 12,\frac 12)$  &$A_1T_1$\\
\hline
$D_4(a_1)$ &$(0,0,1,1,1,-1,-1,1)$  &$T_2$\\
\hline
$A_3+A_1$ &$(-\frac 54,-\frac 14,\frac 34,-\frac 54,-\frac 14,-\frac
  34,-\frac 34,\frac 34)+\nu(\frac 12,\frac 12,\frac 12,\frac 12,\frac 12,-\frac 12,-\frac
  12,\frac 12)$  &$A_1T_1$\\
\hline
$2A_2+A_1$ &$(0,1,-\frac 32,-\frac 12,\frac 12,-\frac 12,-\frac
  12,\frac 12)+\nu(0,0,1,1,1,-1,-1,1)$  &$A_1$\\
\hline
$A_3$ &$(-\frac 32,-\frac 12,\frac 12,\frac
  32,0,0,0,0)+(\frac {\nu_1}2,\frac {\nu_1}2,\frac {\nu_1}2,\frac {\nu_1}2,\frac
  {\nu_2}2,-\frac {\nu_2}2,-\frac {\nu_2}2,\frac {\nu_2}2)$  &$B_2T_1$\\

\hline
$A_2+2A_1$ &$(\frac 54,-\frac 14,\frac 34,-\frac 34,\frac 14,-\frac
  14,-\frac 14,\frac 14)+\nu(-\frac 12,\frac 12,\frac 12,\frac 32,\frac 32,-\frac 32,-\frac 32,\frac 32)$   &$A_1T_1$\\

\hline
$2A_2$ &$(-\frac 12,\frac 12,-\frac 32,-\frac 12,\frac 12,-\frac
  12,-\frac 12,\frac 12)+$  &$G_2$\\
&$(\frac {\nu_2}2,\frac {\nu_2}2,\frac{2\nu_1+\nu_2}2,\frac
  {2\nu_1+\nu_2}2,\frac{2\nu_1+\nu_2}2,-\frac
  {2\nu_1+\nu_2}2,-\frac{2\nu_1+\nu_2}2,\frac {2\nu_1+\nu_2}2)$ &\\
\hline

$A_2+A_1$ &$(-\frac 12,\frac 12,-1,0,-\frac 12,-\frac 12,-\frac
  12,\frac 12)+\nu(\frac 12,\frac 12,\frac 12,\frac 12,\frac 12,-\frac 12,-\frac
  12,\frac 12)$  &$A_2T_1$\\

\hline
$A_2$ &$(0,-1,0,1,0,0,0,0)+(\frac {-\nu_1+\nu_2}2,\frac {-\nu_1+\nu_2}2,$  &$2A_2$\\
      &$\frac
  {-\nu_1+\nu_2}2,\frac {-\nu_1+\nu_2}2,\frac {\nu_1+\nu_2}2,\frac
  {-\nu_1+\nu_2}2,\frac {-\nu_1+\nu_2}2,\frac {\nu_1+\nu_2}2)$ &\\

\hline
$3A_1$ &$(0,1,-\frac 12,\frac 12,0,0,0,0)+(0,0,\nu_1,\nu_2,\nu_1,-\nu_1,-\nu_1,\nu_1)$  &$A_2A_1$\\
     \hline
$2A_1$ &$(-\frac 12,\frac 12,-\frac 12,\frac 12,0,0,0,0)+$
   &$B_3T_1$\\
&$(\frac {-\nu_1+\nu_2}2,\frac {-\nu_1+\nu_2}2,\frac
  {\nu_1+\nu_2}2,\frac {\nu_1+\nu_2}2,\frac {\nu_1}2,-\frac
  {\nu_1}2,-\frac {\nu_1}2,\frac {\nu_1}2)$ &\\

\hline
$A_1$ &$(\frac 12,\frac 12,0,0,0,0,0,0)+(\frac {-\nu_1+\nu_2}2,\frac {\nu_1-\nu_2}2,$
  &$A_5$\\*
      &$\frac
  {-\nu_1+\nu_2}2+\nu_3,\frac {\nu_1-\nu_2}2+\nu_3,\frac {\nu_1+\nu_2}2,-\frac
  {\nu_1+\nu_2}2,-\frac {\nu_1+\nu_2}2,\frac {\nu_1+\nu_2}2)$ &\\

\end{longtable}
\end{center}
\end{small}

\subsubsection{$E_7$}

\begin{small}
\begin{center}
\begin{longtable}{|c|c|c|}
\caption{Table of parameters $(\check\CO,\nu)$ for $E_7$}\label{table:E7}\\
\hline
\multicolumn{1}{|c|}{$\mathbf{\check\CO}$} &
\multicolumn{1}{c|}{$\mathbf\chi$}
& \multicolumn{1}{c|}{$\mathbf{\fz(\check\CO)}$} \\ \hline
\endfirsthead

\multicolumn{3}{c}%
{{ \tablename\ \thetable{} -- continued from previous page}} \\
\hline
\multicolumn{1}{|c|}{$\mathbf{\check\CO}$} &
\multicolumn{1}{c|}{$\mathbf\chi$}
& \multicolumn{1}{c|}{$\mathbf{\fz(\check\CO)}$}
 \\ \hline
\endhead

%\hline \multicolumn{4}{|r|}{{Continued on next page}} \\ \hline
%\endfoot

\hline \hline
\endlastfoot

%$\mathbf\check\CO$ &$\mathbf\chi$ &$\mathbf\nu$ &$\mathbf{\fz(\check\CO)}$\\
%\hline
%\hline

$E_7$ &$(0,1,2,3,4,5,-\frac {17}2,\frac {17}2)$   &$1$\\
\hline
$E_7(a_1)$  &$(0,1,1,2,3,4,-\frac{13}2,\frac {13}2)$   &$1$\\
\hline
$E_7(a_2)$  &$(0,1,1,2,2,3,-\frac {11}2,\frac {11}2)$   &$1$\\
\hline
$E_7(a_3)$  &$(0,0,1,1,2,3,-\frac 92,\frac 92)$   &$1$\\
\hline
$E_6$       &$(0,1,2,3,4,-4,-4,4)+\nu(0,0,0,0,0,1,-\frac 12,\frac 12)$
  &$A_1^\ell$\\
\hline
$D_6$  &$(0,1,2,3,4,5,0,0)+\nu(0,0,0,0,0,0,-1,1)$     &$A_1$\\
\hline
$E_6(a_1)$ &$(0,1,1,2,3,-3,-3,3)+\nu(0,0,0,0,0,1,-\frac 12,\frac 12)$
 &$T_1$\\
\hline
$E_7(a_4)$ &$(0,0,1,1,1,2,-\frac 72,\frac 72)$  &$1$\\
\hline
$D_6(a_1)$ &$(0,1,1,2,3,4,0,0)+\nu(0,0,0,0,0,0,-1,1)$  &$A_1$\\
\hline
$A_6$ &$(-\frac 72,-\frac 52,-\frac 32,-\frac 12,\frac 12,\frac
32,-\frac 32,\frac 32)+\nu(\frac 12,\frac 12,\frac 12,\frac 12,\frac
12,\frac 12,-1,1)$  &$A_1^\ell$\\
\hline
$D_5+A_1$ &$(0,1,2,3,-\frac 52,-\frac 32,-2,2)+\nu(0,0,0,0,1,1,-1,1)$
 &$A_1$\\
\hline
$E_7(a_5)$ &$(0,0,1,1,1,2,-\frac 52,\frac 52)$   &$1$\\
\hline
$D_6(a_2)$ &$(0,1,1,2,2,3,0,0)+\nu(0,0,0,0,0,0,-1,1)$  &$A_1$\\
\hline
$A_5+A_1$ &$(\frac {11}4,-\frac 74,-\frac 34,\frac 14,\frac 54,\frac
94,-\frac 14,\frac 14)+\nu(-\frac 12,\frac 12,\frac 12,\frac 12,\frac
12,\frac 12,-\frac 32,\frac 32)$  &$A_1$\\
\hline
$D_5$ &$(0,1,2,3,-2,-2,-2,2)+\nu_1(0,0,0,0,1,1,-1,1)$
&$2A_1$\\
      &$+\nu_2(0,0,0,0,-1,1,0,0)$ &\\
\hline
$E_6(a_3)$ &$(0,0,1,1,2,-2,-2,2)+\nu(0,0,0,0,0,1,-\frac 12,\frac 12)$
 &$A_1^\ell$\\
\hline
$D_5(a_1)A_1$ &$(0,1,1,2,-2,-1,-\frac 32,\frac
32)+\nu(0,0,0,0,1,1,-1,1)$  &$A_1^\ell$\\
\hline
$(A_5)'$ &$(-\frac 52,-\frac 32,-\frac 12,\frac 12,\frac 32,\frac
52,0,0)+\nu_1(0,0,0,0,0,0,-1,1)$  &$A_1A_1^\ell$\\
&$+\nu_2(\frac 12,\frac 12,\frac 12,\frac 12,\frac 12,\frac 12,0,0)$
&\\
\hline
$A_4+A_2$ &$(0,1,2,-2,-1,0,-1,1)+\nu(0,0,0,1,1,1,-\frac 32,\frac 32)$
 &$A_1^\ell$\\
\hline
$(A_5)''$ &$(\frac 52,-\frac 32,-\frac 12,\frac 12,\frac 32,\frac
52,0,0)+\nu_2(0,0,0,0,0,0,-1,1)$  &$G_2$\\
&$+\nu_1(-\frac 12,\frac 12,\frac 12,\frac 12,\frac 12,\frac 12,-\frac
32,\frac 32)$ &\\
\hline
$D_5(a_1)$ &$(0,1,1,2,3,0,0,0)+\nu_1(0,0,0,0,0,0,-1,1)$
 &$A_1T_1$\\
&$+\nu_2(0,0,0,0,0,1,-\frac 12,\frac 12)$ &\\
\hline
$A_4+A_1$ &$(\frac 94,-\frac 54,-\frac 14,\frac 34,\frac 74,-\frac
14,-\frac 12,\frac 14)+\nu_1(\frac 12,\frac 12,\frac 12,\frac 12,\frac
12,\frac 12,-1,1)$  &$T_2$\\
&$+\nu_2(0,0,0,0,0,1,-\frac 12,\frac 12)$  &\\
\hline
$D_4+A_1$ &$(0,1,2,3,-\frac 12,\frac 12,0,0)+\nu_1(0,0,0,0,-\frac
12,-\frac 12,-\frac 12,\frac 12)$  &$B_2$\\
&$+\nu_2(0,0,0,0,\frac 12,\frac 12,-\frac 12,\frac 12)$  &\\
\hline
$A_3A_2A_1$ &$(0,1,-2,-1,0,1,-\frac 12,\frac
12)+\nu(0,0,1,1,1,1,-2,2)$  &$A_1$\\
\hline
$A_4$ &$(0,-2,-1,0,1,2,0,0)+\nu_1(0,0,0,0,0,0,-1,1)$
 &$A_2T_1$ \\
&$+\nu_2(-\frac 12,\frac 12,\frac 12,\frac 12,\frac 12,\frac 12,-\frac
32,\frac 32)+\nu_3(\frac 12,\frac 12,\frac 12,\frac 12,\frac 12,\frac
12,-1,1)$ &\\
\hline
$A_3+A_2$ &$(0,1,2,-1,0,1,0,0)+\nu_1(0,0,0,0,0,0,-1,1)$
 &$A_1T_1$\\
&$+\nu_2(0,0,0,1,1,1,0,0)$ &\\
\hline
$D_4$ &$(0,1,2,3,\nu_2-\nu_1,\nu_2+\nu_1,-\nu_3,\nu_3)$
 &$C_3$\\
\hline
$D_4(a_1)A_1$ &$(0,1,1,2,-\frac 12,\frac
12,0,0)+(0,0,0,0,\nu_2,\nu_2,-\nu_1,\nu_1)$  &$2A_1$\\
\hline
$A_3+2A_1$ &$(0,1,-\frac 32,-\frac 12,\frac 12,\frac
32,0,0)+(0,0,\nu_2,\nu_2,\nu_2,\nu_2,-\nu_1,\nu_1)$
&$2A_1$\\
\hline
$D_4(a_1)$ &$(0,1,1,2,\nu_2-\nu_3,\nu_2+\nu_3,-\nu_1,\nu_1)$
 &$3A_1$\\
\hline
$(A_3+A_1)'$ &$(0,1,2,0,-\frac 12,\frac
12,0,0)+(0,0,0,2\nu_2,\nu_3,\nu_3,-\nu_1,\nu_1)$
 &$3A_1$\\
\hline
$2A_2+A_1$ &$(\frac 54,-\frac 14,\frac 34,-\frac 54,-\frac 14,\frac
34,-\frac 14,\frac 14)+\nu_1(1,-1,-1,1,1,1,0,0)$
&$2A_1$\\
&$+\nu_2(-\frac 12,\frac 12,\frac 12,\frac 12,\frac 12,\frac 12,-\frac
32,\frac 32)$ &\\
\hline
$(A_3+A_1)''$ &$(\frac 32,-\frac 12,\frac 12,\frac 32,-\frac 12,\frac
12,0,0)+$  &$B_3$\\
&$(-\frac {\nu_1}2,\frac {\nu_1}2,\frac {\nu_1}2,\frac {\nu_1}2,\frac
{\nu_3-\nu_2}2,\frac {\nu_3-\nu_2}2,-\frac {\nu_3+\nu_2}2,\frac
{\nu_3+\nu_2}2)$  &\\
\hline
$A_2+3A_1$ &$(0,1,-1,0,-1,0,-\frac 12,\frac
12)+\nu_1(0,0,1,1,1,1,-2,2)$  &$G_2$\\
&$+\nu_2(0,0,0,0,1,1,-1,1)$  &\\
\hline
%\end{longtable}
%\end{table}

%\begin{table}\label{table:E7,2}
%\caption{Table of nilpotent orbits for $E_7$, continued}
%\noindent\begin{longtable}{|c|c|c|c|}
%\hline
%$\mathbf\check\CO$ &$\mathbf\chi$ &$\mathbf\nu$ &$\mathbf{\fz(\check\CO)}$\\
%\hline
%\hline

$2A_2$ &$(-\frac 12,\frac 12,-\frac 32,-\frac 12,\frac 12,-\frac
12,-\frac 12,\frac 12)+\nu_1(0,0,1,1,1,-1,-1,1)$
 &$G_2A_1$\\
&$+\nu_2(\frac 12,\frac 12,\frac 12,\frac 12,\frac 12,-\frac 12,-\frac
12,\frac 12)+\nu_3(0,0,0,0,0,1,-\frac 12,\frac 12)$ & \\
\hline
$A_3$ &$(0,1,2,\nu_1,\nu_2,\nu_3,\nu_4)$
&$B_3A_1$\\
\hline
$*A_2+2A_1$
&$(0,1,-1,0,1,0,0,0)+(0,0,\nu_2,\nu_2,\nu_2,\nu_3,-\nu_1,\nu_1)$
 &$A_12A_1^\ell$\\
\hline
$A_2+A_1$ &$(1,0,1,0,-\frac 12,\frac 12,0,0)+(0,0,0,0,\nu_2,\nu_2,-\nu_1,\nu_1)$
 &$A_3T_1$\\
&$+\nu_3(0,0,0,1,1,1,-\frac 32,\frac 32)+\nu_4(-\frac 12,\frac
12,\frac 12,\frac 12,\frac 12,\frac 12,-\frac 32,\frac 32)$ &\\
\hline
$4A_1$ &$(0,1,-\frac 12,\frac 12,-\frac 12,\frac
12,0,0)+(0,0,\nu_3,\nu_3,\nu_2,\nu_2,-\nu_1,\nu_1)$
 &$C_3$ \\
\hline
$A_2$
&$(1,0,1,0,0,0,0,0)+(0,0,0,0,\nu_2-\nu_3,\nu_2+\nu_3,-\nu_1,\nu_1)$
 &$A_5$\\
&$+\nu_4(-\frac 12,\frac 12,\frac 12,\frac 12,\frac 12,\frac 12,-\frac
32,\frac 32)+\nu_5(0,0,0,1,1,1,-\frac 32,\frac 32)$ &\\
\hline
$(3A_1)'$ &$(-\frac 12,\frac 12,-\frac 12,\frac 12,-\frac 12,\frac
12,0,0)+(\nu_1,\nu_1,\nu_2,\nu_2,\nu_3,\nu_3,-\nu_4,\nu_4)$
 &$C_3A_1$\\
\hline
$(3A_1)''$ &$(\frac 12,\frac 12,-\frac 12,\frac 12,-\frac 12,\frac
12,0,0)+(-\nu_4,\nu_4,\nu_3,\nu_3,\nu_2,\nu_2,-\nu_1,\nu_1)$
 &$F_4$\\
\hline
$2A_1$ &$(0,1,\nu_1,\nu_2,\nu_3,\nu_4,-\nu_5,\nu_5)$
 &$B_4A_1$\\
\hline
$A_1$ &$(\frac {\nu_1+\nu_2+\nu_3-\nu_4}2,\frac
    {\nu_1+\nu_2-\nu_3+\nu_4}2,\frac {\nu_1-\nu_2+\nu_3+\nu_4}2,\frac
    {-\nu_1+\nu_2+\nu_3+\nu_4}2,$  &$D_6$\\
&$-\frac 12+\frac {-\nu_5+\nu_6}2,\frac 12+\frac
    {-\nu_5+\nu_6}2,-\frac {\nu_5+\nu_6}2,\frac {\nu_5+\nu_6}2)$ &\\
\hline

\end{longtable}
\end{center}
\end{small}

\noindent{\bf $\mathbf{E_7}$ exception:}

\noindent $\mathbf{A_2+2A_1}$. Three regions: $\{0\le\nu_1<\frac
12,0\le\nu_2<1,0\le\nu_3<1, \nu_1+\frac {3\nu_2}2+\frac{\nu_3}2<\frac
32\}$, $\{0\le\nu_1<\frac 12,0\le\nu_2<1,0\le\nu_3<1, -\nu_1+\frac
{3\nu_2}2+\frac{\nu_3}2<\frac 32,\nu_1+\frac
{3\nu_2}2-\frac{\nu_3}2>\frac 32 \}$, and  $\{0\le\nu_1<\frac
12,0\le\nu_2<1,0\le\nu_3<1, \frac{3\nu_2}2+\frac{\nu_3}2>\frac 32,\nu_1+\frac
{3\nu_2}2-\frac{\nu_3}2<\frac 32 \}$.

\subsubsection{$E_8$}

\begin{small}
\begin{center}
\begin{longtable}{|c|c|c|}
\caption{Table of parameters $(\check\CO,\nu)$ for $E_8$}\label{table:E8}\\
\hline
\multicolumn{1}{|c|}{$\mathbf{\check\CO}$} &
\multicolumn{1}{c|}{$\mathbf\chi$}
& \multicolumn{1}{c|}{$\mathbf{\fz(\check\CO)}$} \\ \hline
\endfirsthead

\multicolumn{3}{c}%
{{  \tablename\ \thetable{} -- continued from previous page}} \\
\hline
\multicolumn{1}{|c|}{$\mathbf{\check\CO}$} &
\multicolumn{1}{c|}{$\mathbf\chi$}
& \multicolumn{1}{c|}{$\mathbf{\fz(\check\CO)}$}
 \\ \hline
\endhead

%\hline \multicolumn{4}{|r|}{{Continued on next page}} \\ \hline
%\endfoot

\hline \hline
\endlastfoot

$E_8$ &$(0,1,2,3,4,5,6,23)$   &$1$\\
\hline
$E_8(a_1)$ &$(0,1,1,2,3,4,5,18)$ &$1$\\
\hline
$E_8(a_2)$ &$(0,1,1,2,2,3,4,15)$ &$1$\\
\hline
$E_8(a_3)$ &$(0,0,1,1,2,3,4,13)$ &$1$\\
\hline
$E_8(a_4)$ &$(0,0,1,1,2,2,3,11)$ &$1$\\
\hline
$E_7$   &$(0,1,2,3,4,5,-\frac{17}2,\frac {17}2)+\nu(0,0,0,0,0,0,1,1)$
&$A_1$\\
\hline
$E_8(b_4)$ &$(0,0,1,1,1,2,3,10)$ &$1$\\
\hline
$E_8(a_5)$ &$(0,0,1,1,1,2,2,9)$ &$1$\\
\hline
$E_7(a_1)$
&$(0,1,1,2,3,4,-\frac{13}2,\frac{13}2)+\nu(0,0,0,0,0,0,1,1)$ &$A_1$\\
\hline
$E_8(b_5)$ &$(0,0,1,1,1,2,3,8)$ &$1$\\
\hline
$D_7$ &$(0,1,2,3,4,5,6,0)+\nu(0,0,0,0,0,0,0,2)$ &$A_1$\\
\hline
$E_8(a_6)$ &$(0,0,1,1,1,2,2,7)$ &$1$\\
\hline
$E_7(a_2)$ &$(0,1,1,2,2,3,-\frac{11}2,\frac
{11}2)+\nu(0,0,0,0,0,0,1,1)$ &$A_1$\\
\hline
$E_6+A_1$ &$(0,1,2,3,4,-\frac 92,-\frac 72,4)+\nu(0,0,0,0,0,0,1,1,2)$
&$A_1$\\
\hline
$D_7(a_1)$ &$(0,1,1,2,3,4,5,0)+\nu(0,0,0,0,0,0,0,2)$ &$T_1$\\
\hline
$E_8(b_6)$ &$(0,0,1,1,1,1,2,6)$ &$1$\\
\hline
$E_7(a_3)$ &$(0,0,1,1,2,3,-\frac 92,\frac 92)+\nu(0,0,0,0,0,0,1,1)$
&$A_1$\\
\hline
$E_6(a_1)A_1$  &$(0,1,1,2,3,-\frac 72,-\frac
52,3)+\nu(0,0,0,0,0,1,1,2)$  &$T_1$\\
\hline
$A_7$ &$(-\frac{17}4,-\frac {13}4,-\frac 94,-\frac 54,-\frac 14,\frac
34,\frac 74,\frac 74)+\nu(\frac 12,\frac 12,\frac 12,\frac 12,\frac
12,\frac 12,\frac 12,\frac 52)$ &$T_1$\\
\hline
$E_6$ &$(0,1,2,3,4,-4,-4,4)+\nu_1(0,0,0,0,0,1,1,2)$ &$G_2$\\
      &$+\nu_2(0,0,0,0,0,0,1,1)$ &\\
\hline
$D_6$ &$(0,1,2,3,4,5,\nu_1,\nu_2)$ &$B_2$\\
\hline
$D_5+A_2$ &$(0,1,2,3,-3,-2,-1,2)+\nu(0,0,0,0,1,1,1,3)$ &$T_1$\\
\hline
$E_6(a_1)$ &$(0,1,1,2,3,-3,-3,3)+\nu_2(0,0,0,0,0,1,1,2)$ &$A_2$\\
           &$+\nu_1(0,0,0,0,0,0,1,1)$ &\\
\hline
$E_7(a_4)$ &$(0,0,1,1,1,2,-\frac 72,\frac 72)+\nu(0,0,0,0,0,0,1,1)$
&$A_1$\\
\hline
$A_6+A_1$ &$(\frac {13}4,-\frac 94,-\frac 54,-\frac 14,\frac34,\frac
74,\frac {11}4,\frac 14)+\nu(-\frac 12,\frac 12,\frac 12,\frac
12,\frac 12,\frac 12,\frac 12,\frac 72)$ &$A_1$\\
\hline
$D_6(a_1)$ &$(0,1,1,2,3,4,0,0)+\nu_1(0,0,0,0,0,0,-1,1)$
&$2A_1$\\
           &$+(0,0,0,0,0,0,1,1)$ &\\
\hline
$A_6$ &$(-3,-2,-1,0,1,2,3,0)+\nu_2(\frac 12,\frac 12,\frac 12,\frac 12,\frac 12,\frac
12,\frac 12,\frac 12)$ &$2A_1$\\
&$+\nu_1(-\frac 12,-\frac 12,-\frac 12,-\frac
12,-\frac 12,-\frac 12,\frac 72)$ &\\
\hline
$E_8(a_7)$ &$(0,0,0,1,1,1,1,4)$ &$1$\\
\hline
$D_5+A_1$ &$(0,1,2,3,4,-\frac 12,\frac 12,0)+\nu_1(0,0,0,0,0,0,0,2)$
&$2A_1$\\
         &$+\nu_2(0,0,0,0,0,1,1,0)$ &\\
\hline
$E_7(a_5)$ &$(0,0,1,1,1,2,-\frac 52,\frac 52)+\nu(0,0,0,0,0,0,1,1)$
&$A_1$\\
\hline
$E_6(a_3)A_1$ &$(0,0,1,1,2,-\frac52,-\frac 32,2)+\nu(0,0,0,0,0,1,1,2)$
&$A_1$\\
\hline
$D_6(a_2)$ &$(0,1,1,2,2,3,-\nu_1+\nu_2,\nu_1+\nu_2)$ &$2A_1$\\
\hline
$D_5(a_1)A_2$ &$(0,1,1,2,-\frac 52,-\frac 32,-\frac 12,\frac
32)+\nu(0,0,0,0,1,1,1,3)$ &$A_1$\\
\hline
$A_5+A_1$ &$(\frac 14,-\frac{11}4,-\frac 74,-\frac 34,\frac 14,\frac
54,\frac 94,\frac 14)+\nu_2(-1,0,0,0,0,0,0,1)$ &$2A_1$\\
        &$+\nu_1(\frac 32,\frac 12,\frac 12,\frac 12,\frac 12,\frac
        12,\frac 12,\frac 32)$ &\\
\hline
$A_4+A_3$ &$(0,1,2,-\frac 52,-\frac 32,-\frac 12,\frac
12,1)+\nu(0,0,0,1,1,1,1,4)$ &$A_1$\\
\hline
$D_5$ &$(0,1,2,3,4,\nu_1,\nu_2,\nu_3)$ &$B_3$\\
\hline
$E_6(a_3)$ &$(0,0,1,1,2,-2,-2,2)+\nu_1(0,0,0,0,0,1,1,2)$
&$G_2$\\
&$+\nu_2(0,0,0,0,0,0,1,1)$ &\\
\hline
$D_4+A_2$ &$(0,1,2,3,-1,0,1,0)+\nu_2(0,0,0,0,1,1,1,3)$ &$A_2$\\
         &$+\nu_1(0,0,0,0,0,0,0,2)$ &\\
\hline
$*A_4A_2A_1$ &$(0,1,-\frac 52,-\frac 32,-\frac 12,\frac 12,\frac
32,\frac 12)+\nu(0,0,1,1,1,1,1,5)$ &$A_1$\\
\hline
$*D_5(a_1)A_1$ &$(0,1,1,2,3,-\frac 12+\nu_2,\frac 12+\nu_2,2\nu_1)$
&$A_1^\ell A_1$\\
\hline
$A_5$ &$(\frac 52,-\frac 32,-\frac 12,\frac 12,\frac 32,\frac
52,0,0)+\nu_1(-\frac 12,\frac 12,\frac 12,\frac 12,\frac 12,\frac
12,-\frac 32,\frac 32)$ &$G_2A_1$\\
&$+\nu_2(0,0,0,0,0,0,-1,1)+\nu_3(0,0,0,0,0,0,1,1)$ &\\
\hline
$*A_4+A_2$ &$(-\frac 12,\frac 12,-\frac 52,-\frac 32,-\frac 12,\frac
12,\frac 32,\frac 12)+\nu_2(1,1,0,0,0,0,0,0)$ &$2A_1$\\
&$+\nu_1(0,0,1,1,1,1,1,5)$ &\\
\hline
$A_4+2A_1$ &$(0,1,-2,-1,0,1,2,0)+\nu_1(0,0,0,0,0,0,0,2)$
&$A_1T_1$\\
&$\nu_2(0,0,1,1,1,1,1,0)$ &\\
\hline
$D_5(a_1)$ &$(0,1,1,2,3,\nu_3,\nu_2,\nu_1)$ &$A_3$\\
          &                                &\\
\hline
$2A_3$ &$(0,1,2,-\frac 32,-\frac 12,\frac 12,\frac
32,0)+\nu_2(0,0,0,\frac 12,\frac 12,\frac 12,\frac 12,1)$
&$B_2$\\
\hline
$A_4+A_1$ &$(0,1,2,-\frac 32,-\frac
12,-1,-1,1)+\nu_2(0,0,0,0,0,1,1,2)$ &$A_2T_1$\\
& $+\nu_1(0,0,0,0,0,0,1,1)+\nu_3(0,0,0,1,1,1,1,4)$
 &\\
\hline
$D_4(a_1)A_2$ &$(0,1,1,2,-1,0,1,0)+\nu_1(0,0,0,0,1,1,1,3)$
&$A_2$\\
             &$+\nu_2(0,0,0,0,0,0,0,2)$ &\\
\hline
$D_4+A_1$ &$(0,1,2,3,-\frac 12,\frac 12,0,0)+$ &$C_3$\\
&$(0,0,0,0,\nu_1,\nu_1,-\nu_2+\nu_3,\nu_2+\nu_3)$ &\\
\hline
$A_3A_2A_1$ &$(0,1,-2,-1,0,1,-\frac 12,\frac
12)+\nu_1(0,0,1,1,1,1,-2,2)$ &$2A_1$\\
&$+\nu_2(0,0,0,0,0,0,1,1)$ &\\
\hline
$A_4$ &$(0,-2,-1,0,1,2,0,0)+$ &$A_4$\\*
&$(\nu_4,-\nu_1+\nu_2,\nu_3,\nu_3,\nu_3,\nu_3,\nu_3,\nu_1+\nu_2)$ &\\
\hline
$A_3+A_2$ &$(0,1,2,-1,0,1,0,0)+(0,0,0,\nu_3,\nu_3,\nu_3,\nu_1,\nu_2)$
&$B_2T_1$\\
\hline
$D_4(a_1)A_1$ &$(0,1,1,2,-\frac 12,\frac
12,0,0)+(0,0,0,0,\nu_1,\nu_1,-\nu_2+\nu_3,\nu_2+\nu_3)$
&$3A_1$\\
\hline
$A_3+2A_1$ &$(0,1,-\frac 32,-\frac 12,\frac 12,\frac
32,0,0)+(0,0,\nu_1,\nu_1,\nu_1,\nu_1,\nu_2,\nu_3)$
&$A_1B_2$\\
\hline
$2A_2+2A_1$ &$(0,1,-\frac 32,-\frac 12,\frac 12,-1,0,\frac
12)+\nu_1(0,0,-\frac 12,-\frac 12,-\frac 12,1,1,\frac 12)$ &$B_2$\\
&$+\nu_2(0,0,\frac 12,\frac 12,\frac 12,0,0,\frac 32)$ &\\
\hline
$D_4$
&$(0,1,2,3,\nu_3-\nu_4,\nu_3+\nu_4,\nu_1-\nu_2,\nu_1+\nu_2)$
&$F_4$\\
\hline
$D_4(a_1)$ &$(0,1,1,2,\nu_4,\nu_3,\nu_2,\nu_1)$  &$D_4$\\
\hline
$A_3+A_1$ &$(0,1,2,-\frac 12,\frac
12,0,0,0)+(0,0,0,\nu_1,\nu_1,\nu_2,\nu_3,\nu_4)$
&$A_1B_3$\\
\hline
$2A_2+A_1$ &$(0,1,-\frac 32,-\frac 12,\frac 12,-\frac 12,-\frac
12,\frac 12)+\nu_1(0,0,1,1,1,-1,-1,1)$ &$A_1G_2$\\
&$+\nu_2(0,0,0,0,0,1,1,2)+\nu_3(0,0,0,0,0,0,1,1)$ &\\
\hline
$2A_2$ &$(-\frac 12,\frac 12,-\frac 32,-\frac 12,\frac 12,-\frac
12,-\frac 12,\frac 12)+\nu_1(0,0,1,1,1,-1,-1,1)$
&$2G_2$\\
&$+\nu_2(\frac 12,\frac 12,\frac 12,\frac 12,\frac 12,-\frac 12,-\frac
12,\frac 12)+\nu_3(0,0,0,0,0,1,1,2)$ &\\
&$+\nu_4(0,0,0,0,0,0,1,1)$ &\\
\hline
$*A_2+3A_1$ &$(0,1,-1,0,-1,0,-\frac 12,\frac
12)+\nu_1(0,0,1,1,1,1,-2,2)$ &$G_2A_1$\\
&$+\nu_2(0,0,0,0,1,1,-1,1)+\nu_3(0,0,0,0,0,0,1,1)$ &\\
\hline
$A_3$ &$(0,1,2,\nu_1,\nu_2,\nu_3,\nu_4,\nu_5)$ &$B_5$\\
\hline
$*A_2+2A_1$
&$(0,1,-1,0,1,0,0,0)+(0,0,\nu_1,\nu_1,\nu_1,\nu_2,\nu_3,\nu_4)$
 &$A_1B_3$\\
\hline
$A_2+A_1$ &$(1,0,1,0,-\frac 12,\frac 12,0,0)+$ & $A_5$\\
&$(-\nu_5,\nu_5,\nu_5,\nu_4,\nu_3,\nu_3,-\nu_2+\nu_1,\nu_2+\nu_1)$
&\\
\hline
$*4A_1$ &$(0,1,-\frac 12,\frac 12,-\frac 12,\frac 12,0,0)+$
&$C_4$\\
&$(0,0,\nu_1,\nu_1,\nu_2,\nu_2,-\nu_3+\nu_4,\nu_3+\nu_4)$&\\
\hline
$A_2$ &$(\frac
{\nu_1-\nu_2-\nu_3+\nu_4}2,\frac{-\nu_1+\nu_2-\nu_3+\nu_4}2,\frac{-\nu_1-\nu_2+\nu_3+\nu_4}2,$
&$E_6$\\
&$\frac{\nu_1+\nu_2+\nu_3+\nu_4}2,-1+\frac
{\nu_5-\nu_6}2,\frac{\nu_5-\nu_6}2,1+\frac{\nu_5-\nu_6}2,\frac{\nu_5+3\nu_6}2)$
&\\
\hline
$3A_1$ &$(\frac 12,\frac 12,-\frac 12,\frac 12,-\frac 12,\frac
12,0,0)+(-\nu_4,\nu_4,\nu_3,\nu_3,\nu_2,\nu_2,-\nu_1,\nu_1)$
&$F_4A_1$\\
&$+\nu_5(0,0,0,0,0,0,1,1)$ &\\
\hline
$2A_1$ &$(0,1,\nu_1,\nu_2,\nu_3,\nu_4,\nu_5,\nu_6)$
&$B_6$\\
\hline
$A_1$
&$(\frac{\nu_1+\nu_2+\nu_3-\nu_4}2,\frac{\nu_1+\nu_2-\nu_3+\nu_4}2,\frac
{\nu_1-\nu_2+\nu_3+\nu_4}2,\frac {-\nu_1+\nu_2+\nu_3+\nu_4}2,$ &$E_7$\\
&$\frac
{-\nu_5-\nu_6+2\nu_7}2,-\frac 12+\frac {-\nu_5+\nu_6}2,\frac
12+\frac{-\nu_5+\nu_6}2, \frac {\nu_5+\nu_6+2\nu_7}2)$ &\\

\end{longtable}
\end{center}
\end{small}

\noindent{\bf $\mathbf {E_8}$ exceptions:}

\noindent$\mathbf{A_4+A_2+A_1}.$ $\{0\le\nu<\frac 3{10}\}.$

\noindent$\mathbf{D_5(a_1)+A_1}.$ Two regions: $\{0\le\nu_2<\frac 12,2\nu_1+\nu_2<\frac
32\},$ and $\{0\le\nu_1<1,2\nu_1-\nu_2>\frac 32\}.$

\noindent$\mathbf{A_4+A_2}.$ Two regions: $\{0\le\nu_2<\frac
12,5\nu_1+\nu_2<2\},$ and $\{0\le\nu_1<\frac 12,5\nu_1-\nu_2>2\}.$

\noindent$\mathbf{A_2+3 A_1}$. Four regions: $\{3\nu_1+2\nu_2<1,
0\le\nu_3<\frac 12\}$, $\{2\nu_1+\nu_2<1<3\nu_1+\nu_2,0\le\nu_3<\frac
12,3\nu_1+2\nu_2+\nu_3<\frac 32\}$,
$\{2\nu_1+\nu_2<1<3\nu_1+\nu_2,0\le\nu_3<\frac
12,3\nu_1+\nu_2+\nu_3<\frac 32<3\nu_1+2\nu_2-\nu_3\}$, and
$\{2\nu_1+\nu_2<1<3\nu_1+\nu_2,0\le\nu_3<\frac 12,
3\nu_1+2\nu_2-\nu_3<\frac 32<3\nu_1+\nu_2+\nu_3\}.$

\noindent$\mathbf{A_2+2A_1}.$ Seven regions:
$\{0\le\nu_1<1,\nu_3+\nu_4<1, 3\nu_1+\nu_2+\nu_3+\nu_4<3\}$,
$\{0\le\nu_1<1,\nu_3+\nu_4<1,3\nu_1+\nu_2-\nu_3+\nu_4<3<3\nu_1-\nu_2+\nu_3+\nu_4\}$,
$\{0\le\nu_1<1,\nu_3+\nu_4<1, 3\nu_1-\nu_2-\nu_3+\nu_4>3\}$,
$\{0\le\nu_1<1,\nu_3+\nu_4<1, 3\nu_1+\nu_2+\nu_3-\nu_4>3 \}$,
$\{0\le\nu_1<1,\nu_2+\nu_4>1,\nu_2+\nu_3<1,\nu_4<1,3\nu_1+\nu_2+\nu_3+\nu_4<3 \}$,
$\{0\le\nu_1<1,\nu_2+\nu_4>1,\nu_2+\nu_3<1,\nu_4<1,
3\nu_1-\nu_2-\nu_3+\nu_4>3 \}$, and
$\{0\le\nu_1<1,\nu_2+\nu_4>1,\nu_2+\nu_3<1,\nu_4<1,
3\nu_1+\nu_2+\nu_3-\nu_4>3\}$.

\noindent$\mathbf{4A_1}$. Two regions:
$\{0\le\nu_1\le\nu-2\le\nu_3\le\nu_4<\frac 12\}$ and
$\{\nu_1+\nu_4<1,\nu_2+\nu_3<1,\nu_2+\nu_4>1, -\nu_1+\nu_3+\nu_4<\frac
32<\nu_1+\nu_3+\nu_4\}.$

%%%%%%%%%%%%%%%%%%%%%%%%%%%%%%%%%%%%%%%%%%%%%%%%%%%%%%%%%%%%%55
\ifx\undefined\bysame
\newcommand{\bysame}{\leavevmode\hbox to3em{\hrulefill}\,}
\fi

\end{document}